\documentclass[12pt]{article}
\usepackage{amssymb,amsmath}

\newtheorem{theorem}{Theorem}
\newtheorem{proposition}{Proposition}
\newtheorem{definition}{Definition}
\newtheorem{corollary}{Corollary}

\newcommand{\Ker}{\operatorname{Ker}}
\renewcommand{\Im}{\operatorname{Im}}
\newcommand{\Hom}{\operatorname{Hom}}
\newcommand{\id}{\operatorname{id}}
\renewcommand{\times}{\operatorname{times}}
\renewcommand{\deg}{\operatorname{deg}}
\renewcommand{\int}{\operatorname{int}}

\begin{document}

\title{

\author{Tornike Kadeishvili\\
A. Razmadze Mathematical Institute of Tbilisi State University}

$A_\infty$-algebra Structure in Cohomology and its Applications}
\date{ }

\maketitle

%\footnotetext[1]{ICTP Map Summer School, 11-29 August 2008, Trieste }

%%%%%%%%%%%%%%%%%%%%%%%%%%%%%%%%%%%%%%%%%%%%%%%%%%%%%%%%%%%%%%%%%%%%%%%%%%%%%%%%%%%%%%%%%%%%%%%%%%%%%%%%%%%%%%%%%
%%%%%%%%%%%%%%%%%%%%%%%%%%%%%%%%%%%%%%%%%%%%%%%%%%%%%%%%%%%%%%%%%%%%%%%%%%%%%%%%%%%%%%%%%%%%%%%%%%%%%%%%%%%%%%%%%
%%%%%%%%%%%%%%%%%%%%%%%%%%%%%%%%%%%%%%%%%%%%%%%%%%%%%%%%%%%%%%%%%%%%%%%%%%%%%%%%%%%%%%%%%%%%%%%%%%%%%%%%%%%%%%%%%
%%%%%%%%%%%%%%%%%%%%%%%%%%%%%%%%%%%%%%%%%%%%%%%%%%%%%%%%%%%%%%%%%%%%%%%%%%%%%%%%%%%%%%%%%%%%%%%%%%%%%%%%%%%%%%%%%

%%%%%%%%%%%%%%%%%%%%%%%%%%%%%%%%%%%%%%%%%%%%%%%%%%%%%%%%%%%%%%%%%%%%%%%%%%%%%%%%%%%%%%%%%%%%%%%%%%%%%%%%%%%%%%%%%
%%%%%%%%%%%%%%%%%%%%%%%%%%%%%%%%%%%%%%%%%%%%%%%%%%%%%%%%%%%%%%%%%%%%%%%%%%%%%%%%%%%%%%%%%%%%%%%%%%%%%%%%%%%%%%%%%
%%%%%%%%%%%%%%%%%%%%%%%%%%%%%%%%%%%%%%%%%%%%%%%%%%%%%%%%%%%%%%%%%%%%%%%%%%%%%%%%%%%%%%%%%%%%%%%%%%%%%%%%%%%%%%%%%
%%%%%%%%%%%%%%%%%%%%%%%%%%%%%%%%%%%%%%%%%%%%%%%%%%%%%%%%%%%%%%%%%%%%%%%%%%%%%%%%%%%%%%%%%%%%%%%%%%%%%%%%%%%%%%%%%

\begin{center}

LECTURE NOTES\\  \vspace{5mm}
Algebra, Topology and Analysis: $C^*$ and $A_\infty$-Algebras \\
\vspace{5mm}

Summer School/Conference\\

\vspace{5mm}

 Ivane Javakhishvili Tbilisi State University Press, 2023

 \end{center}
\vspace{5mm}

%%%%%%%%%%%%%%%%%%%%%%%%%%%%%%%%%%%%%%%%%%%%%%%%%%%%%%%%%%%%%%%%%%%%%%%%%%%%%%%%%%%%%%%%%%%%%%%%%%%%%%%%%%%%%%%%%
%%%%%%%%%%%%%%%%%%%%%%%%%%%%%%%%%%%%%%%%%%%%%%%%%%%%%%%%%%%%%%%%%%%%%%%%%%%%%%%%%%%%%%%%%%%%%%%%%%%%%%%%%%%%%%%%%
%%%%%%%%%%%%%%%%%%%%%%%%%%%%%%%%%%%%%%%%%%%%%%%%%%%%%%%%%%%%%%%%%%%%%%%%%%%%%%%%%%%%%%%%%%%%%%%%%%%%%%%%%%%%%%%%%
%%%%%%%%%%%%%%%%%%%%%%%%%%%%%%%%%%%%%%%%%%%%%%%%%%%%%%%%%%%%%%%%%%%%%%%%%%%%%%%%%%%%%%%%%%%%%%%%%%%%%%%%%%%%%%%%%

%%%%%%%%%%%%%%%%%%%%%%%%%%%%%%%%%%%%%%%%%%%%%%%%%%%%%%%%%%%%%%%%%%%%%%%%%%%%%%%%%%%%%%%%%%%%%%%%%%%%%%%%%%%%%%%%%
%%%%%%%%%%%%%%%%%%%%%%%%%%%%%%%%%%%%%%%%%%%%%%%%%%%%%%%%%%%%%%%%%%%%%%%%%%%%%%%%%%%%%%%%%%%%%%%%%%%%%%%%%%%%%%%%%
%%%%%%%%%%%%%%%%%%%%%%%%%%%%%%%%%%%%%%%%%%%%%%%%%%%%%%%%%%%%%%%%%%%%%%%%%%%%%%%%%%%%%%%%%%%%%%%%%%%%%%%%%%%%%%%%%
%%%%%%%%%%%%%%%%%%%%%%%%%%%%%%%%%%%%%%%%%%%%%%%%%%%%%%%%%%%%%%%%%%%%%%%%%%%%%%%%%%%%%%%%%%%%%%%%%%%%%%%%%%%%%%%%%

%{Introduction}

The main method of algebraic topology is to assign to a
topological space a certain algebraic object (model) and to study
this relatively simple algebraic object instead of complex
geometric one.

Examples of such models are chain and cochain complexes, homology
and homotopy groups, cohomology algebras, etc.

The main problem here is to find models that classify spaces up to
some equivalence relation, such as homeomorphism, homotopy
equivalence, rational homotopy equivalence, etc.

Usually such models are not {\it complete}: the equivalence of
models does not guarantee the equivalence of spaces. They can just
distinguish spaces.

The models which carry richer algebraic structure contain more
information about the space. For example the model "cohomology
algebra" allows to distinguish spaces which cannot be
distinguished by the model "cohomology groups".

Here we are going to present one more additional algebraic structure on cohomology, which was constructed in
\cite{Kad76}, \cite{Kad80}, namely, we show that on cohomology $H^*(X,R)$ there exists Stasheff's $A_\infty$-algebra
structure (minimality theorem). This structure consists of a collection
of operations
$$
\{m_i: H^*(X,R)\otimes ... (i\ \times)...\otimes H^*(X,R)\to
H^*(X,R),\ i=1,2,3,...\}.
$$
In fact this structure extends the usual structure of the cohomology
algebra: the first operation $m_1:H^*(X,R)\to
H^*(X,R)$, the differential, is trivial (minimality) and 
 $m_2:H^*(X,R)\otimes H^*(X,R)\to
H^*(X,R)$ coincides with the cohomology multiplication.

Stasheff's $A_\infty$ algebra is a sort of Strong Homotopy Associative Algebra, the operation $m_3$
is a homotopy which measures the nonassociativity of the product $m_2$. So the existence of an $A_\infty$ algebra structure on a \emph{strictly associative} cohomology algebra $H^*(X,R)$
looks a bit strange, but although the product on $H^*(X,R)$ is associative,
there appears a structure of a (generally nondegenerate) \emph{minimal}
$A_{\infty}$-algebra, which can be considered as an $A_{\infty}$
\emph{deformation} of the classical cohomology $(H^*(X,R),\mu^*)$, \cite{Kad07}.

The cohomology algebra equipped with this additional structure
$$
(H^*(X,R), \{m_i:H^*(X,R)^{\otimes i}\to
H^*(X,R),\ i=2,3,...\},\ m_2=\mu^*\} ),
$$
which we call cohomology $A_{\infty}$-algebra, carries more
information about the space than the cohomology algebra. For
example just the cohomology algebra $H^*(X,R)$ does not determine the
cohomology of the loop space $H^*(\Omega X,R)$, but the cohomology
$A_{\infty}$-algebra $(H^*(X,R),\{m_i\})$ does. Dually, the
Pontriagin ring $H_*(G)$ does not determine the homology $H_*(B_G)$ of
the classifying space, but the homology $A_{\infty}$-algebra
$(H_*(G),\{m_i\})$ does.

These $A_{\infty}$-algebras have several applications in the
cohomology theory of fibre bundles too, see \cite{Kad80}.

But this invariant also is not complete. One cannot expect the
existence of a more or less simple complete algebraic invariant in
general, but for the rational homotopy category there are
various complete homotopy invariants (algebraic models):

(i) The model of Quillen  $L_X$, which is a
differential graded Lie algebra;

(ii) The minimal model of Sullivan $M_X$, which is a
commutative graded differential algebra;

(iii) The filtered model of Halperin and Stasheff
$\Lambda X$, which is a filtered commutative graded differential
algebra.

The rational cohomology algebra $H^*(X,Q)$ is not a complete
invariant even for rational spaces: two spaces might have
isomorphic cohomology algebras, but different rational homotopy
types.
Bellow we also present the main result of \cite{Kad88}: There is the
notion of $C_{\infty}$-algebra which is the commutative version of
 Stasheff's notion of $A_{\infty}$-algebra, and in \cite{Kad88}
we have shown that in the rational case on cohomology $H^*(X,Q)$ arises
a structure of $C_{\infty}$-algebra $(H^*(X,Q),\{m_i\})$. The main
application of this structure is the following: it completely
determines the rational homotopy type, that is, 1-connected spaces
$X$ and $X'$ have the same rational homotopy type if and only if
their cohomology $C_{\infty}$-algebras $(H^*(X,Q),\{m_i\})$ and
$(H^*(X',Q),\{m'_i\})$ are isomorphic.

We present also several applications of this complete rational
homotopy invariant to some problems of rational homotopy theory.

The $C_\infty$-algebra structure in cohomology and the applications of this
structure in rational homotopy theory were already presented in the
hardly available small book \cite{Kad93} (see also the preprint
\cite{Kad88h}).

Applications of cohomology $C_\infty$-algebra in rational homotopy
theory are inspired by the existence of Sullivan's commutative
cochains $A(X)$ in this case. The cohomology $C_\infty$-algebra
$(H^*(X,Q),\{m_i\})$ carries the same amount of information as
$A(X)$ does. Actually these two objects are equivalent in the
category of $C_\infty$-algebras.

We want to remark that for simplicity in these lectures signs are ignored.
Of course they can be reconstructed using the Koszul sign rule.

The organization is as follows.

In Section 1 the notions of chain and cochain complexes
are presented. In Section 2 the differential algebras
and coalgebras are defined. In Section 3 the bar and cobar
constructions are introduced. Twisting cochains and Berikashvili's functor $D$ are presented in Section 4.
In Section 5 the Stasheff's $A_\infty$-algebras are discussed. Hochschild cochains which are used for description of
$A_\infty$-algebras are presented in Section 6.  Section 7 is dedicated to our central topic, the Minimality Theorem.
In the next Section 8 its applications are given. And the last Section 9 is dedicated to
applications of the cohomology $C_\infty$ algebra in rational homotopy theory.

%%%%%%%%%%%%%%%%%%%%%%%%%%%%%%%%%%%%%%%%%%%%%%%%%%%%%%%%%%%%%%%%%%%%%%%%%%%%%%%%%%%%%%%%%%%%%%%%%%%%%%%%%%%%
%%%%%%%%%%%%%%%%%%%%%%%%%%%%%%%%%%%%%%%%%%%%%%%%%%%%%%%%%%%%%%%%%%%%%%%%%%%%%%%%%%%%%%%%%%%%%%%%%%%%%%%%%%%%
%%%%%%%%%%%%%%%%%%%%%%%%%%%%%%%%%%%%%%%%%%%%%%%%%%%%%%%%%%%%%%%%%%%%%%%%%%%%%%%%%%%%%%%%%%%%%%%%%%%%%%%%%%%%

%%%%%%%%%%%%%%%%%%%%%%%%%%%%%%%%%%%%%%%%%%%%%%%%%%%%%%%%%%%%%%%%%%%%%%%%%%%%%%%%%%%%%%%%%%%%%%%%%%%%%%%%%%%%
%%%%%%%%%%%%%%%%%%%%%%%%%%%%%%%%%%%%%%%%%%%%%%%%%%%%%%%%%%%%%%%%%%%%%%%%%%%%%%%%%%%%%%%%%%%%%%%%%%%%%%%%%%%%
%%%%%%%%%%%%%%%%%%%%%%%%%%%%%%%%%%%%%%%%%%%%%%%%%%%%%%%%%%%%%%%%%%%%%%%%%%%%%%%%%%%%%%%%%%%%%%%%%%%%%%%%%%%%

\section{Differential Graded Modules}

\subsection{Chain and Cochain Complexes} \label{chain}

\subsubsection{Graded Modules}
We work over a commutative associative ring with unit $R$.

A \emph{graded module} is a collection of $R$-modules
$$
M_*=\{..., M_{-1},\ M_0,\ M_1,\ ...\ ,\ M_n,\ M_{n+1},\ ...\ \}\ .
$$
A \emph{morphism of graded modules} $M_*\to M'_*$ is a collection
of homomorphisms $\{f_i:M_i\to M'_i,\ i\in Z\}$.

Sometimes we use the following notion: a \emph{morphism of graded
modules of degree }$n$ is a collection of homomorphisms
$\{f_i:M_i\to M'_{i+n},\ i\in Z\}$. So a morphism of graded
modules has the degree $0$.

\subsubsection{Chain Complexes}
\begin{definition}.
A differential graded ($dg$) module (or a chain complex) is a
sequence of
 $R$ modules and homomorphisms
$$
...\stackrel{d_{-1}}{\leftarrow }C_{-1}\stackrel{d_0}{\leftarrow
}C_0
\stackrel{d_1}{\leftarrow }C_1\stackrel{d_2}{\leftarrow }...\stackrel{%
d_{n-1}}{\leftarrow }C_{n-1}\stackrel{d_n}{\leftarrow }C_n\stackrel{d_{n+1}}{%
\leftarrow }C_{n+1}\stackrel{d_{n+2}}{\leftarrow }...\
$$
such that $d_id_{i+1}=0.$
\end{definition}

Elements of $C_n$ are called $n$-\emph{dimensional chains}; the
homomorphisms $d_i $ are called \emph{boundary operators}, or {\em
differentials}; elements of $Z_n=\Ker\ d_n\subset C_n $ are called
$n$-\emph{dimensional cycles} and elements of $ B_n=\Im\
d_{n+1}\subset C_n $ are called $n$-\emph{dimensional boundaries}.

It follows from the condition $d_id_{i+1}=0\,$ that $B_n\subset
Z_n.$

\begin{definition}.
The $n${\em -th homology module} $H_n(C_{*})$ of a dg module $%
(C_{*},d_{*})$ is defined as the quotient $Z_n/B_n.$
\end{definition}

A sequence $ C_{n-1}\stackrel{d_n}{\leftarrow
}C_n\stackrel{d_{n+1}}{\leftarrow }C_{n+1} $ is \emph{exact}, that
is, $B_n=Z_n,$ iff $H_n(C_{*})=0.$ Thus homology measures the
deviation from exactness.

\subsubsection{Cochain Complexes}

The notion of \emph{cochain complex} differs from the notion of
chain complex by the direction of the differential
$$
...\stackrel{d^{-1}}{\rightarrow }C^{-1}\stackrel{d^0}{\rightarrow
}C^0
\stackrel{d^1}{\rightarrow }C^1\stackrel{d^2}{\rightarrow }...\stackrel{%
d^{n-1}}{\rightarrow }C^{n-1}\stackrel{d^n}{\rightarrow }C^n\stackrel{d^{n+1}}{%
\rightarrow }C^{n+1}\stackrel{d^{n+2}}{\rightarrow }...\ \ .
$$
Corresponding terms here are cochains, cocycles $Z^n=\Ker \ d_{n+1}\subset C^n $, coboundaries $ B^n=\Im \  d^{n}\subset C^n $, cohomology $H^n(C^{*})=Z^n/B^n$.

Changing indices $C^n=C_{-n},\ d^n=d_{-n}$ we  convert a chain
complex $(C_*,d_*)$ to a cochain complex $(C^*,d^*)$.

\subsubsection{Dual Cochain Complex}

For a chain complex $(C_*,d_*)$ and an $R$-module $A$ the
\emph{dual cochain complex} $C^*=(\Hom(C_*,A),\delta^*)$ is defined
as
$$
C^n=\Hom(C_n,A),\ \delta^*), \ \delta^*(\phi)=\phi d.
$$

\subsubsection{Chain maps}

\begin{definition}.
A {\em chain map} of chain complexes $f:(C_{*},d_{*})\rightarrow
(C_{*}^{\prime },d_{*}^{\prime })$
 is defined as a sequence of homomorphisms
$\{f_i:C_i\rightarrow C_i^{\prime }\} $ such that $d_n^{\prime
}f_n=f_{n-1}d_n$.
\end{definition}

This condition means the commutativity of the diagram%
$$
\begin{array}{ccccccccc}
... & \stackrel{d_{n-1}}{\leftarrow } & C_{n-1} &
\stackrel{d_n}{\leftarrow }
& C_n & \stackrel{d_{n+1}}{\leftarrow } & C_{n+1} & \stackrel{d_{n+2}}{\leftarrow } & ... \\
&  & \downarrow f_{n-1}&  & \downarrow f_n &  & \downarrow f_{n+1} \\
... & \stackrel{d_{n-1}^{\prime }}{\leftarrow } & C_{n-1}^{\prime
} &
\stackrel{d_n^{\prime }}{\leftarrow } & C_n^{\prime } & \stackrel{%
d_{n+1}^{\prime }}{\leftarrow } & C'_{n+1} &
\stackrel{d'_{n+2}}{\leftarrow } &...\ .
\end{array}
$$

\begin{proposition}.
The composition of chain maps is a chain map.
\end{proposition}

Chain complexes and chain maps form a category, which we denote by
$DGMod$.

\begin{proposition}.
If $\left\{ f_i\right\} :(C_{*},d_{*})\rightarrow (C_{*}^{\prime
},d_{*}^{\prime })$ is a chain map, then $f_n$ sends cycles to
cycles and
boundaries to boundaries, i.e., $f_n(Z_n)\subset Z_n^{\prime }$ and $%
f_n(B_n)\subset B_n^{\prime }.$
\end{proposition}

\begin{proposition}.
A chain map $\left\{ f_i\right\} :(C_{*},d_{*})\rightarrow
(C_{*}^{\prime
},d_{*}^{\prime })$ induces the well defined homomorphism of homology groups%
$$
f_n^{*}:H_n(C_{*})\rightarrow H_n(C_{*}^{\prime }).
$$
\end{proposition}

Homology is a functor from the category of $dg$ modules to the
category of graded modules
$$
H:DGMod\to GMod.
$$

\subsubsection{$Hom$ Complex}

For two chain complexes $C,\ C'$ define the chain complex
$(Hom(C,C'),D)$ as
$$
Hom(C,C')_n=Hom^n(C,C')
$$
where $Hom^m(C,C')=\{\phi:C_*\to C'_{*+n},\}$ is the module of
homomorphisms of degree $n$, and the differential
$D:Hom^n(C,C')\to Hom^{n-1}(C,C')$ is given by $
D(\phi)=d'\phi+(-1)^{\deg\ \phi}\phi d. $

\subsubsection{Tensor Product}

For two chain complexes $A$ and $B$ the tensor product $A\otimes
B$ is defined as the following chain complex:
$$
(A\otimes B)_n=\sum_{p+q=n}A_p\otimes B_q,
$$
with differential $d_\otimes:(A\otimes B)_n\to (A\otimes B)_{n-1}$
given by
$$
d_{\otimes}(a_p\otimes b_q)=d_p(a_p)\otimes b_q+(-1)^{p}a_p\otimes
d'_q(b_q).
$$

If $f:A\to A'$ and $g:B\to B'$ are chain maps then there is a
chain map
$$
f\otimes g:A\otimes B\to A'\otimes B'
$$
defined as $ (f\otimes g)(a\otimes b)=f(a)\otimes g(b). $

\subsubsection{Chain Homotopy}

\begin{definition}.
Two chain maps $\left\{ f_i\right\} ,\{g_i\}:(C_{*},d_{*})\rightarrow
(C_{*}^{\prime },d_{*}^{\prime })$ are called {\em chain
homotopic}, if there exists a sequence of homomorphisms
$D_n:C_n\rightarrow C_{n+1}^{\prime }, $
$$
\begin{array}{ccccccccc}
... & \stackrel{d_{n-1}}{\longleftarrow } & C_{n-1} &
\stackrel{d_n}{\longleftarrow }
& C_n & \stackrel{d_{n+1}}{\longleftarrow } & C_{n+1} & \stackrel{d_{n+2}}{\longleftarrow } & ... \\
&  & f_{n-1}\downarrow \downarrow g_{n-1}& \searrow D_{n-1} &
f_n\downarrow \downarrow g_n &
 \searrow D_n & f_{n+1}\downarrow \downarrow g_{n+1} \\
... & \stackrel{d_{n-1}^{\prime }}{\leftarrow } & C_{n-1}^{\prime
} &
\stackrel{d_n^{\prime }}{\longleftarrow } & C_n^{\prime } & \stackrel{%
d_{n+1}^{\prime }}{\longleftarrow } & C'_{n+1} &
\stackrel{d'_{n+2}}{\longleftarrow } &...\ .
\end{array}
$$
such that $ f_n-g_n=d_{n+1}^{\prime }D_n+D_{n-1}d_n. $ In this
case we write $f \sim_{D}g$.
\end{definition}

\begin{proposition}.
Chain homotopy is an equivalence relation:
$$
\begin{array}{cl}
(a)& \ \ f\sim_0 f;\\
(b)& \ \ f \sim_{D}g \Longrightarrow g \sim_{-D}f;\\
(c)& \ \ f \sim_{D}g,\ g \sim_{D'}h \Longrightarrow f
\sim_{D+D'}h.
\end{array}
$$
\end{proposition}

\begin{proposition}.
Chain homotopy is compatible with compositions:
$$
\begin{array}{cl}
(a)& \ \ f\sim_D g \Longrightarrow hf\sim_{hD}hg;\\
(b)& \ \ f \sim_{D}g \Longrightarrow fk \sim_{Dk}gk.\\
\end{array}
$$
\end{proposition}

Thus there is a category $hoDGMod$ whose objects are chain
complexes and morphisms are chain homotopy classes
$$
Hom_{hoDGMod}(C,C')=[C,C']=Hom_{DGMod}(C,C')/\sim.
$$

\begin{proposition}.
If two chain maps $\left\{ f_i\right\}
,\{g_i\}:(C_{*},d_{*})\rightarrow (C_{*}^{\prime },d_{*}^{\prime
})$ are chain homotopic, then the induced
homomorphisms of homology groups coincide: $f_n^{*}=g_n^{*}:H_n(C_{*})%
\rightarrow H_n(C_{*}^{\prime }).$

Thus we have the commutative diagram of functors
$$
\begin{array}{ccccc}
DGMod & & \longrightarrow & & hoDGMod \\
  & H\searrow & & \swarrow H & \\
  & & GMod. & &
\end{array}
$$
\end{proposition}

\subsubsection{Chain Equivalence}

Chain complexes $C$ and $C'$ are called \emph{chain equivalent}
$C\sim C'$, if there exist chain maps
$$
f:C\stackrel{\longleftarrow}{\to}C':g
$$
such that $gf\sim \id_C,\ \ fg\sim \id_{C'}$. This means that $C$
and $C'$ are isomorphic in $hoDGMod$.

A chain complex $C$ is called \emph{contractible} if $C\sim 0$,
equivalently if $id_C\sim 0:C\to C$.

\begin{proposition}.
Each contractible $C$ is acyclic, i.e., $H_i(C)=0$ for all $i$.
\end{proposition}

\begin{proposition}. If all $C_i$ are free and $C$ is acyclic, then
$C$ is contractible.
\end{proposition}

\subsection{Algebraic and Topological Examples}

\subsubsection{Algebraic Example}

Let $(A,\mu:A\otimes A\to A)$ be an associative algebra, then
$$
C(A)=(A \stackrel{\mu}{\leftarrow } A\otimes A
\stackrel{\mu\otimes id-id\otimes \mu}{\leftarrow } A\otimes
A\otimes A \stackrel{\mu\otimes id\otimes id-id\otimes \mu\otimes
id+id\otimes id\otimes \mu}{\leftarrow }...)
$$
is a chain complex: the associativity condition guarantees that
$dd=0$.

If $A$ has a unit $e\in A$ then this complex is contractible, that
is $id:C(A)\to C(A)$ and $0:C(A)\to C(A)$ are homotopic: the
suitable chain homotopy is given by $D(a_1\otimes ...\otimes
a_n)=(e\otimes a_1\otimes ...\otimes a_n)$. This immediately
implies that $C(A)$ is acyclic, that is $H_i(C(A))=0$ for all
$i>0$.

This example is a particular case of more general chain complex
called the bar construction, see later.

\subsubsection{Simplicial Complexes}

Simplicial complex is a formal construction, which models
topological spaces.

\begin{definition}.
A {\em simplicial complex} is a set $V$ with a given family of
finite subsets, called simplices, so that the following conditions
are satisfied:

(1) all points (called \emph{vertices}) of $V$ are simplices;

(2) any nonempty subset of a simplex is a simplex.
\end{definition}

A simplex consisting of $(n+1)$ points is called $n${\em %
-dimensional simplex}. The $0$-dimensional simplexes, i.e., the
points of $V$ are called {\em vertices}.

\begin{definition}.
A {\em simplicial map} of simplicial complexes $V\rightarrow
V^{\prime }$ is a map of vertices $f:V\rightarrow V^{\prime }$
such that the image of any simplex of $V$ is a simplex in
$V^{\prime }.$
\end{definition}

\begin{proposition}.
The composition of simplicial maps is a simplicial map.
\end{proposition}

Simplicial complexes and simplicial maps form a category which we
denote as $SC$.

\begin{proposition}.
To any simplicial set $V$ corresponds a topological space $| V| $
(called its {\em realization}) and to any simplicial map
$f:V\rightarrow V^{\prime }$ corresponds a continuous map of
realizations $| f| :|V| \rightarrow | V^{\prime }| .$
\end{proposition}

So the realization is a functor from the category of simplicial
complexes  to the category of topological spaces
$$
|- |:SC\to Top.
$$

\subsubsection{Homology Modules of a Simplicial Complex} \label{HSimpl}

In this section we consider \emph{ordered } simplicial complexes:
we assume that the set of vertices $V$ is ordered by a certain order.

We assign to such an ordered simplicial complex the following chain complex $%
(C_{*}(V),d_{*})$: Let $C_n(V)$ be the free $R$-module, generated
by all ordered $n$-simplices $\sigma
_n=(v_{k_0},v_{k_1},...,v_{k_n}),$ where
$v_{k_0}<v_{k_1}<...<v_{k_n}$; the differential $
d_n:C_n(V)\rightarrow C_{n-1}(V) $ on a generator $\sigma
_n=(v_{k_0},v_{k_1},...,v_{k_n})\in C_n(V)$ is given by
$$
d_n(v_{k_0},v_{k_1},...,v_{k_n})=\sum_{i=0}^n(-1)^i(v_{k_0},...,\widehat{v}%
_{k_i},...,v_{k_n}),
$$
where $(v_{k_0},...,\widehat{v}_{k_i},...,v_{k_n})$ is the
$(n-1)$-simplex obtained by omitting $v_{k_i},$ and is
extended on the whole $C_n(V)$ linearly.

\begin{proposition}.
The composition $d_{n-1}d_n$ is zero, thus
$(C_{*}(V),d_{*})$ is a chain complex.
\end{proposition}

\begin{definition}.
The $n${\em -th homology group} $H_n(V)$ of an ordered simplicial
set $V$ is defined as the $n$-th homology group $H_n(C_{*}(V)).$
\end{definition}

%\begin{proposition}.
%Any monotonic simplicial map $f:V\rightarrow V^{\prime }$ induces
%a chain map $C_{*}(f):(C_{*}(V),d_{*})\rightarrow (C_{*}(V^{\prime
%}),d_{*}),$ and
%consequently it induces a homomorphism of homology groups $%
%f_n^{*}:H_n(V)\rightarrow H_n(V^{\prime }).$
%\end{proposition}

%\paragraph{Cohomology Modules of a Simplicial Complex.}

\subsubsection{Cohomology Modules of a Simplicial Complex}
Let $A$ be an $R$-module. The cochain complex of $V$ with
coefficients in $A$ is defined as the dual to the chain complex
$C_*(V)$: $ C^*(V,A)=Hom(C_*(V),A). $ The $n$-th cohomology module
of $V$ with coefficients in $A$ is just the $n$-th homology of
this cochain complex.

Below we show that the cohomology $H^*(V,A)$ is more interesting than the
homology $H_*(V)$, since cohomology possesses richer algebraic
structure: it is a ring.

\section{Differential Graded Algebras and Coalgebras}

\subsection{Differential Graded Algebras}

\subsubsection{Graded algebras} \label{Galg}

A \emph{graded algebra} is a graded module
$$
A_*=\{..., A_{-1},\ A_0,\ A_1,\ ...\ ,\ A_n,\ A_{n+1},\ ...\ \}
$$
equipped with an associative multiplication $\mu:A_*\otimes A_*\to A_{*}$ of degree 0, i.e., $\mu(\mu\otimes id)=\mu(id\otimes\mu)$ and $\mu:A_p\otimes A_q\to A_{p+q}$. We denote $a\cdot b=\mu(a\otimes b)$.

For a graded algebra $A_*$ the component $A_0$ is an
associative algebra.

A \emph{morphism of graded algebras }$f:(A,\mu)\to (A',\mu')$ is
a morphism of graded modules $\{f_k:A_k\to A'_k\}$ which is
multiplicative, that is, $f\mu=\mu(f\otimes f)$, i.e., $f(a\cdot b)=f(a)\cdot f(b)$.

%We'll use also the notion of $(f-g)$ \emph{k-derivation}: for the given two graded algebra morphisms
%$f,g:(A,\mu)\to (A',\mu')$
%this is a homomorphism  $D:(A,\mu)\to (A',\mu')$ of degree $k$, that is $D:A_n\to A_{n+k}$, which satisfies the condition
%$D(a\cdot b)=D(a)\cdot g(b)+f(a)\cdot D(b)$.

Let $f,g:(A,\mu)\to (A',\mu')$ be two morphisms of graded algebras. An
$(f,g)$-\emph{derivation} of degree $k$ is defined as a morphism of graded modules
of degree $k$ $D:A_*\to A'_{*+k}$, i.e., a collection of
homomorphisms $\{D_i:A_i\to A'_{i+k},\ i\in Z\}$ which satisfies
the condition
$$
D(a\cdot b)=D(a)\cdot g(b)+(-1)^{k\cdot |a|}f(a)\cdot D(b).
$$

An essential particular case of this notion is a
 \emph{k-derivation} $D: (A,\mu)\to (A,\mu)$ which is an $(id,id)$ \emph{k-derivation}, i.e., it satisfies the condition
 $D(a\cdot b)=D(a)\cdot b+a\cdot D(b)$.

 It is easy to see that if $D, D':A_*\to A'_{*+k}$ are two $(f,g)$ derivation, then their sum $D+D'$ is also an $(f,g)$ derivation.

 Moreover, for graded algebra morphisms $h:B\to A,\ l:A'\to C$ and an $(f,g)$ derivations $D:A\to A'$ the composition
 $Dh$ is an $(fh,gh)$ derivation and $lD$ is an $(lf,lg)$-derivation.

\subsubsection{Differential graded algebras} \label{DGA}

\begin{definition}. A differential graded algebra (dga in short)
$(A,d,\mu)$ is a dg module $(A,d)$ equipped additionally with a
multiplication
$$
\mu:A\otimes A\to A
$$
so that $(A,\mu)$ is a graded algebra, and the multiplication
$\mu$ is a chain map, that is, the differential $d$ and $\mu$ are connected by the condition
$$
d(a\cdot b)=da\cdot b +(-1)^{|a|}a\cdot db.
$$
This condition means simultaneously that $\mu$ is a chain map, and
that $d$ is an $(id,id)$-\emph{derivation} of degree $-1$.

A morphism of dga-s $f:(A,d,\mu)\to (A,d,\mu)$ is defined as a
multiplicative chain map:
$$
df=fd,\ \ f(a\cdot b)=f(a)\cdot f(b).
$$
\end{definition}

We denote the obtained category  as $DGAlg$.

For a dg algebra $(A,d,\mu)$ its homology $H_*(A)$ is a graded
algebra with the following multiplication:
$$
H_*(A)\otimes H_*(A)\stackrel{\phi}{\longrightarrow} H_*(A\otimes
A)\stackrel{H_*(\mu)}{\longrightarrow}H_*(A),
$$
where $\phi:H_*(A)\otimes H_*(A)\to H_*(A\otimes A)$ is the
standard map
$$
\phi(h_1\otimes h_2)=cl(z_{h_1}\otimes z_{h_2}).
$$

In other words the multiplication on $H(A)$ is defined as follows:
For $h_1,h_2\in H(A)$ the product $h_1\cdot h_2$ is the homology
class of the cycle $z_{h_1}\cdot z_{h_2}$.

Furthermore, a dga map induces a multiplicative map of homology
graded algebras.

Thus homology is a functor from the category of dg algebras to the
category of graded algebras.

\subsubsection{Derivation homotopy}

Two dg algebra maps $f,g:A\to A'$ are called homotopic if there
exists a chain homotopy $D:A\to A'$, $ f-g=dD+Dd $ which, in
addition is a $(f,g)$-derivation, that is
$$
D(a\cdot b)=D(a)\cdot g(b)+(-1)^{|a|}f(a)\cdot D(b).
$$
Note that generally this is not an equivalence relation.

\subsection{Differential Graded Colgebras}

\subsubsection{Graded coalgebras} \label{Gcoalg}

A graded coalgebra $(C,\Delta)$ is a graded module
$$
C=\{\ ...,\ C_{-1},\ C_0,\ C_1,\ ...\ ,\ C_n,\ C_{n+1},\ ...\ \}
$$
equipped with a comultiplication
$$
\Delta:C\otimes C\to C
$$
which is coassociative, that is $(\Delta\otimes
id)\Delta=(id\otimes \Delta)\Delta$, i.e., the following diagram commutes:
$$
\begin{array}{rcccl}
  &C & \stackrel{\Delta}{\longrightarrow}&C\otimes C&\\
 \Delta & \downarrow &  & \downarrow &\Delta\otimes id \\

 &  C\otimes C &  \stackrel{id\otimes
  \Delta}{\longrightarrow}& C\otimes C\otimes C&\ .
\end{array}
$$

A morphism of graded coalgebras $f:(C,\Delta)\to (C',\Delta')$ is
a morphism of graded modules $\{f_k:C_k\to C'_k\}$ which is
comultiplicative, that is $\Delta'f=(f\otimes f)\Delta$, i.e., the following diagram commutes:
$$
\begin{array}{rcccl}
  &C & \stackrel{f}{\longrightarrow}&C'&\\
 \Delta & \downarrow &  & \downarrow &\Delta' \\
 &  C\otimes C &  \stackrel{f\otimes f}{\longrightarrow}& C'\otimes C' & .
\end{array}
$$

If $f,g:C\to C'$ are two morphisms of graded coalgebras, then an
$(f,g)$-coderivation of degree $k$ is defined as a collection of
homomorphisms $\{D_i:C_i\to C_{i+k}\}$ which satisfies the
condition $ \Delta'D=(f\otimes D+D\otimes g)\Delta, $ i.e., the following diagram commutes:
$$
\begin{array}{rcccl}
  &C & \stackrel{D}{\longrightarrow}&C'&\\
 \Delta & \downarrow &  & \downarrow &\Delta' \\
 &  C\otimes C &  \stackrel{f\otimes D+D\otimes g}{\longrightarrow}& C'\otimes
 C'&.
\end{array}
$$

%%%%%%%%%%%%%%%%%%%%%%%%%%%%%%%%%%%%%%%%%%%%%%%%%%%%%%%%%%%%

%%%%%%%%%%%%%%%%%%%%%%%%%%%%%%%%%%%%%%%%%%%%%%%%%%%%%%%%%%%%

%%%%%%%%%%%%%%%%%%%%%%%%%%%%%%%%%%%%%%%%%%%%%%%%%%%%%%%%%%%%

An essential particular case of this notion is that of a
 \emph{k-coderivation} $D: (C,\Delta)\to (C,\Delta)$, which is an $(id,id)$- \emph{k-coderivation}, i.e., it satisfies the condition
$ \Delta D=(id\otimes D+D\otimes id)\Delta$.

 It is easy to see that if $D, D':C_*\to C'_{*+k}$ are two $(f,g)$- coderivations, then their sum $D+D'$ is also an $(f,g)$-coderivation.

 Moreover, for graded coalgebra morphisms $h:A\to C,\ l:C'\to B$ and an
 $(f,g)$-coderivation $D:C\to C'$ the composition
 $Dh$ is an $(fh,gh)$-coderivation and $lD$ is an $(lf,lg)$ coderivation.

\subsubsection{Differential Graded Coalgebras}  \label{DGC}

\begin{definition}. A differential graded coalgebra (dgc in short)
$(C,d,\Delta)$ is a dg module $(C,d)$ equipped additionally with a
comultiplication $\Delta:C\to C\otimes C$ so that $(C,\Delta)$ is
a graded coalgebra and the comultiplication $\Delta$ and the
differential $d$ are related by the condition
$$
\Delta d=(d\otimes id + id\otimes d)\Delta .
$$
This condition means simultaneously that $\Delta$ is a chain map,
and that $d$ is a $(id,id)$-coderivation of degree $-1$.

A morphism of dgc-s $f:(C,d,\Delta)\to (C',d',\Delta')$ is defined
as a morphism of graded coalgebras which is a chain map.
\end{definition}

We denote the obtained category as $DGCoalg$.

Generally for a dg coalgebra $(C,d,\mu)$ its homology $H_*(C)$
\emph{is not} a graded coalgebra:
$$
H_*(C)\otimes H_*(C)\stackrel{\phi}{\longrightarrow} H_*(C\otimes
C)\stackrel{H(\Delta)}{\longleftarrow}H_*(C),
$$
the map $\phi:H_*(C)\otimes H_*(C)\to H_*(C\otimes C)$ has the wrong
direction, but if all $H_i(C)$ are free then $\phi$ is invertible
and $H_*(C)$ is a graded coalgebra.

\subsubsection{Coderivation homotopy}

In the category of dg coalgebras there is the following notion of homotopy: two dg coalgebra maps
$f,g:(C,d_C,\Delta_C)\to (C',d_{C'},\Delta_{C'})$
are homotopic, if there exists $D:C\to C'$ of degree $+1$ such that
$d_{C'}D+Dd_C=f-g$, i.e., the chain maps $f$ and $g$ are chain homotopic, and additionally the homotopy $D$ is a $f-g$-coderivation, that is
$\Delta_{C'}D=(f\otimes D+D\otimes g)\Delta_C.$

\subsubsection{Dual cochain algebra}

Let $(C_*,d_C,\Delta_C)$ be a dg coalgebra and let
$(A,\mu_A:A\otimes A\to A)$ be a (nondifferential and nongraded) algebra. Then the dual cochain complex
$(C^*=Hom(C_*,A), \delta^*),\ \delta^*(\phi)=\phi d_C$, becomes a dg algebra (of cochain type, i.e., the degree of the differential is +1)
with the multiplication (cup product)
$$
\mu^*:Hom(C_*,A)\otimes Hom(C_*,A)\to Hom(C_*,A)
$$
given by
$$
\phi\smile \psi=\mu^*(\phi\otimes\psi)=\mu_A(\phi\otimes \psi)\Delta_C.
$$
 The obtained object $(C^*=Hom(C_*,A),\ \delta^*(\phi)=\phi d_C,\ \smile)$
is called the cochain dg algebra of the chain dg coalgebra $(C_*,d_C,\Delta_C)$ with coefficients in $A$.

%%%%%%%%%%%%%%%%%%%%%%%%%%%%%%%%%%%%%%%%%%%%%%%%%%%%
%%%%%%%%%%%%%%%%%%%%%%%%%%%%%%%%%%%%%%%%%%%%%%%%%%%%
%%%%%%%%%%%%%%%%%%%%%%%%%%%%%%%%%%%%%%%%%%%%%%%%%%%%
%%%%%%%%%%%%%%%%%%%%%%%%%%%%%%%%%%%%%%%%%%%%%%%%%%%%
%%%%%%%%%%%%%%%%%%%%%%%%%%%%%%%%%%%%%%%%%%%%%%%%%%%%

\subsection{Applications in Topology}

\subsubsection{Alexander-Whitney diagonal}

Let $V$ be an ordered simplicial complex (\ref{HSimpl}) and let $(C_*(V),d_*)$ be its chain
complex. There exists a comultiplication
$$
\Delta:C_*(V)\to C_*(V)\otimes C_*(V),
$$
the so called Alexander-Whitney diagonal, which turns $(C_*(V),d_*,\Delta)$ into a
dg coalgebra. This diagonal is defined by
$$
\Delta(v_{i_0},...,v_{i_n})=\sum_{k=0}^n
(v_{i_0},...,v_{i_k})\otimes (v_{i_k},...,v_{i_n}).
$$

\subsubsection{Cohomology algebra}

Let again $(A,\mu_A:A\otimes A\to A)$ be a (nondifferential and nongraded) associative algebra.
The Alexander-Whitney diagonal of $C_*(V)$ induces on the dual cochain complex $C^*(V,A)=Hom(C^*(V),A),\delta^*)$
the cup product
$$
\smile :C^*(V)\otimes C^*(V)\to C^*(V)
$$
which for $\phi\in C^p(V,A),\ \psi\in C^q(V,A))$ looks as
$$
\phi\smile \psi
(v_{i_0},...,v_{i_{p+q}})=\phi(v_{i_0},...,v_{i_p})\cdot
\psi(v_{i_p},...,v_{i_{p+q}}).
$$
This turns $C^*(V,A),\delta^*,\smile)$ into a dg algebra.

This structure induces on the cohomology $H^*(V,A)$ a structure of
graded algebra.

The cohomology algebra $H^*(V,A)$ is a more powerful invariant than the
cohomology groups: the two spaces $X=S^1\times S^1$ and $Y=S^1\vee S^1\vee
S^2$ have the same cohomology groups
$$
H^0=R,\ \ H^1=R\cdot a\oplus R\cdot b,\ \ H^2=R\cdot c,
$$
with generators $a,\ b$ in dimension $1$ and $c$ in dimension $2$,
but they have different cohomology algebras, namely $a\cdot b=0$ in
$H^*(Y)$ and $a\cdot b=c$ in $H^*(X)$.

\section{Bar and Cobar Functors} \label{bar}

Here we describe two classical adjoint functors

$$
B:DGAlg \stackrel{\longrightarrow}{\leftarrow} DGCoalg:\Omega,
$$
the bar functor $B:DGAlg\to DGCoalg$ from the category of dg
algebras to the category of dg coalgebras, and the cobar functor
$\Omega:DGCoalg\to DGAlg$ in the opposite direction.

We start with the definitions of free (in the category of graded algebras) and cofree  (in the category of graded coalgebras) objects.

\subsection{Tensor Algebra and Tensor Coalgebra}

\subsubsection{Tensor Algebra}

Let $V=\{V_i\}$ be a graded $R$-module. The tensor algebra
generated by $V$ is defined as
$$
T(V)=R\oplus V\oplus (V\otimes V)\oplus (V\otimes V\otimes V)\oplus
...=\sum_{i=0}^\infty V^{\otimes i}
$$
with grading $ dim (a_1\otimes ... \otimes a_m)=dim\ a_1+...+dim\
a_m, $ and with \emph{multiplication} $\mu:T(V)\otimes T(V)\to T(V)$ given by
$$
\mu((a_1\otimes...\otimes a_m)\otimes (a_{m+1}\otimes ...\otimes
a_{m+n}))=a_1\otimes ...\otimes a_{m+n}.
$$
The unit element for this multiplication is $1\in R=V^{\otimes
0}$.

By $i_k$ we denote the obwious inclusion $i_k:V^{\otimes k}\to
T(V)$.

\vspace{5mm}

\subsubsection{ Universal Property of $T(V)$}  \label{UnivTV}
The tensor algebra $T(V)$ is the free object in the category of graded
algebras: for an arbitrary graded algebra $(A,\mu_A)$ and a map of graded
modules $\alpha :V\to A$ there exists a unique morphism of graded
algebras $f_{\alpha}:T(V)\to A$ such that $f_\alpha (v)=\alpha(v)$
(i.e., $f_{\alpha}i_1=\alpha$).

This morphism $f_\alpha$ (which is called \emph{multiplicative extension}
of $\alpha$) is defined as $f_\alpha (a_1\otimes ... \otimes a_m)=\alpha
(a_1)\cdot ...\cdot \alpha (a_m). $ Or, equivalently $f_{\alpha}$
is described by:
$$
f_\alpha i_k=\sum_k \mu_A^k(\alpha\otimes...\otimes \alpha),
$$
 where $\mu_A^k:A^{\otimes k}\to A$ is the $k$-fold iteration of the multiplication
$\mu_A :A\otimes A\to A$, namely, $\mu_A^1=id,\ \mu_A^2=\mu_A,\
\mu_A^k=\mu_A(\mu_A^{k-1}\otimes id)$.

So, to summarize, we have the following universal property
$$
\begin{array}{rcl}

V & \stackrel{i_1}{\longrightarrow } & T(V)\\
   \alpha&   \searrow & \downarrow f_{\alpha}=\sum_k \mu_A^k(\alpha\otimes...\otimes \alpha)\\
  &           & A.
\end{array}
$$

\vspace{5mm}

 \subsubsection{ Universal Property for Derivations}  \label{UnivTVder}

%%%%%%%%%%%%%%%%%%%%%%%%%%%%%%%%%%%%%%%%%%%%%%%%%%%%%%%%%%%%%%%%%%%%%%%%%%%%%%%%%
%%%%%%%%%%%%%%%%%%%%%%%%%%%%%%%%%%%%%%%%%%%%%%%%%%%%%%%%%%%%%%%%%%%%%%%%%%%%%%%%%
%%%%%%%%%%%%%%%%%%%%%%%%%%%%%%%%%%%%%%%%%%%%%%%%%%%%%%%%%%%%%%%%%%%%%%%%%%%%%%%%%
%%%%%%%%%%%%%%%%%%%%%%%%%%%%%%%%%%%%%%%%%%%%%%%%%%%%%%%%%%%%%%%%%%%%%%%%%%%%%%%%%
%%%%%%%%%%%%%%%%%%%%%%%%%%%%%%%%%%%%%%%%%%%%%%%%%%%%%%%%%%%%%%%%%%%%%%%%%%%%%%%%%
%%%%%%%%%%%%%%%%%%%%%%%%%%%%%%%%%%%%%%%%%%%%%%%%%%%%%%%%%%%%%%%%%%%%%%%%%%%%%%%%%
%%%%%%%%%%%%%%%%%%%%%%%%%%%%%%%%%%%%%%%%%%%%%%%%%%%%%%%%%%%%%%%%%%%%%%%%%%%%%%%%%
%%%%%%%%%%%%%%%%%%%%%%%%%%%%%%%%%%%%%%%%%%%%%%%%%%%%%%%%%%%%%%%%%%%%%%%%%%%%%%%%%
%%%%%%%%%%%%%%%%%%%%%%%%%%%%%%%%%%%%%%%%%%%%%%%%%%%%%%%%%%%%%%%%%%%%%%%%%%%%%%%%%
%%%%%%%%%%%%%%%%%%%%%%%%%%%%%%%%%%%%%%%%%%%%%%%%%%%%%%%%%%%%%%%%%%%%%%%%%%%%%%%%%
%%%%%%%%%%%%%%%%%%%%%%%%%%%%%%%%%%%%%%%%%%%%%%%%%%%%%%%%%%%%%%%%%%%%%%%%%%%%%%%%%
%%%%%%%%%%%%%%%%%%%%%%%%%%%%%%%%%%%%%%%%%%%%%%%%%%%%%%%%%%%%%%%%%%%%%%%%%%%%%%%%%

The tensor algebra has an analogous universal property also for
derivations: for a graded algebra $(A,\mu)$, two homomorphisms $\alpha, \alpha':V\to A$ of degree $0$ and a homomorphism $\beta:V\to A$ of degree $k$
there exist: morphisms of graded algebras $f_\alpha, f_{\alpha'}:T(V)\to A$ and a unique $(f_\alpha, f_{\alpha'})$ derivation of degree $k$
$$
D_\beta:T(V)\to A
$$
such that $f_\alpha (v)=\alpha (v)$, $f_\alpha'(v)=\alpha' (v)$ and  $D(v)=\beta (v)$, i.e., the following diagram commutes.
$$
\begin{array}{rcl}

V & \stackrel{i_1}{\longrightarrow } & T(V)\\
   \beta&   \searrow & \downarrow D \\
  &           & A.
\end{array}
$$

The derivation $D$ is defined as
$$
D (a_1\otimes ... \otimes a_n)=
\sum_{k=1}^n \mu^n(\alpha(a_1)\otimes ... \otimes \alpha(a_{k-1})\otimes \beta(a_k)\otimes \alpha'(a_{k+1})\otimes  ... \otimes \alpha'(a_n)).
$$
Or, equivalently $ D \cdot i_n=\sum_{k} \mu^n(\alpha^{\otimes
(k-1)}\otimes \beta \otimes \alpha'^{\otimes (n-k)})i_n. $

%%%%%%%%%%%%%%%%%%%%%%%%%%%%%%%%%%%%%%%%%%%%%%%%%%%%%%%%%%%%%%%%%%%%%%%%%%%%%%%%%
%%%%%%%%%%%%%%%%%%%%%%%%%%%%%%%%%%%%%%%%%%%%%%%%%%%%%%%%%%%%%%%%%%%%%%%%%%%%%%%%%
%%%%%%%%%%%%%%%%%%%%%%%%%%%%%%%%%%%%%%%%%%%%%%%%%%%%%%%%%%%%%%%%%%%%%%%%%%%%%%%%%
%%%%%%%%%%%%%%%%%%%%%%%%%%%%%%%%%%%%%%%%%%%%%%%%%%%%%%%%%%%%%%%%%%%%%%%%%%%%%%%%%
%%%%%%%%%%%%%%%%%%%%%%%%%%%%%%%%%%%%%%%%%%%%%%%%%%%%%%%%%%%%%%%%%%%%%%%%%%%%%%%%%
%%%%%%%%%%%%%%%%%%%%%%%%%%%%%%%%%%%%%%%%%%%%%%%%%%%%%%%%%%%%%%%%%%%%%%%%%%%%%%%%%
%%%%%%%%%%%%%%%%%%%%%%%%%%%%%%%%%%%%%%%%%%%%%%%%%%%%%%%%%%%%%%%%%%%%%%%%%%%%%%%%%
%%%%%%%%%%%%%%%%%%%%%%%%%%%%%%%%%%%%%%%%%%%%%%%%%%%%%%%%%%%%%%%%%%%%%%%%%%%%%%%%%
%%%%%%%%%%%%%%%%%%%%%%%%%%%%%%%%%%%%%%%%%%%%%%%%%%%%%%%%%%%%%%%%%%%%%%%%%%%%%%%%%
%%%%%%%%%%%%%%%%%%%%%%%%%%%%%%%%%%%%%%%%%%%%%%%%%%%%%%%%%%%%%%%%%%%%%%%%%%%%%%%%%
%%%%%%%%%%%%%%%%%%%%%%%%%%%%%%%%%%%%%%%%%%%%%%%%%%%%%%%%%%%%%%%%%%%%%%%%%%%%%%%%%
%%%%%%%%%%%%%%%%%%%%%%%%%%%%%%%%%%%%%%%%%%%%%%%%%%%%%%%%%%%%%%%%%%%%%%%%%%%%%%%%%

\subsubsection{Tensor Coalgebra}  \label{TcV}

Here is the dualization of the previous notion.

Let $V=\{V_i\}$ be a graded $R$-module. The \emph{tensor coalgebra}
cogenerated by $V$ is defined (again) as
$$
T^c(V)=R\oplus V\oplus V\otimes V\oplus V\otimes V\otimes V\oplus
...=\sum_{i=0}^\infty V^{\otimes i}
$$
with the same grading $ dim (a_1\otimes ... \otimes a_m)=dim \ a_1+...+dim\
a_m, $ but now with \emph{comultiplication} $\Delta:T^c(V)\to T^c(V)\otimes
T^c(V)$ given by
$$
\Delta (a_1\otimes. . . \otimes a_n)=\sum_{i=0}^{n}(a_1\otimes. .
. \otimes a_i)\otimes (a_{i+1}\otimes. . . \otimes a_n),
$$
here $(\ )=1\in R=V^{\otimes 0}$.

 By $p_k$ we denote the obwious projection $p_k:T^c(V)\to
V^{\otimes k}$.

\vspace{5mm}

 \subsubsection{ Universal Property of $T^c(V)$}  \label{UnivTcV}
In order to formulate the about universal property in this case we have to
introduce some dimensional restrictions in this case.

Let $V=\{...,0,0,V_1,V_2,...\}$ be a \emph{connected} graded
module, that is, $V_i=0$ for $i\leq 0$.

The tensor coalgebra of such $V$ is the cofree object in the
category of \emph{connected} graded coalgebras: for a map of graded
modules $\alpha :C\to V$ there exists a unique morphism of graded
coalgebras $f_{\alpha}:C\to T^c(V)$ such that
$p_1f_{\alpha}=\alpha$, i.e., the following diagram commutes
$$
\begin{array}{rcl}

V & \stackrel{p_1}{\longleftarrow } & T^c(V)\\
   \alpha&   \nwarrow & \uparrow f_{\alpha}\\
  &           & C.
\end{array}
$$

The coalgebra map $f_\alpha$ (which is called comultiplicative
coextension of $\alpha$) is defined as $ f_\alpha=\sum_k
(\alpha\otimes...\otimes \alpha)\Delta^k, $ where $\Delta^k:C\to
C^{\otimes k}$ is the $k$-th iteration of the comultiplication
$\Delta :C\to C\otimes C$, i.e., $\Delta^1=id,\ \Delta^2=\Delta,\
\Delta^k=(\Delta^{k-1}\otimes id)\Delta$.

%%%%%%%%%%%%%%%%%%%%%%%%%%%%%%%%%%%%%%%%%%%%%%%%%%%%
%%%%%%%%%%%%%%%%%%%%%%%%%%%%%%%%%%%%%%%%%%%%%%%%%%%%
%%%%%%%%%%%%%%%%%%%%%%%%%%%%%%%%%%%%%%%%%%%%%%%%%%%%
%%%%%%%%%%%%%%%%%%%%%%%%%%%%%%%%%%%%%%%%%%%%%%%%%%%%
%%%%%%%%%%%%%%%%%%%%%%%%%%%%%%%%%%%%%%%%%%%%%%%%%%%%

%%%%%%%%%%%%%%%%%%%%%%%%%%%%%%%%%%%%%%%%%%%%%%%%%%%%%%%%%%%%%%%%%%%%%%%%%%%%%%%%%
%%%%%%%%%%%%%%%%%%%%%%%%%%%%%%%%%%%%%%%%%%%%%%%%%%%%%%%%%%%%%%%%%%%%%%%%%%%%%%%%%
%%%%%%%%%%%%%%%%%%%%%%%%%%%%%%%%%%%%%%%%%%%%%%%%%%%%%%%%%%%%%%%%%%%%%%%%%%%%%%%%%
%%%%%%%%%%%%%%%%%%%%%%%%%%%%%%%%%%%%%%%%%%%%%%%%%%%%%%%%%%%%%%%%%%%%%%%%%%%%%%%%%
%%%%%%%%%%%%%%%%%%%%%%%%%%%%%%%%%%%%%%%%%%%%%%%%%%%%%%%%%%%%%%%%%%%%%%%%%%%%%%%%%
%%%%%%%%%%%%%%%%%%%%%%%%%%%%%%%%%%%%%%%%%%%%%%%%%%%%%%%%%%%%%%%%%%%%%%%%%%%%%%%%%
%%%%%%%%%%%%%%%%%%%%%%%%%%%%%%%%%%%%%%%%%%%%%%%%%%%%%%%%%%%%%%%%%%%%%%%%%%%%%%%%%
%%%%%%%%%%%%%%%%%%%%%%%%%%%%%%%%%%%%%%%%%%%%%%%%%%%%%%%%%%%%%%%%%%%%%%%%%%%%%%%%%
%%%%%%%%%%%%%%%%%%%%%%%%%%%%%%%%%%%%%%%%%%%%%%%%%%%%%%%%%%%%%%%%%%%%%%%%%%%%%%%%%

\vspace{5mm}
 \subsubsection{ Universal Property for Coderivations}  \label{UnivTcVder}

The tensor coalgebra has a similar universal property also for
coderivations.

Namely, for a graded coalgebra $(C,\Delta)$, two homomorphisms of degree $0$  $\alpha, \alpha' :C\to V$  and a homomorphism of degree $k$
$\beta :C\to V$ there exist: morphisms of graded coalgebras
$f_\alpha, f_{\alpha'}:C\to T(V)$ and a unique
$(f_\alpha, f_{\alpha'})$-coderivation of degree $k$
$\partial_{\beta}:C\to T^c(V)$ such that
$p_1\partial_{\beta}=\beta$, i.e., commutes the following diagram:
$$
\begin{array}{rcl}

V & \stackrel{p_1}{\longleftarrow } & T^c(V)\\
   \beta&   \nwarrow & \uparrow \partial_{\beta}\\
  &           & C.
\end{array}
$$
The coderivation $\partial_\beta$ is defined as

$$
\partial_\beta=\sum_{n=0}^\infty\sum_{k=1}^n (\alpha^{\otimes (k-1)} \otimes\beta\otimes \alpha '^{(n-k)}\Delta^k.
$$

%%%%%%%%%%%%%%%%%%%%%%%%%%%%%%%%%%%%%%%%%%%%%%%%%%%%%%%%%%%%%%%%%%%%%%%%%%%%%%%%%
%%%%%%%%%%%%%%%%%%%%%%%%%%%%%%%%%%%%%%%%%%%%%%%%%%%%%%%%%%%%%%%%%%%%%%%%%%%%%%%%%
%%%%%%%%%%%%%%%%%%%%%%%%%%%%%%%%%%%%%%%%%%%%%%%%%%%%%%%%%%%%%%%%%%%%%%%%%%%%%%%%%
%%%%%%%%%%%%%%%%%%%%%%%%%%%%%%%%%%%%%%%%%%%%%%%%%%%%%%%%%%%%%%%%%%%%%%%%%%%%%%%%%
%%%%%%%%%%%%%%%%%%%%%%%%%%%%%%%%%%%%%%%%%%%%%%%%%%%%%%%%%%%%%%%%%%%%%%%%%%%%%%%%%
%%%%%%%%%%%%%%%%%%%%%%%%%%%%%%%%%%%%%%%%%%%%%%%%%%%%%%%%%%%%%%%%%%%%%%%%%%%%%%%%%
%%%%%%%%%%%%%%%%%%%%%%%%%%%%%%%%%%%%%%%%%%%%%%%%%%%%%%%%%%%%%%%%%%%%%%%%%%%%%%%%%
%%%%%%%%%%%%%%%%%%%%%%%%%%%%%%%%%%%%%%%%%%%%%%%%%%%%%%%%%%%%%%%%%%%%%%%%%%%%%%%%%
%%%%%%%%%%%%%%%%%%%%%%%%%%%%%%%%%%%%%%%%%%%%%%%%%%%%%%%%%%%%%%%%%%%%%%%%%%%%%%%%%

%%%%%%%%%%%%%%%%%%%%%%%%%%%%%%%%%%%%%%%%%%%%%%%%%%%%%%%%
\subsection{ Shuffle Comultiplication and Shuffle Comultiplication --- Bialgebra Structures on $T(V)$ and $T'(V)$}  \label{Shuffle}

In fact, $T(V)$ and $T^c(V)$ coincide as graded modules, but the
multiplication of $T(V)$ and the comultiplication of $T^c(V)$ are
not compatible with each other, so they do not define a graded
bialgebra structure on the graded module  $T(V)=T^c(V)$.

Nevertheless there exists the {\it shuffle} comultiplication
$$
\nabla_{sh}:T(V)\to T(V)\otimes T(V)
$$
introduced by Eilenberg
and MacLane \cite{EMacL53}, which turns $(T(V),\Delta_{sh},\mu)$ into a graded
bialgebra.

And dually, there exists the {\it shuffle} multiplication
$\mu_{sh}:T^c(V)\otimes T^c(V)\to T^c(V)$ which turns $(T^c(V),\Delta,\mu_{sh})$ into a graded
bialgebra.

%%%%%%%%%%%%%%%%%%%%%%%%%%%%%%%%%%%%%%%%%%%%%%%%%%%%%%%%

\subsubsection{ Shuffle Comultiplication -- Bialgebra Structure on Tensor Algebra}  \label{ShuffleTV}

There exists on this free graded algebra $(T(V),\mu )$ a comultiplication
$\nabla_{sh}:T(V)\to T(V)\otimes T(V)$ which is a morphism of graded algebras, and, consequently
turns $(T(V), \mu, \nabla_{sh})$ into a graded bialgebra.

Namely, the \emph{shuffle comultiplication}  $\nabla_{sh}:T(V)\to T(V)\otimes T(V)$ is
a graded coalgebra map induced
by the universal property of $T(V)$ (\ref{UnivTV})  by $\alpha:T(V)\to T(V)\otimes
T(V)$ given by
$\alpha(v)=v\otimes 1+1\otimes v,\ \alpha(1)=1\otimes 1$.

This comultiplication is associative and
in fact is given by
$$\begin{array}{c}

\nabla(v_1\otimes\ ...\ \otimes v_n)=1\otimes (v_1\otimes\ ...\ \otimes v_n)+(v_1\otimes\ ...\ \otimes v_n)\otimes 1+\\

\Sigma_{p}\Sigma_{\sigma \in sh(p,n-p)} (v_{\sigma(1)}\otimes\ ...\ \otimes v_{\sigma(p)})\otimes (v_{\sigma(p+1)}\otimes\ ...\ \otimes v_{\sigma(n)})
\end{array}
$$
where $sh(p,n-p)$ consists of all $(p, n-p)$-shuffles, that is all permutations of ${1,\ ...\ ,n}$
such that $\sigma(1)<\  ... \ < \sigma(p)$ and $\sigma(p+ 1) < \ ...\ < \sigma( n)$.

%%%%%%%%%%%%%%%%%%%%%%%%%%%%%%%%%%%%%%%%%%%%%%%%%%%%%%%%

 \subsubsection{ Shuffle Multiplication -- Bialgebra Stucture on Tensor Coalgebra} \label{ShuffleTcV}

%%%%%%%%%%%%%%%%%%%%%%%%%%%%%%%%%%%%%%%%%%%%%%%%%%%%%%%%

The {\it shuffle} multiplication $\mu_{sh}:T^c(V)\otimes T^c(V)\to
T^c(V)$, introduced by Eilenberg and MacLane \cite{EM53}, turns
$(T^c(V),\Delta,\mu_{sh})$ into a graded bialgebra.

This multiplication is defined as the graded coalgebra map induced
by the universal property of $T^c(V)$ (\ref{UnivTcV}) by $\alpha:T^c(V)\otimes
T^c(V)\to V$ given by $\alpha(v\otimes 1)=\alpha(1\otimes v)=v$
and $\alpha=0$ otherwise. This multiplication is associative and
in fact is given by
\begin{equation}
\label{shuffle}
 \mu_{sh} ((v_1\otimes ... \otimes v_n)\otimes
(v_{n+1}\otimes ... \otimes v_{n+m}))= \sum \pm v_{\sigma
(1)}\otimes ... \otimes v_{\sigma (n+m)},
\end{equation}
where summation is taken over all $(m,n)$-shuffles, that is, over
all permutations of the set $(1,2,...,n+m)$ which satisfy the
condition: $i<j$ if $1\leq \sigma (i)<\sigma (j)\leq n$ or
$n+1\leq \sigma (i)<\sigma (j)\leq n+m$.

Particularly,
$$
\begin{array}{l}
\mu_{sh}((v_1)\otimes (v_2))=v_1\otimes v_2 \pm v_2\otimes v_1,\\

\mu_{sh}((v_1)\otimes (v_2\otimes v_3)=v_1\otimes v_2 \otimes v_3 \pm v_2\otimes v_1\otimes v_3 \pm v_2\otimes v_3\otimes v_1,\\
...
\end{array}
$$

%%%%%%%%%%%%%%%%%%%%%%%%%%%%%%%%%%%%%%%%%%%%%%%%%%%%
%%%%%%%%%%%%%%%%%%%%%%%%%%%%%%%%%%%%%%%%%%%%%%%%%%%%
%%%%%%%%%%%%%%%%%%%%%%%%%%%%%%%%%%%%%%%%%%%%%%%%%%%%
%%%%%%%%%%%%%%%%%%%%%%%%%%%%%%%%%%%%%%%%%%%%%%%%%%%%
%%%%%%%%%%%%%%%%%%%%%%%%%%%%%%%%%%%%%%%%%%%%%%%%%%%%

\subsection{Bar and Cobar Constructions}

\subsubsection{Cobar Construction}  \label{Cobar}

Let $(C,d,\Delta)$ be a dg coalgebra with $C_i=0,\ i\leq 1$ and
let $s^{-1}C$ be the desuspension of $C$, that is,
$(s^{-1}C)_k=C_{k+1}$.

The cobar construction $\Omega C$ is defined as the tensor algebra
$T(s^{-1}C)$. We use the following notation for elements of this
tensor coalgebra:
$$
s^{-1}a_1\otimes ...\otimes s^{-1}a_n=[a_1,...,a_n].
$$
So the dimension of $[a_1,...,a_n]$ is $\sum_i\ dim\ a_i - n$.

The differential $d_{\Omega }:\Omega C\to \Omega C$ is defined as
$$
\begin{array}{c}
d_{\Omega }[a_1,...,a_n]=\\ \sum_i\pm
[a_1,...,a_{i-1},da_i,a_{i+1},...,a_n]+\sum_i\pm
[a_1,...,a_{i-1},\Delta (a_i), a_{i+1},...,a_n].
\end{array}
$$
In fact $d_{\Omega }=\partial_\beta$ where $\partial_\beta$ is the derivation
defined by the above universal property (\ref{UnivTVder})
for
$$
\beta [a]=da+\Delta a.
$$
Besides, the properties of $d$ and $\Delta$ from the definition of a dg coalgebra (\ref{DGC}) guarantee that the restriction
$d_{\Omega }d_{\Omega }|V=d_{\Omega }d_{\Omega }i_1$ is $0$ and this, by the universal property (\ref{UnivTVder}), implies $d_{\Omega }d_{\Omega }=0$.
Thus $\Omega C\in DGAlg$.

%%%%%%%%%%%%%%%%%%%%%%%%%%%%

\subsubsection{Bar Construction} \label{Bar}

Let $(A,d, \mu=\cdot)$ be a dg algebra with $A_i=0$ for $i\leq 1$ and let
$sA$ be the suspension of $A$, that is $(sA)_k=A_{k-1}$.

As a graded coalgebra the bar construction $BA$ is defined as the tensor coalgebra:
$T^c(sA)$. We use the following notation for elements of this
tensor coalgebra
$$
sa_1\otimes ...\otimes sa_n=[a_1,...,a_n].
$$
So the dimension of $[a_1,...,a_n]$ is $\sum_i\ dim\ a_i + n$.

The differential $d_B:BA\to BA$ is defined as
$$
d_B[a_1,...,a_n]=\sum_i\pm
[a_1,...,a_{i-1},da_i,a_{i+1},...,a_n]+\sum_i\pm [a_1,...,a_i\cdot
a_{i+1},...,a_n].
$$

In fact, $d_B=D_\beta$ where $D_\beta$ is the coderivation
defined by the above universal property (\ref{UnivTcVder})
for
$$
\beta [a_1,...,a_n]=\left\{
\begin{array}{crr} [\ da_1]& for &
n=1;\\
\ [\ a_1 \cdot a_2\ ] & for & n=2\\
0 &   for & n>2.
\end{array}\right.
$$

Besides, the properties of $d$ and $\mu$ from the definition of a dg algebra (\ref{DGA}) guarantee that the projection
$p_1d_Bd_B$ is $0$ and this, by the universal property (\ref{UnivTcVder}), implies $d_Bd_B=0$. Thus $BA\in DGCoalg$.

\subsubsection{Bar Construction of a Commutative dg Algebra} \label{BarCommut}

Assume now that $(A,d, \mu=\cdot)$ is a commutative dg algebra.
How the commutativity reflects on the bar construction $BA$?

By definition, the differential $d_B:BA\to BA$ is a coderivatian with respect
to the standard graded coalgebra structure of the tensor coalgebra  $BA=T^c(sA)$.
As it was mentioned above in (\ref{ShuffleTcV}), $BA$ carries also the shuffle product
$\mu_{sh}:T^c(sV)\otimes T^c(sA)\to T^c(sA)$ which turns it into a graded bialgebra.

\begin{proposition}.
If a dg algebra $(A,d, \mu=\cdot)$ is commutative, the differential $d_\beta:BA\to BA$
is not only a coderivation with respect to comultiplication $\Delta$ but also a derivation with respect
to the shuffle product $\mu_{sh}$. So in this case the bar construction $(BA,\Delta, \mu_{sh})$ is a dg bialgebra.
\end{proposition}

\noindent \textbf{Proof.}
The map $ \Phi:BA\otimes  BA\to BA$ defined as
$$
\Phi=d_\beta\mu_{sh}-\mu_{sh}(d_\beta\otimes id+id\otimes d_\beta)
$$
is a coderivation, see the arguments in (\ref{Gcoalg}).
Thus, according to universal the property (\ref{UnivTcVder})
of $T^c(s^{-1}A)$, the map $\Phi$ is trivial if and only if
$p_1\Phi=0$. This is equivalent to the commutativity of $A$.

\subsubsection{Adjunction}  \label{AdjBarCobar}

Let $(C,d,\Delta)$ be a dg coalgebra and  $(A,d,\mu)$ a dg algebra. A twisting cochain
\cite{Brown} is a homomorphism $\tau:C\to A$ of degree +1
satisfying Browns' condition
\begin{equation}
 d\tau+\tau d=\tau\smile \tau,
\end{equation}
where $\tau\smile \tau'=\mu_A(\tau\otimes\tau')\Delta $. We denote
by $T(C,A)$ the set of all twisting cochains $\tau:C\to A$.

There are universal twisting cochains $C\to \Omega C$ and $BA\to
A$, namely, the obvious inclusion and projection, respectively. Here are
essential consequences of the condition (\ref{Brown}):

\noindent (i) The multiplicative extension $f_{\tau}:\Omega C\to
A$ is a map of dg algebras, so there is a bijection
$T(C,A)\leftrightarrow Hom_{DG-Alg}(\Omega C,A)$;

\noindent (ii) The comultiplicative coextension $f_{\tau}:C\to BA$
is a map of dg coalgebras, so there is a bijection
$T(C,A)\leftrightarrow Hom_{DG-Coalg}(C,BA)$.

Thus we have two bijections
$$
Hom_{DG-Alg}(\Omega C,A) \leftrightarrow T(C,A)\leftrightarrow
Hom_{DG-Coalg}(C,BA).
$$

Besides, there are two weak equivalences (homology isomorphisms)
$$
\alpha_A:\Omega B(A)\to A,\ \ \ \ \ \  \beta_C:C\to B\Omega (C).
$$

\section{Twisting Cochains and Functor $D$} \label{Twist}

\subsection{Brown's twisting cochains}

\subsubsection{Definition of twisting cochain}

Let $(K,d_K:K_*\to K_{*-1},\nabla_K:K\to K\otimes K)$ be a dg coalgebra and let $(A,d_A:A_*\to A_{*-1},\mu_A:A\otimes A\to A)$ be a dg algebra. Then
$C^*(K,A)=Hom(K,A)$ with differential $\delta \alpha=\alpha d_K+d_A\alpha$ and multiplication $\alpha \smile \beta=\mu_A(\alpha\otimes \beta)\nabla_K$ is a dg algebra.

\begin{definition}.
A Brown twisting cochain is a homomorphism
$$
\phi:K_*\to A_{*-1},
$$
i.e., $deg\phi=-1$, satisfying
\begin{equation}
\label{Brown} d\phi+\phi d=\phi\smile \phi.
\end{equation}
\end{definition}
Brown's condition sometimes is called also \emph{Maurer-Cartan equation}. It can be rewritten as
$d_A\phi+\phi d_K=\mu_A(\phi\otimes \phi)\nabla_K$.

The set of all twisting cochains $\phi:K\to A$ we denote as $Tw(K,A)$.

\subsubsection{Twisted tensor product}

Let $(M,d_M,\nu:A\otimes M\to M)$ be a dg $A$-module. Then any twisting cochain $\phi:K\to A$ determines a homomorphism
$d_\phi:K\otimes M\to K\otimes M$ by
$$
d_\phi(k\otimes m)=d_Kk\otimes m+k\otimes d_Mm+(k\otimes m)\cap \phi
$$
where $(k\otimes m)\cap \phi=(id_K\otimes \nu)(id_K\otimes \phi \otimes id_M)(\nabla_K\otimes id_M)(k\otimes m)$.
Brown's condition $d\phi =\phi\smile \phi$ implies that $d_\phi d_\phi=0$.

\begin{definition}.
The twisted tensor product $(K\otimes_\phi M, d_\phi)$ is defined as the chain complex $(K\otimes M,\ d_\phi)$.
\end{definition}

\subsubsection{Application: model of fibration}

Let $(E,p,B,F,G)$ be a fibre bundle with base $B$, fibre $F$, structure group $G$. So $K=C_*(B)$ is a dg coalgebra, $A=C_*(G)$ is a dg algebra,
$M=C_*(F)$ is a dg module over $A$.
Then, by Brown's theorem there exists a twisting cochain $\phi:K=C_*(B)\to A=C_*(G)$ such that the twisted tensor product
 $$
 (K\otimes_\phi M, d_\phi)=(C_*(B)\otimes_{\phi:C_*(B)\to C_{*-1}(G)} C_*(F), d_\phi)
 $$
gives the homology of fibre space $H_*(E)$. It is clear that $\phi$ determines all differentials of the Serre spectral sequence also.

Brown's twisting cochain $\phi$ is not determined uniquely, so it would be useful to have a possibility to choose one convenient for computations.

%%%%%%%%%%%%%%%%%%%%%%%%%%%%%%%%%%%%%%%%%%%%%%%%%%%%%%%%%%%%%%%%%%%%%%%%%%%%%%%%%%%%%%%%%%%%%%%%%%%%%%%%%%%%%%%%%%%%%%%%%%%%%%%

\subsubsection{Twisting cochains and the Bar and cobar constructions}

By the universal property of the tensor coalgebra from (B1.3) any homomorphism $\phi:K\to A$ of degree $-1$ induces a graded coalgebra map $f_\phi:K\to T^c(sA)$ by
$$
f_\phi=\sum_i (\phi\otimes \ ...\ \otimes \phi)\nabla_K^i.
$$
If, in addition, $\phi:K\to A$ is a twisting cochain, that is, it satisfies  Brown's condition $\delta\phi=\phi\smile \phi$, then $f_\phi:K\to B(A)$ is a chain map, i.e., is a map of dg coalgebras: since the tensor coalgebra is cofree one has the equality $d_Bf_\phi=f_\phi d_K$ if and only if equal the projections $pd_Bf_\phi$ and $=pf_\phi d_K$ are equal, and this is exactly Brown's condition $\delta\phi=\phi\smile \phi$.

Conversely, any dg coalgebra map  $f:K\to BA$ is $f_\phi$ for $\phi=p f:K\to BA\to A$.
In fact we have a bijection $Mor_{DGCoalg}(K,BA)\leftrightarrow Tw(K,A)$.

%%%%%%%%%%%%%%%%%%%%%%%%%%%%%%%%%%%%%%%%%%%%%%%%%%%%%%%%%%%%%%%%%%%%%%%%%%%%%%%%%%%%%%%%%%%%%%%%%%%%%%%%%%%%%%%%%%%%%%%%%%%%%%%%%%%%%%%%%%%%%%%

Dually,
by the universal property of the tensor algebra (B1.2), any homomorphism $\phi:K\to A$ of degree $-1$ induces a graded algebra map $g_\phi:T(s^{-1}K)\to A$ by
$$
g_\phi =\sum_k \mu^k(\alpha\otimes...\otimes \alpha).
$$
If, in addition, $\phi:K\to A$ is a twisting cochain, that is, it satisfies  Brown's condition $\delta\phi=\phi\smile \phi$, then $g_\phi:\Omega K\to A$ is a chain map, i.e., it is a map of dg algebras: since of freeness of tensor algebra one has the equality $d_Ag_\phi=g_\phi d_\Omega$
if and only if equal the restrictions $d_Ag_\phi i=g_\phi d_\Omega i$ and this is exactly the Brown's condition $\delta\phi=\phi\smile \phi$.

Conversely, any dg algebra map  $g:\Omega K\to A$ is $g_\phi$ for $\phi=gi:K\to \Omega K\to A$.
In fact we have a bijection $Mor_{DGAlg}(\Omega K,A)\leftrightarrow Tw(K,A)$

%%%%%%%%%%%%%%%%%%%%%%%%%%%%%%%%%%%%%%%%%%%%%%%%%%%%%%%%%%%%%%%%%%%%%%%%%%%%%%%%%%%%%%%%%%%%%%%%%%%%%%%%%%%%%%%%%%%%%%%%%%%%%%%%%%%%%%%%%%%%%%%

So we have bijections
$$
\begin{array}{ccccc}
& & Tw(K,A) & &  \\

 & \approx\swarrow & & \searrow \approx& \\

Mor_{DGAlg}(\Omega K,A) &   &   &   &  Mor_{DGCoalg}(K,BA)
\end{array}
$$

%%%%%%%%%%%%%%%%%%%%%%%%%%%%%%%%%%%%%%%%%%%%%%%%%%%%%%%%%%%%%%%%%%%%%%%%%%%%%%%%%%%%%%%%%%%%%%%%%%%%%%%%%%%%%%%%%%%%%%%%%%%%
%%%%%%%%%%%%%%%%%%%%%%%%%%%%%%%%%%%%%%%%%%%%%%%%%%%%%%%%%%%%%%%%%%%%%%%%%%%%%%%%%%%%%%%%%%%%%%%%%%%%%%%%%%%%%%%%%%%%%%%%%%%%
%%%%%%%%%%%%%%%%%%%%%%%%%%%%%%%%%%%%%%%%%%%%%%%%%%%%%%%%%%%%%%%%%%%%%%%%%%%%%%%%%%%%%%%%%%%%%%%%%%%%%%%%%%%%%%%%%%%%%%%%%%%%
%%%%%%%%%%%%%%%%%%%%%%%%%%%%%%%%%%%%%%%%%%%%%%%%%%%%%%%%%%%%%%%%%%%%%%%%%%%%%%%%%%%%%%%%%%%%%%%%%%%%%%%%%%%%%%%%%%%%%%%%%%%%
%%%%%%%%%%%%%%%%%%%%%%%%%%%%%%%%%%%%%%%%%%%%%%%%%%%%%%%%%%%%%%%%%%%%%%%%%%%%%%%%%%%%%%%%%%%%%%%%%%%%%%%%%%%%%%%%%%%%%%%%%%%%

\subsubsection{(Co)universal twisting cochains}

The standard inclusion $i:K\to \Omega K$ satisfies Brown's condition. So it is a (so called universal) twisting cochain. Thus it induces two chain maps
$$
f_i:\Omega K\to \Omega K\ \ \ and\ \ \ g_i:K\to B\Omega K.
$$
The first one is the identity map and the second one a quasi isomorphism (homology isomorphism).

Dually, the standard projection $p:BA\to A $ satisfies the Brown's condition, so it is a (so called couniversal) twisting cochain.
 Thus it induces two chain maps
$$
f_p:\Omega BA\to A\ \ \ and\ \ \ g_p:BA\to BA.
$$
The first one is a quasi isomorphism (homology isomorphism) and the second one is the identity map.

%%%%%%%%%%%%%%%%%%%%%%%%%%%%%%%%%%%%%%%%%%%%%%%%%%%%%%%%%%%%%%%%%%%%%%%%%%%%%%%%%%%%%%%%%%%%%%%%%%%%%%%%%%%%%%%%%%%%%%%%%%%%
%%%%%%%%%%%%%%%%%%%%%%%%%%%%%%%%%%%%%%%%%%%%%%%%%%%%%%%%%%%%%%%%%%%%%%%%%%%%%%%%%%%%%%%%%%%%%%%%%%%%%%%%%%%%%%%%%%%%%%%%%%%%
%%%%%%%%%%%%%%%%%%%%%%%%%%%%%%%%%%%%%%%%%%%%%%%%%%%%%%%%%%%%%%%%%%%%%%%%%%%%%%%%%%%%%%%%%%%%%%%%%%%%%%%%%%%%%%%%%%%%%%%%%%%%
%%%%%%%%%%%%%%%%%%%%%%%%%%%%%%%%%%%%%%%%%%%%%%%%%%%%%%%%%%%%%%%%%%%%%%%%%%%%%%%%%%%%%%%%%%%%%%%%%%%%%%%%%%%%%%%%%%%%%%%%%%%%
%%%%%%%%%%%%%%%%%%%%%%%%%%%%%%%%%%%%%%%%%%%%%%%%%%%%%%%%%%%%%%%%%%%%%%%%%%%%%%%%%%%%%%%%%%%%%%%%%%%%%%%%%%%%%%%%%%%%%%%%%%%%

\subsection{Berikashvili's Functor $D$}

\subsubsection{Equivalence of twisting cochains}

Two twisting cochains $\phi,\psi:K\to A$ are equivalent (Berikashvili \cite{Berik1}) if there exists
$c:K\to A,\ deg\ c=0,\ c=\Sigma_i c_i,\ c_i:C_i\to A_i$  with $c_0=c|C_0=0$, such that
 \begin{equation}
\label{Berikeq}
\psi=\phi+\delta c+\psi\smile c+c\smile \phi,
\end{equation}
notation $\phi\sim_c \psi$.

This is an equivalence relation:
$$
\phi\sim_{c=0}\phi;\ \ \ \phi\sim_c\phi',\phi'\sim_{c'}\phi'' \Rightarrow \phi\sim_{c+c'+c'\smile c}\phi''
$$
and $\phi\sim_c\phi'\Rightarrow \phi'\sim_{c'}\phi$ where $c'$ can be solved from $c+c'+c'\smile c=0$ inductively.

%%%%%%%%%%%%%%%%%%%%%%%%%%%%%%%%%%%%%%%%%%%%%%%%%%%%%%%%%%%%%%%%%%%%%%%%%%%%%%%%%%%%%%%%%%%%%%%%%%%
%%%%%%%%%%%%%%%%%%%%%%%%%%%%%%%%%%%%%%%%%%%%%%%%%%%%%%%%%%%%%%%%%%%%%%%%%%%%%%%%%%%%%%%%%%%%%%%%%%%
%%%%%%%%%%%%%%%%%%%%%%%%%%%%%%%%%%%%%%%%%%%%%%%%%%%%%%%%%%%%%%%%%%%%%%%%%%%%%%%%%%%%%%%%%%%%%%%%%%%
%%%%%%%%%%%%%%%%%%%%%%%%%%%%%%%%%%%%%%%%%%%%%%%%%%%%%%%%%%%%%%%%%%%%%%%%%%%%%%%%%%%%%%%%%%%%%%%%%%%
%%%%%%%%%%%%%%%%%%%%%%%%%%%%%%%%%%%%%%%%%%%%%%%%%%%%%%%%%%%%%%%%%%%%%%%%%%%%%%%%%%%%%%%%%%%%%%%%%%%

This notion of equivalence allows to \emph{perturb} twisting cochains.
Let
$$
\phi=\phi_2+\phi_3+\ ...\ +\phi_n+\ ...\ :K\to A,\ \ \ \phi_n:K_n\to A_{n-1}
$$
be a twisting cochain, and let's take an arbitrary cochain $c=c_n:K_n\to A_n$. Then there exists a twisting cochain $F_{c_n}\phi=\psi:K\to A$ such that $\phi\sim_{c_n}\psi$. Actually the components of the perturbed twisting cochain
$$
F_{c_n}\phi=\psi=\psi_2+\psi_3+\ ...\ +\psi_n+\ ...
$$
can be solved from (\ref{Berikeq}) inductively, and the solution particularly gives that the perturbation $F_{c_n}\phi$ does not change the first components, i.e.
 $\psi_i=\phi_i$ for $i<n$ and $\psi_n=\phi_n+d_Ac_n$.

The main benefit of this notion is that the equivalent twisting cochains define isomorphic twisted tensor products:

\begin{theorem}. If $\phi\sim_c \psi$ then
$$
(K\otimes_\phi M, d_\phi) \stackrel{F_c}{\longrightarrow}  (K\otimes_\psi M, d_\psi)
$$
given by $F_c(k\otimes m)=(k\otimes m)+(k\otimes m)\cap c$   is an isomorphism of dg modules.
\end{theorem}

The inverse $(K\otimes_\psi M, d_\psi) \stackrel{F_{c'}}{\longrightarrow}  (K\otimes_\phi M, d_\phi)$
is defined by $c':K\to A$ which, as above, can be solved inductively from $c+c'+c'\smile c=0$.

Let
$$
Tw(K,A)=\{\phi:K\to A,\ \delta\phi=\phi\circ \phi\}.
$$
be the set of all twisting cochains.

\begin{definition}.
Berikasvili's functor $D(K,A)$ is defined as the factorset $D(K,A)=\frac{Tw(K,A)}{\sim}$.
\end{definition}

Berikashvili's relation of equivalence of twisting cochains allows to perturb a given twisting cochain so as to get the simplest
one to simplify calculations.

\subsubsection{Equivalence of twisting cochains and homotopy of induced maps}

How does the equivalence of twistimg cochains affect induced maps $K\to BA$ and $\Omega K\to A$?

%%%%%%%%%%%%%%%%%%%%%%%%%%%%%%%%%%%%%%%%%%%%%%%%%%%%%%%%%%%%%%%%%%%%%%%%%%%%%%%%%%%%%%%%%%%%%%%%%%%
%%%%%%%%%%%%%%%%%%%%%%%%%%%%%%%%%%%%%%%%%%%%%%%%%%%%%%%%%%%%%%%%%%%%%%%%%%%%%%%%%%%%%%%%%%%%%%%%%%%
%%%%%%%%%%%%%%%%%%%%%%%%%%%%%%%%%%%%%%%%%%%%%%%%%%%%%%%%%%%%%%%%%%%%%%%%%%%%%%%%%%%%%%%%%%%%%%%%%%%
%%%%%%%%%%%%%%%%%%%%%%%%%%%%%%%%%%%%%%%%%%%%%%%%%%%%%%%%%%%%%%%%%%%%%%%%%%%%%%%%%%%%%%%%%%%%%%%%%%%
%%%%%%%%%%%%%%%%%%%%%%%%%%%%%%%%%%%%%%%%%%%%%%%%%%%%%%%%%%%%%%%%%%%%%%%%%%%%%%%%%%%%%%%%%%%%%%%%%%%

%%%%%%%%%%%%%%%%%%%%%%%%%%%%%%%%%%%%%%%%%%%%%%%%%%%%%%%%%%%%%%%%%%%%%%%%%%%%%%%%%%%%%%%%%%%%%%%%%%%
%%%%%%%%%%%%%%%%%%%%%%%%%%%%%%%%%%%%%%%%%%%%%%%%%%%%%%%%%%%%%%%%%%%%%%%%%%%%%%%%%%%%%%%%%%%%%%%%%%%
%%%%%%%%%%%%%%%%%%%%%%%%%%%%%%%%%%%%%%%%%%%%%%%%%%%%%%%%%%%%%%%%%%%%%%%%%%%%%%%%%%%%%%%%%%%%%%%%%%%
%%%%%%%%%%%%%%%%%%%%%%%%%%%%%%%%%%%%%%%%%%%%%%%%%%%%%%%%%%%%%%%%%%%%%%%%%%%%%%%%%%%%%%%%%%%%%%%%%%%
%%%%%%%%%%%%%%%%%%%%%%%%%%%%%%%%%%%%%%%%%%%%%%%%%%%%%%%%%%%%%%%%%%%%%%%%%%%%%%%%%%%%%%%%%%%%%%%%%%%

\begin{theorem}.
If $\phi\sim_c\psi$ then $f_\phi$ and $f_\psi$ are homotopic as dg coalgebra maps: the coderivation chain homotopy $D(c):K\to BA$ is given by
\begin{equation}\label{Dc}
D(c)=\sum_{i,j} (\psi\otimes\ ...\ (j\ \times)\ ...\ \otimes \psi\otimes c\otimes\phi\otimes\ ...\ \otimes\phi)\nabla_K^i.
\end{equation}
\end{theorem}
\textbf{Proof.} The homotopy $D(c)$ satisfies $f_\phi-f_\psi=d_BD(c)+D(c)d_K$ because the projection of this condition on $A$ gives $pf_\phi-pf_\psi=pd_BD(c)+pD(c)d_K$ and this is exactly
the condition (\ref{Berikeq}). Besides, $D(c)$ is a $f_\phi-f_\psi$-coderivation, that is,
$$\nabla_{B}D(c)=(f_\psi\otimes D(c)+D(c)\otimes f_\phi)\nabla_K.$$
This also follows from the universal property of the tensor coalgebra and the extension rule (\ref{Dc}).

The converse is also true: if $f_\phi$ and $f_\psi$ are homotopic by a coderivation homotopy $D:K\to BA$ then $\phi \sim_c\psi$ by $c=p_1D$.

The dual statement is also true: $g_\phi,\ g_\psi:\omega K\to A$ are homotopic in the category of dg algebras if and only if $\phi\sim \psi$.

So we have bijections
$$
\begin{array}{ccccc}
& & D(K,A) & &  \\

 & \approx\swarrow & & \searrow \approx& \\

[\Omega K,A]_{DGAlg} &   &   &   &  [(K,BA)]_{DGCoalg}
\end{array}
$$
this means that $B$ and $\Omega$ are adjoint functors.
%%%%%%%%%%%%%%%%%%%%%%%%%%%%%%%%%%%%%%%%%%%%%%%%%%%%%%%%%%%%%%%%%%%%%%%%%%%%%%%%%%%%%%%%%%%%%%%%%%%
%%%%%%%%%%%%%%%%%%%%%%%%%%%%%%%%%%%%%%%%%%%%%%%%%%%%%%%%%%%%%%%%%%%%%%%%%%%%%%%%%%%%%%%%%%%%%%%%%%%
%%%%%%%%%%%%%%%%%%%%%%%%%%%%%%%%%%%%%%%%%%%%%%%%%%%%%%%%%%%%%%%%%%%%%%%%%%%%%%%%%%%%%%%%%%%%%%%%%%%
%%%%%%%%%%%%%%%%%%%%%%%%%%%%%%%%%%%%%%%%%%%%%%%%%%%%%%%%%%%%%%%%%%%%%%%%%%%%%%%%%%%%%%%%%%%%%%%%%%%
%%%%%%%%%%%%%%%%%%%%%%%%%%%%%%%%%%%%%%%%%%%%%%%%%%%%%%%%%%%%%%%%%%%%%%%%%%%%%%%%%%%%%%%%%%%%%%%%%%%

%%%%%%%%%%%%%%%%%%%%%%%%%%%%%%%%%%%%%%%%%%%%%%%%%%%%%%%%%%%%%%%%%%%%%%%%%%%%%%%%%%%%%%%%%%%%%%%%%%%
%%%%%%%%%%%%%%%%%%%%%%%%%%%%%%%%%%%%%%%%%%%%%%%%%%%%%%%%%%%%%%%%%%%%%%%%%%%%%%%%%%%%%%%%%%%%%%%%%%%
%%%%%%%%%%%%%%%%%%%%%%%%%%%%%%%%%%%%%%%%%%%%%%%%%%%%%%%%%%%%%%%%%%%%%%%%%%%%%%%%%%%%%%%%%%%%%%%%%%%
%%%%%%%%%%%%%%%%%%%%%%%%%%%%%%%%%%%%%%%%%%%%%%%%%%%%%%%%%%%%%%%%%%%%%%%%%%%%%%%%%%%%%%%%%%%%%%%%%%%
%%%%%%%%%%%%%%%%%%%%%%%%%%%%%%%%%%%%%%%%%%%%%%%%%%%%%%%%%%%%%%%%%%%%%%%%%%%%%%%%%%%%%%%%%%%%%%%%%%%

%%%%%%%%%%%%%%%%%%%%%%%%%%%%%%%%%%%%%%%%%%%%%%%%%%%%%%%%%%%%%%%%%%%%%%%%%%%%%%%%%%%%%%%%%%%%%%%%%%%
%%%%%%%%%%%%%%%%%%%%%%%%%%%%%%%%%%%%%%%%%%%%%%%%%%%%%%%%%%%%%%%%%%%%%%%%%%%%%%%%%%%%%%%%%%%%%%%%%%%
%%%%%%%%%%%%%%%%%%%%%%%%%%%%%%%%%%%%%%%%%%%%%%%%%%%%%%%%%%%%%%%%%%%%%%%%%%%%%%%%%%%%%%%%%%%%%%%%%%%
%%%%%%%%%%%%%%%%%%%%%%%%%%%%%%%%%%%%%%%%%%%%%%%%%%%%%%%%%%%%%%%%%%%%%%%%%%%%%%%%%%%%%%%%%%%%%%%%%%%
%%%%%%%%%%%%%%%%%%%%%%%%%%%%%%%%%%%%%%%%%%%%%%%%%%%%%%%%%%%%%%%%%%%%%%%%%%%%%%%%%%%%%%%%%%%%%%%%%%%

\subsubsection{Lifting of Twisting Cochains}

Any dg algebra map $f:A\to A'$ induces the map $Tw(K,A)\to Tw(K,A')$: if $\phi:K\to M$ is a twisting cochain so is the composition $f \phi:K\to A\to A'$.
Moreover, if $\phi  \sim_c\phi'$ then $f \phi\sim_{f c}f \phi'$.
Thus we have a map $D(f):D(K,A)\to D(K,A')$.

\begin{theorem}.\label{BerThm} {(Berikashvili \cite{Berik1}). Let $(K,d_K,\nabla_K)$ be a dg colagebra with free $K_i$ and let $(A,d_A,\mu_A)$ be a connected dg algebra.
If $f:A\to A'$ is a weak equivalence of connected dg algebras (i.e., a homology isomorphism), then
$$
D(f):D(K,A)\to D(K,A')
$$
is a bijection.}
\end{theorem}

Particularly this theorem means that $[K,BA]\to [K,BA']$ is a bijection.

%%%%%%%%%%%%%%%%%%%%%%%%%%%%%%%%%%%%%%%%%%%%%%%%%%%%%%%%%%%%%%%%%%%%%%%%%%%%%%%%%%%%%%%%%%%%%%%%%%%
%%%%%%%%%%%%%%%%%%%%%%%%%%%%%%%%%%%%%%%%%%%%%%%%%%%%%%%%%%%%%%%%%%%%%%%%%%%%%%%%%%%%%%%%%%%%%%%%%%%
%%%%%%%%%%%%%%%%%%%%%%%%%%%%%%%%%%%%%%%%%%%%%%%%%%%%%%%%%%%%%%%%%%%%%%%%%%%%%%%%%%%%%%%%%%%%%%%%%%%
%%%%%%%%%%%%%%%%%%%%%%%%%%%%%%%%%%%%%%%%%%%%%%%%%%%%%%%%%%%%%%%%%%%%%%%%%%%%%%%%%%%%%%%%%%%%%%%%%%%
%%%%%%%%%%%%%%%%%%%%%%%%%%%%%%%%%%%%%%%%%%%%%%%%%%%%%%%%%%%%%%%%%%%%%%%%%%%%%%%%%%%%%%%%%%%%%%%%%%%

%%%%%%%%%%%%%%%%%%%%%%%%%%%%%%%%%%%%%%%%%%%%%%%%%%%%%%%%%%%%%%%%%%%%%%%%%%%%%%%%%%%%%%%%%%%%%%%%%%%
%%%%%%%%%%%%%%%%%%%%%%%%%%%%%%%%%%%%%%%%%%%%%%%%%%%%%%%%%%%%%%%%%%%%%%%%%%%%%%%%%%%%%%%%%%%%%%%%%%%
%%%%%%%%%%%%%%%%%%%%%%%%%%%%%%%%%%%%%%%%%%%%%%%%%%%%%%%%%%%%%%%%%%%%%%%%%%%%%%%%%%%%%%%%%%%%%%%%%%%
%%%%%%%%%%%%%%%%%%%%%%%%%%%%%%%%%%%%%%%%%%%%%%%%%%%%%%%%%%%%%%%%%%%%%%%%%%%%%%%%%%%%%%%%%%%%%%%%%%%
%%%%%%%%%%%%%%%%%%%%%%%%%%%%%%%%%%%%%%%%%%%%%%%%%%%%%%%%%%%%%%%%%%%%%%%%%%%%%%%%%%%%%%%%%%%%%%%%%%%

%%%%%%%%%%%%%%%%%%%%%%%%%%%%%%%%%%%%%%%%%%%%%%%%%%%%%%%%%%%%%%%%%%%%%%%%%%%%%%%%%%%%%%%%%%%%%%%%%%%
%%%%%%%%%%%%%%%%%%%%%%%%%%%%%%%%%%%%%%%%%%%%%%%%%%%%%%%%%%%%%%%%%%%%%%%%%%%%%%%%%%%%%%%%%%%%%%%%%%%
%%%%%%%%%%%%%%%%%%%%%%%%%%%%%%%%%%%%%%%%%%%%%%%%%%%%%%%%%%%%%%%%%%%%%%%%%%%%%%%%%%%%%%%%%%%%%%%%%%%
%%%%%%%%%%%%%%%%%%%%%%%%%%%%%%%%%%%%%%%%%%%%%%%%%%%%%%%%%%%%%%%%%%%%%%%%%%%%%%%%%%%%%%%%%%%%%%%%%%%
%%%%%%%%%%%%%%%%%%%%%%%%%%%%%%%%%%%%%%%%%%%%%%%%%%%%%%%%%%%%%%%%%%%%%%%%%%%%%%%%%%%%%%%%%%%%%%%%%%%

Below we'll need the surjectivity part of this theorem whose proof we sketch here.

\begin{theorem}.\label{kadlifting} Let $(K,d_K,\nabla_K)$ be a dg colagebra with free $K_i$s and let
$f:(A,d_A,\mu_A)\to (A',d_{A'},\mu_{A'})$ be a weak equivalence of connected dg algebras. Then for an arbitrary twisting cochain
$$
\phi=\phi_2+\phi_3+\ ...\ +\phi_n+\ ...\ :K\to A'
$$
there exists a twisting cochain
$$
\psi=\psi_2+\psi_3+\ ...\ +\psi_n+\ ...\ :K\to A
$$
such that $\phi \sim f\psi$.
\end{theorem}

\noindent \textbf{Proof.} Start with a twisting cohain
$\phi=\phi_2+\phi_3+\ ...\ +\phi_n+\ ...\ :K\to A'.$
Brown's defining condition (\ref{Brown}) gives $d_A\phi_2=0$, and since $f:A\to A'$ is a homology isomorphism there exist $\psi_2:K_2\to A_1$ and $c'_2:K_2\to A'_2$ such that $d_A\psi_2=0$ and $f\psi_2=\phi_2+d_{A'}c'_2$ (we assume all $K_n$ are free modules). Perturbing $\phi$ by this $c'_2$ we obtain a new twisting cochain $F_{c'_2}\phi$ for which
$(F_{c'_2}\phi)_2=\phi_2+d_{A'}c'_2=f\psi_2$. So we can assume that $\phi_2=f\psi_2$.

Assume now that we already have
$\psi_2,\ \psi_3,\  ...\ ,\psi_{n-1}$ which satisfy (\ref{Brown}) in appropriate dimensions and $f\psi_k=\phi_k,\ k=2,3,...,n-1$. We need the next component $\psi_n:K_n\to A_{n-1}$ such that
$$
d_A\psi_n=\psi_{n-1}d_K+\sum_{i=2}^{n-2}\psi_i\smile\psi_{n-i}
$$
and $c'_n:K_n\to A'_n$ such that $f\psi_n=\phi_n+d_{A'}c'_n$. Then perturbing $\phi$ by $c'_n$ we obtain new $\phi_n$ for which $f\psi_n=\phi_n$, and this will complete the proof.

Let us write
$$
\begin{array}{c}
U_n=\psi_{n-1}d_K+\sum_{i=2}^{n-2}\psi_i\smile\psi_{n-i}:K_n\to A_{n-2},\\
U'_n=\phi_{n-1}d_K+\sum_{i=2}^{n-2}\phi_i\smile\phi_{n-i}:K_n\to A'_{n-2}.
\end{array}
$$
So we have $U'_n =d_{A'}\phi_n$ and we want $\psi_n:K_n\to A_{n-1},\ \ c'_n:K_n\to A'_n$ such that $U_n=d_A\psi_n$ and $f\psi_n=\phi_n+d_{A'}c'_n$.

First, it is not hard to check that $d_AU_n=0$, that is, $U_n$ maps $K_n$ to cycles $Z(A_{n-2})\subset A_{n-2}$.

Then
$$
\begin{array}{c}
fU_n=f(\psi_{n-1}d_K+\sum_{i=2}^{n-2}\psi_i\smile\psi_{n-i})
=f\psi_{n-1}d_K+\sum_{i=2}^{n-2}f\psi_i\smile f\psi_{n-i}\\
=\phi_{n-1}d_K+\sum_{i=2}^{n-2}\phi_i\smile\phi_{n-i}=U'_n=d_A\phi_n,
\end{array}
$$
thus, $fU_n$ maps $K_n$ to boundaries $B(A_{n-2})\subset A_{n-2}$.
Since $f:A\to A'$ is a homology isomorphism, there is $\overline{\psi}_n:K_n\to A_{n-1}$ such that
$d_A\overline{\psi}_n=U_n.$

Now we must take care of the condition $f\psi_n=\phi_n$. It is clear  that $d_Af\overline{\psi_n}=d_A\phi_n$. Let us denote by
$z'_n=f\overline{\psi_n}-\phi_n$; this is the homomorphism which maps $K_n$ to cycles $Z(A'_{n-1})$.  Again, since $f:A\to A'$ is a homology isomorphism there exist $z_n:K_n\to Z(A_{n-1})$ and $c'_n:K_n\to A'_{n}$ such that $fz_n=z'_n-d_{A'}c'_n$. Let us define $\psi_n=\overline{\psi_n}-z_n$. Then $d_A\psi_n=d_A\overline{\psi_n}$ and
$$
f\psi_n=f\overline{\psi_n}-fz_n=(\phi_n+z'_n)-(z'_n-d_{A'}c'_n)=\phi_n+d_{A'}c'_n.
$$

Perturbing $\phi$ by this $c'_n$ we obtain $F_{c'_n}\phi$ with
$(F_{c'_n}\phi)_n=\phi_n+d_A'c'_n=f\psi_n$.

\section{Stasheff's $A_\infty$-algebras }\label{Ainf}

\subsection{Category of $A_\infty$-algebras}

The notion of $A_{\infty}$-algebra was introduces by J. Stasheff
\cite{Sta63}. This notion generalizes the notion of differential
graded algebra and in fact it is so called \emph{strong homotopy associative}
algebra where the strict associativity is replaced with \emph{associativity up to higher
coherent homotopies}.

\subsubsection{Notion of $A_\infty$-algebra}  \label{Ainf}

\begin{definition}.\label{ainf}
An $A_{\infty}$-algebra is a graded module $M=\{M^k\}_{k\in Z}$
equipped with a sequence of operations
$$
\{m_i:M\otimes ...(i\ times)...\otimes M\to M, i=1,2,3,. . . \}
$$
satisfying the conditions $ m_i((\otimes^iM)^q)\subset M^{q-i+2}
$, that is $deg\ m_i=2-i$, and
\begin{equation}
\label{ainfalg}
\begin{array}{l}
\sum_{k=0}^{i-1} \sum_{j=1}^{i-k} \pm \\
m_{i-j+1}(a_1\otimes . . . \otimes a_k\otimes  m_j(a_{k+1}\otimes
. . . \otimes a_{k+j})\otimes . . . \otimes a_i)=0 .
\end{array}
\end{equation}
\end{definition}

For $i=1$ this condition reads
$$
m_1m_1=0.
$$
For $i=2$ this condition reads
$$
m_1m_2(a_1\otimes a_2)\pm m_2(m_1(a_1)\otimes a_2)\pm
m_2(a_1\otimes m_1(a_2))=0.
$$
For $i=3$ this condition reads
$$
\begin{array}{c}
m_1m_3(a_1\otimes a_2\otimes a_3)\pm \\
m_3(m_1(a_1)\otimes a_2\otimes a_3)\pm m_3(a_1\otimes
m_1(a_2)\otimes a_3)\pm m_3(a_1\otimes a_2\otimes m_1(a_3))\pm \\
m_2(m_2(a_1\otimes a_2)\otimes a_3)\pm m_2(a_1\otimes
m_2(a_2\otimes a_3)=0 .
\end{array}
$$
These three conditions mean that for an $A_{\infty}$-algebra
$(M,\{m_i\})$ the first two operations form a {\it nonassociative} dga
$(M,m_1,m_2)$ with differential $m_1$ and multiplication $m_2$
which is associative just up to homotopy and the suitable homotopy
is the operation $m_3$.

\emph{Spacial case: }An $A_\infty$-algebra $(M,\{m_1,m_2,m_3=0,m_4=0,...\})$ is a strictly associative dg algebra.

%%%%%%%%%%%%%%%%%%%%%%%%%%%%%%%%%%%%%%%%%%%%%%%%%%%%%%%%%%%%%%%%%%%%%%%%%%%
%%%%%%%%%%%%%%%%%%%%%%%%%%%%%%%%%%%%%%%%%%%%%%%%%%%%%%%%%%%%%%%%%%%%%%%%%%%
%%%%%%%%%%%%%%%%%%%%%%%%%%%%%%%%%%%%%%%%%%%%%%%%%%%%%%%%%%%%%%%%%%%%%%%%%%%
%%%%%%%%%%%%%%%%%%%%%%%%%%%%%%%%%%%%%%%%%%%%%%%%%%%%%%%%%%%%%%%%%%%%%%%%%%%
%%%%%%%%%%%%%%%%%%%%%%%%%%%%%%%%%%%%%%%%%%%%%%%%%%%%%%%%%%%%%%%%%%%%%%%%%%%

%%%%%%%%%%%%%%%%%%%%%%%%%%%%%%%%%%%%%%%%%%%%%%%%%%%%%%%%%%%%%%%%%%%%%%%%%%%
%%%%%%%%%%%%%%%%%%%%%%%%%%%%%%%%%%%%%%%%%%%%%%%%%%%%%%%%%%%%%%%%%%%%%%%%%%%
%%%%%%%%%%%%%%%%%%%%%%%%%%%%%%%%%%%%%%%%%%%%%%%%%%%%%%%%%%%%%%%%%%%%%%%%%%%
%%%%%%%%%%%%%%%%%%%%%%%%%%%%%%%%%%%%%%%%%%%%%%%%%%%%%%%%%%%%%%%%%%%%%%%%%%%
%%%%%%%%%%%%%%%%%%%%%%%%%%%%%%%%%%%%%%%%%%%%%%%%%%%%%%%%%%%%%%%%%%%%%%%%%%%

\subsubsection{Morphism of $A_\infty$-algebras}

This is the notion from (\cite{Kad80}).

\begin{definition}.
 A  morphism of $A_{\infty}$-algebras
$$ \{f_i\}:
(M,\{m_i\})\to (M',\{m'_i\})
$$
is a sequence $
 \{f_i:\otimes^iM\to M', i=1,2,. . . ,\ deg\ f_1=1-i \}
$ such that
 \begin{equation}
\label{morphism}
\begin{array}{l}
\sum_{k=0}^{i-1} \sum_{j=1}^{i-k} \pm \\
f_{i-j+1}(a_1\otimes . . . \otimes a_k \otimes m_j(a_{k+1}\otimes
. . . \otimes a_{k+j})\otimes . . .  \otimes a_i)=\\
\sum_{t=1}^{i} \sum_{k_1+...+k_t=i}\pm \\
m'_t(f_{k_1}(a_1\otimes. . . \otimes a_{k_1})\otimes. . . \otimes
f_{k_t}(a_{i-k_t+1}\otimes. . . \otimes a_{i})).
\end{array}
\end{equation}
\end{definition}

In particular for $n=1$ this condition gives $f_1m_1(a)=m'_1f_1(a)$, i.e.,
$f_1:(M,m_1)\to (M',m'_1)$ is a chain map; for $n=2$ it gives
$$
\begin{array}{l}
f_1m_2(a_1\otimes a_2)+m'_2(f_1(a_1)\otimes f_1(a_2))=\\
m'_1f_2(a_1\otimes a_2)+
f_2(m_1a_1\otimes a_2)+f_2(a_1\otimes m_1a_2),
\end{array}
$$
thus $f_1:(M,m_1,m_2)\to (M',m'_1,m'_2)$ is multiplicative just up to the chain homotopy $f_2$.

The composition of $A_{\infty}$ morphisms
$$
\{h_i\}:(M,\{m_i\})\stackrel{\{f_i\}}{\longrightarrow}
(M',\{m'_i\})\stackrel{\{g_i\}}{\longrightarrow} (M'',\{m''_i\})
$$
is defined as
\begin{equation}
\label{composition}
\begin{array}{l}
 h_n(a_1\otimes...\otimes
a_n)=\sum_{t=1}^{n}\sum_{k_1+...+k_t=n}\\
g_{n}(f_{k_1}(a_{1}\otimes...\otimes a_{k_1})\otimes ...\otimes
f_{k_t}(a_{n-k_t+1}\otimes...\otimes a_{n})).
\end{array}
\end{equation}
The bar construction argument (see (\ref{bartilda}) below) allows
to show that so defined composition satisfies the
condition (\ref{morphism}).

\emph{Special case:} A morphism of $A_\infty$-algebras $\{f_1,f_2=0,f_3=0,...\}:(M,\{m_1,m_2,m_3=0,m_4=0,...\})\to (M,\{m_1,m_2,m_3=0,m_4=0,...\})$ is ordinary map of $DG$-algebras. In fact, the category of dg algebras is a subcategory of the category of $A_\infty$-algebras.

%%%%%%%%%%%%%%%%%%%%%%%%%%
%%%%%%%%%%%%%%%%%%%%%%%%%%
%%%%%%%%%%%%%%%%%%%%%%%%%%

\subsubsection{Bar construction of an $A_{\infty}$-algebra}
\label{bartilda}

 Let  $(M,\{m_i\})$ be an $A_{\infty}$-algebra. The
structure maps $m_i$ define the map  $\beta:T^c(s^{-1}M)\to
s^{-1}M$ by $\beta [a_1,...,a_n]=[s^{-1}m_n(a_1\otimes ...\otimes
a_n)]$. Extending this $\beta$ as a coderivation (see the universal property (\ref{UnivTcVder})) we obtain
$d_\beta:T^c(s^{-1}M)\to T^c(s^{-1}M)$ which in fact looks as
$$
d_\beta[a_1,...,a_n]=\sum_k\pm [a_1,...,a_k,m_j(a_{k+1}\otimes
...\otimes a_{k+j}),a_{k+j+1},... a_n].
$$
The defining condition (\ref{ainfalg}) of an $A_{\infty}$-algebra guarantees that $d_\beta d_\beta = 0$: the composition of coderivations
$d_\beta d_\beta:T^c(s^{-1}M)\to T^c(s^{-1}M)$ is a coderivation and the defining condition (\ref{ainfalg})
is nothing else than that $p_1d_\beta d_\beta = 0$ and this is equivalent to $d_\beta d_\beta = 0$ since the tensor coalgebra is cofree.

The obtained dg coalgebra
$(T^c(s^{-1}M),d_\beta,\Delta)$ is called {\it bar construction}
of the $A_{\infty}$-algebra $(M,\{m_i\})$ and is denoted by
$B(M,\{m_i\})$.

For an $A_{\infty}$-algebra of type $(M,\{m_1,m_2,0,0,...\})$, i.e., for a strictly associative dg algebra, this
bar construction coincides with the ordinary bar construction of
this dg algebra.

\subsubsection{Bar interpretation of a morphim of $A_\infty$-algebras}

A morphism of $A_{\infty}$-algebras $ \{f_i\}: (M,\{m_i\})\to
(M',\{m'_i\}) $ defines a dg coalgebra map of bar constructions
$$F=B(\{f_i\}):B(M,\{m_i\})\to B(M',\{m'_i\}) $$
as follows: The collection $\{f_i\}$ defines the map $\alpha:T^c(s^{-1}M)\to
s^{-1}M$ by $ \alpha [a_1,...,a_n]=[s^{-1}f_n(a_1\otimes
...\otimes a_n)]. $ Extending this $\alpha$ as a coalgebra map
(\ref{UnivTV})) we
obtain $F:T^c(s^{-1}M)\to T^c(s^{-1}M)$ which in fact looks like
$$
F[a_1,...,a_n]=\sum \pm [f_{k_1}(a_{1}\otimes...\otimes a_{k_1}),
..., f_{k_t}(a_{n-k_t+1}\otimes...\otimes a_{n})].
$$
The defining condition (\ref{morphism}) of an $A_{\infty}$-morphism
is nothing else than $p_1d_{\beta'} F =p_1 F d_\beta $, and this is equivalent to $d_{\beta'} F =F d_\beta$,
so $F$ is a chain map.

Now we are able to show that the composition of
$A_{\infty}$-morphisms is well defined: to the composition of morphisms
(\ref{composition}) corresponds the composition of dg coalgebra
maps
$$B((M,\{m_i\}))\stackrel{B(\{f_i\})}{\longrightarrow}
B((M',\{m'_i\}))\stackrel{B(\{g_i\})}{\longrightarrow}
B((M'',\{m''_i\}))$$ which is a dg coalgebra map. Thus for
the projection $p_1B(\{g_i\})B(\{f_i\})$, i.e., for
the collection $\{h_i\}$, the condition (\ref{morphism}) is
satisfied.

%%%%%%%%%%%%%%%%%%%%%%%%%%%%%%%%%%%%%%%%%%%%%%%%%%%%%%%%%%%%%%%%%%%%%%%%%%%%%%%%%%%%%%%%%%%%%%%%%%%%%%%%

%%%%%%%%%%%%%%%%%%%%%%%%%%%%%%%%%%%%%%%%%%%%%%%%%%%%%%%%%%%%%%%%%%%%%%%%%%%%%%%%%%%%%%%%%%%%%%%%%%%%%%%%

%%%%%%%%%%%%%%%%%%%%%%%%%%%%%%%%%%%%%%%%%%%%%%%%%%%%%%%%%%%%%%%%%%%%%%%%%%%%%%%%%%%%%%%%%%%%%%%%%%%%%%%%

\subsubsection{Homotopy in the category of $A_\infty$-algebras}

Two morphisms of $A_\infty$-algebras
$$
\{f_i\},\{g_i\}:(M,\{m_i\})\rightarrow (M^\prime,\{m_i^\prime \})
$$
we call \emph{homotopic} if there exists a collection of homomorphisms
$$
\{h_i:(\otimes ^iM)\rightarrow M^{\prime },\quad i=1,2,...,\quad \deg
h_i=-i\},
$$
which satisfy the following condition

\begin{equation}
\label{Ainfhomot}
\begin{array}{c}
f_n(a_1\otimes ...\otimes a_n)-g_n(a_1\otimes ...\otimes a_n)\\
= \sum_{i+j=n+1}\sum_{k=0}^{n-j}h_i(a_1\otimes ...\otimes
a_k\otimes m_j(a_{k+1}\otimes ...\otimes a_{k+j})\otimes ...\otimes a_n)\\
+\sum_{k_1+...+k_t=n}
m'_t(f_{k_1}(a_1\otimes...\otimes a_{k_1})\otimes ... \otimes f_{k_{i-1}}(a_{k_1+...k_{i-2}+1}\otimes...\otimes a_{k_1+...k_{i-1}}) \\
\otimes h_{k_{i}}(a_{k_1+...+k_{i-1}+1}\otimes ... \otimes a_{k_1+...+k_{i}})\otimes g_{k_{i+1}}(a_{k_1+...k_{i}+1}\otimes...\otimes  a_{k_1+...k_{i+1}}) \\
\otimes g_{k_{t}}(a_{k_1+...k_{t-1}+1}\otimes...\otimes  a_{k_1+...+k_{t}}).
\end{array}
\end{equation}

In particular for $n=1$ this condition means
$$
f_1(a)-g_1(a)=m'_1h_1(a)+h_1(m_1a),
$$
that is, the chain maps $f_1,g_1:(M,m_1)\to (M',m'_1)$ ar chain homotopic.

For $n=2$ this condition means
\begin{equation}
\label{Ainfhomot2}
\begin{array}{c}
f_2(a_1\otimes a_2)-g_2(a_1\otimes a_2)=m'_1h_2(a_1\otimes a_2)+m'_2(f_1(a_1)\otimes h_1(a_2))\\
+m'_2(h_1(a_1)\otimes g_1(a_2))+

h_1(m_2(a_1\otimes a_2)+h_2(m_1a_1\otimes a_2)+h_2(a_1\otimes m_1a_2).
\end{array}
\end{equation}

\subsubsection{Bar interpretation of homotopy}

 The collections $\{f_i\}, \{g_i\},\ \{h_i\}$ define a homomorphism
$$
D:B(M,\{m_i\})\rightarrow B(M^{\prime },\{m_i^{\prime }\})
$$
by
$$
\begin{array}{c}
D(a_1\otimes...\otimes a_n)\\ = \sum_{k_1+...+k_t=n}
f_{k_1}(a_1\otimes...\otimes a_{k_1})\otimes ... \otimes f_{k_{i-1}}(a_{k_1+...k_{i-2}+1}\otimes...\otimes a_{k_1+...k_{i-1}}) \\
\otimes h_{k_{i}}(a_{k_1+...+k_{i-1}+1}\otimes ... \otimes a_{k_1+...+k_{i}})\otimes g_{k_{i+1}}(a_{k_1+...k_{i}+1}\otimes...\otimes  a_{k_1+...k_{i+1}})\\  \otimes
g_{k_{t}}(a_{k_1+...k_{t-1}+1}\otimes...\otimes  a_{k_1+...+k_{t}},
\end{array}
$$
which is a $(B(\{f_i\}),B(\{g_i\}))$-coderivation.

Besides, the condition (\ref{Ainfhomot}) means nothing else than
$$
p_1(B(\{f_i\})-B(\{g_i\}))=p_1(d_{m'}D+Dd_m),
$$
and this, again since the tensor coalgebra is cofree, gives
$$
B(\{f_i\})-B(\{g_i\})=d_{m'}D+Dd_m,
$$
that is, the dg coalgebra maps $B(\{f_i\})$ and $B(\{g_i\})$ are homotopic in the category of dg coalgebras.

%%%%%%%%%%%%%%%%%%%%%%%%%%%%%%%%%%%%%%%%%%%%%%%%%%%%%%%%%%%%%%%%%%%%%%%%%%%%%%%%%%%%%%%%%%%%%%%%%%%%%%%%%%%%%%%%%%%%%%%%%%%%%%%%%%

\subsubsection{Category $DASH$}

The category $DASH$ (Differential Algebras and Strongly Homotopy Multiplicative Maps) was first considered in Halperin and Stasheff's article \cite{HS}. The object are dg algebras, and a morphism  $\{f_i\}:(A,d,\mu)\to (A',d',\mu')$ is defined as a collection of homomorphisms
$$
\{f_i:\otimes^i A\to A;, i=1, 2, ..., deg f_i=1-i\}
$$
which satisfies the following conditions
\begin{equation}
\label{DASHfmor}
\begin{array}{c}
\sum_i f_{n}(a_1\otimes...\otimes a_{i-1}\otimes da_i\otimes a_{i+1}\otimes...\otimes a_n)\\ +

\sum_i f_{n-1}(a_1\otimes...\otimes a_{i-1}\otimes a_i\cdot a_{i+1}\otimes a_{i+2}\otimes...\otimes a_n)\\ =

\sum_i f_{i}(a_1\otimes...\otimes a_{i})\cdot f_{n-i}(a_{i+1}\otimes...\otimes a_n).

\end{array}
\end{equation}
In particular, for $n=1$ this condition gives $f_1d(a)=df_1(a)$, i.e.,
$f_1:(A,d)\to (A',d'_1)$ is a chain map; for $n=2$ it gives
$$
\begin{array}{l}
f_1(a_1\cdot a_2)-f_1(a_1)\cdot f_1(a_2)\\ =
d'f_2(a_1\otimes a_2)+
f_2(da_1\otimes a_2)+f_2(a_1\otimes da_2),
\end{array}
$$
thus $f_1:(A,d,\mu)\to (A',d',\mu')$ is multiplicative up to the homotopy $f_2$.

The existence of higher components $\{f_i,\ i=1,2,3,4,...\}$ is the reason why such a morphism is called a \emph{Strongly Homotopy Multiplicative Map} and the category is called DASH: \emph{Differential Algebras and Strongly Homotopy multiplicative maps}.

%%%%%%%%%%%%%%%%%%%%%%%%%%%%%%%%%%%%%%%%%%%%%%%%%%%%%%%%%%%%%%%%%%%%%%%%%%%%%%%%%%%%

In fact this is a \emph{full subcategory} of the category of $A_\infty$-algebras whose objects are ordinary dg-algebras $(A,d,\mu)=(A,\{m_1=d,m_2=\mu,m_3=0,m_4=0,...\})$: the defining condition of an $A_\infty$-morphism (\ref{morphism}) here looks like (\ref{DASHfmor}).

%%%%%%%%%%%%%%%%%%%%%%%%%%%%%%%%%%%%%%%%%%%%%%%%%%%%%%%%%%%%%%%%%%%%%%%%%%%%%%%%%%%%

Thus we have hierarchy of categories
$$
DGAlg\ \subset \ DASH\ \subset A_\infty,
$$
all three with notions of homotopies, and we have the commutative diagram of functors
$$
\begin{array}{ccc}

DGAlg&\ \subset \ DASH\ \subset &A_\infty\\

B\searrow & B\downarrow & \swarrow B \\

& DGCoalg.

\end{array}
$$

\subsubsection{Minimality}
\label{minimal}

Let $ \{f_i\}: (M,\{m_i\})\to (M',\{m'_i\})$ be a morphism of
$A_{\infty}$-algebras. It follows from (\ref{morphism}) that the
first component $f_1:(M,m_1)\to (M',m'_1)$ is a chain map.

Let us define a weak equivalence of $A_{\infty} $-algebras as a
morphism $\{f_i\}$ for which $B(\{f_i\})$ is a weak equivalence
(homology isomorphism) of dg coalgebras. The standard spectral
sequence argument allows to prove the following
\begin{proposition}.
\label{mapweak}
 A morphism of $A_{\infty}$-algebras is a weak
equivalence if and only if its first component $f_1:(M,m_1)\to
(M',m'_1)$ is a weak equivalence of chain complexes.
\end{proposition}
%%%%%%%%%%%%%%%%%%%%%%%%%%%%%%%

Furthermore, it is easy to see that if $ \{f_i\}: (M,\{m_i\})\to (M',\{m'_i\})$ is an isomorphism of $A_\infty$ algebras then
its first component $f_1:(M,m_1)\to (M',m'_1)$ is an isomorphism of chaian complexes. Conversely, if $f_1$ is an isomorphism, then
$\{f_i\}$ is an isomorphism of $A_\infty$ algebras:
The components of the opposite morphism $
\{g_i\}: (M',\{m'_i\})\to (M,\{m_i\})$ can be solved inductively
from the equation $\{g_i\}\{f_i\}=\{id_M,0,0,...\}$. Thus we have the
\begin{proposition}.
\label{mapiso}
 A morphism of $A_{\infty}$-algebras is an
isomorphism if and only if its first component $f_1:(M,m_1)\to
(M',m'_1)$ is an isomorphism of dg modules.
\end{proposition}
\begin{definition}.
An $A_{\infty}$-algebra $(M,\{m_i\})$ we call minimal if $m_1=0$.
\end{definition}
In this case $(M,m_2)$ is a \emph{ strictly} associative graded
algebra.

The above propositions easily imply:
\begin{proposition}.
\label{miniso}
 Each weak equivalence of minimal $A_{\infty}$-algebras is an isomorphism.
\end{proposition}

\noindent {\bf Proof.} Suppose
$f=\{f_i\}:(M,\{m_i\})\rightarrow (M^{\prime },\{m_i^{\prime }\})$
is a weak equivalence of $A_\infty$-algebras. Then by (\ref{mapweak}) the chain map
$f_1:(M,m_1=0)\rightarrow (M^{\prime },m_1^{\prime }=0)$ induces an isomorphism of homology modules
$f^*_1:H(M,m_1=0)\rightarrow H(M^{\prime },m_1^{\prime }=0)$, but $H(M,m_1=0)=M$ and
$H(M',m_1=0)=M'$. So $f_1$ is an isomorphism and by of (\ref{mapiso}) $f$ is an isomorphism too.

This fact motivates the word {\it minimal} in this notion since Sullivan's minimal model has a similar property.

\subsection{$A_\infty$ Twisting Cochains}

Here we are going to generalize the above material about twisting cochains and the functor $D$ from the case of dg algebras to the case of $A_\infty$ algebras, see \cite{Kad86}, see also \cite{Smi80}.

Start with a dg coalgebra $(C,d_C,\nabla:C\to C\otimes C)$ and an $A_\infty$-algebra $(M,\{m_i\})$.
We'll recall the following multicooperations
$$
\begin{array}{c}
\nabla^i:C\to C\otimes ... (i\ times)...\otimes C\to C,\ i=1,2,\ ...\ ,\\ \nabla^1=id,\ \nabla^2=\nabla,\ \nabla^i=(\nabla^{i-1}\otimes id)\nabla.
\end{array}
$$

\subsubsection{The notion of $A_\infty$-twisting cochain}
An $A_\infty$-twisting cochain \cite{Kad86}, \cite{Kad87} we have defined as a homomorphism $\phi:C\to M$ of degree $-1$ (that is $\phi:C^*\to M^{*-1}$) satisfying the condition
 \begin{equation}
\label{AinfTw}
\phi d_C=\sum^\infty_{i=1} m_i(\phi\otimes...\otimes\phi)\nabla^i.
\end{equation}
Let $Tw(C,M)$ be the set of all such twisting cochains.

Note that if for $(M,\{m_i\})$ one has $m_{>2}=0$, i.e., it is an ordinary dg algebra $(M,\{m_1,m_2,0,0,\ ...\ \})$ then this notion coincides with Brown' classical notion.

%\subsubsection{Bar Interpretation}

\subsubsection{Bar interpretation}
Any $A_\infty$ twisting cochain $\phi:C\to M$ induces a dg-coalgebra map $B(\phi):C\to BM$ by
$$
B(\phi)=\sum_i (\phi\otimes \ ...\ \otimes \phi)\nabla^i.
$$
So defined $B(\phi)$ automatically is a coalgebra map and the condition (\ref{AinfTw}) guarantees that it is a chain map.

Conversely, any dg coalgebra map  $f:C\to BM$ is $B(\phi)$ for $\phi=p\circ f:C\to BM\to M$.
In fact we have a bijection
$$
Mor_{DGCoalg}(C,BM)\leftrightarrow Tw(C,M).
$$

\subsubsection{Equivalence of $A_\infty$-twisting cochains.}
Two $A_\infty$-twisting cochains $\phi,\psi:C\to M$ are called equivalent \cite{Kad86} if there exists $c:C\to M,\ deg\ c=0$, such that
 \begin{equation}
\label{AinfTwEq}
\phi-\psi=cd_C+\sum_{i=1}^\infty \sum_{k=1}^{i-1} m_i(\phi\otimes...(k\ times)...\phi\otimes c\otimes \psi\otimes...\psi)\nabla^i,
\end{equation}
notation $\phi\sim_c \psi$.

Again, if for $(M,\{m_i\})$ one has $m_{>2}=0$, i.e., it is an ordinary dg algebra $(M,\{m_1,m_2,0,0,\ ...\ \})$ then this notion coincides with Berikashvili's described above equivalence.

\subsubsection{Bar interpretation}
If $\phi\sim_c\psi$ then $B(\phi)$ and $B(\psi)$ are homotopic in the category of dg coalgebras: a chain homotopy $D(c):C\to BA$ is given by
\begin{equation}\label{Dc}
D(c)=\sum_{i=1}^\infty \sum_{k=1}^{i-1} (\phi\otimes...(k\ times)...\phi\otimes c\otimes \psi\otimes...\psi)\nabla^i.
\end{equation}
So defined $D$ automatically is a $B(\phi),B(\psi)$-coderivation and the condition (\ref{AinfTwEq}) guarantees that it realizes chain homotopy, that is, $B(\phi)-B(\psi)=d_BD(c)+D(c)d_C$.

%%%%%%%%%%%%%%%%%%%%%%%%%%%%%%%%%%%%%%%%%%%%%%%%%%%%%%%%%%%%%%%%%%%%%%%%%%%%%%%%%%%%%%%%%%%%%%%%%%%
%%%%%%%%%%%%%%%%%%%%%%%%%%%%%%%%%%%%%%%%%%%%%%%%%%%%%%%%%%%%%%%%%%%%%%%%%%%%%%%%%%%%%%%%%%%%%%%%%%%
%%%%%%%%%%%%%%%%%%%%%%%%%%%%%%%%%%%%%%%%%%%%%%%%%%%%%%%%%%%%%%%%%%%%%%%%%%%%%%%%%%%%%%%%%%%%%%%%%%%
%%%%%%%%%%%%%%%%%%%%%%%%%%%%%%%%%%%%%%%%%%%%%%%%%%%%%%%%%%%%%%%%%%%%%%%%%%%%%%%%%%%%%%%%%%%%%%%%%%%
%%%%%%%%%%%%%%%%%%%%%%%%%%%%%%%%%%%%%%%%%%%%%%%%%%%%%%%%%%%%%%%%%%%%%%%%%%%%%%%%%%%%%%%%%%%%%%%%%%%

%%%%%%%%%%%%%%%%%%%%%%%%%%%%%%%%%%%%%%%%%%%%%%%%%%%%%%%%%%%%%%%%%%%%%%%%%%%%%%%%%%%%%%%%%%%%%%%%%%%
%%%%%%%%%%%%%%%%%%%%%%%%%%%%%%%%%%%%%%%%%%%%%%%%%%%%%%%%%%%%%%%%%%%%%%%%%%%%%%%%%%%%%%%%%%%%%%%%%%%
%%%%%%%%%%%%%%%%%%%%%%%%%%%%%%%%%%%%%%%%%%%%%%%%%%%%%%%%%%%%%%%%%%%%%%%%%%%%%%%%%%%%%%%%%%%%%%%%%%%
%%%%%%%%%%%%%%%%%%%%%%%%%%%%%%%%%%%%%%%%%%%%%%%%%%%%%%%%%%%%%%%%%%%%%%%%%%%%%%%%%%%%%%%%%%%%%%%%%%%
%%%%%%%%%%%%%%%%%%%%%%%%%%%%%%%%%%%%%%%%%%%%%%%%%%%%%%%%%%%%%%%%%%%%%%%%%%%%%%%%%%%%%%%%%%%%%%%%%%%

\subsubsection{The Functor $\tilde{D}$.}
The functor $\tilde{D}(C,M)$ (\cite{Kad86}) is defined as the factorset $$\tilde{D}(C,A)=\frac{Tw(C,M)}{\sim}.$$

%%%%%%%%%%%%%%%%%%%%%%%%%%%%%%%%%%%%%%%%%%%%%%%%%%%%%%%%%%%%%%%%%%%%%%%%%%%%%%%%%%%%%%%%%%%%%%%%%%%
%%%%%%%%%%%%%%%%%%%%%%%%%%%%%%%%%%%%%%%%%%%%%%%%%%%%%%%%%%%%%%%%%%%%%%%%%%%%%%%%%%%%%%%%%%%%%%%%%%%
%%%%%%%%%%%%%%%%%%%%%%%%%%%%%%%%%%%%%%%%%%%%%%%%%%%%%%%%%%%%%%%%%%%%%%%%%%%%%%%%%%%%%%%%%%%%%%%%%%%
%%%%%%%%%%%%%%%%%%%%%%%%%%%%%%%%%%%%%%%%%%%%%%%%%%%%%%%%%%%%%%%%%%%%%%%%%%%%%%%%%%%%%%%%%%%%%%%%%%%
%%%%%%%%%%%%%%%%%%%%%%%%%%%%%%%%%%%%%%%%%%%%%%%%%%%%%%%%%%%%%%%%%%%%%%%%%%%%%%%%%%%%%%%%%%%%%%%%%%%

Any $A_\infty$-algebra map $f=\{f_i\}:M\to M'$ induces the map $Tw(C,M)\to Tw(C,M')$: if $\phi:C\to M$ is an $A_\infty$-twisting cochain then so is
\begin{equation}\label{AinfTwMap}
f\phi=p'B(f)B(\phi)=\sum_i f_i(\phi\otimes...\otimes \phi)\nabla^i:C\to B(M)\to B(M')\to M'.
\end{equation}
Furthermore, if $\phi  \sim_c\psi$ then $f \phi\sim_{C}f \psi$ with
$$
C=\sum_{i=1}^\infty \sum_{k=1}^{i-1} (\phi\otimes...(k\ times)...\phi\otimes c\otimes \psi\otimes...\psi)\nabla^i.
$$

Thus we have a map $\tilde{D}(f):\tilde{D}(C,M)\to \tilde{D}(C,M')$.

Dually, any dg coalgebra map $g:C'\to C$ induces a map $Tw(C,M)\to Tw(C',M)$: if $\phi:C\to M$ is an $A_\infty$-twisting cochain so is the composition $\phi g:C'\to C\to M$.
Moreover, if $\phi  \sim_c\psi$ then $\phi g\sim_{ cg} \psi g$.
Thus we have a map $D(g):D(C,M)\to D(C',M)$.

%Let $(K,d_K,\nabla_K)$ be a dg colagebra with free $K_i$s and $(A,d_A,\mu_A)$ be a connected dg algebra.

%%%%%%%%%%%%%%%%%%%%%%%%%%%%%%%%%%%%%%%%%%%%%%%%%%%%%%%%%%%%%%%%%%%%%%%%%%%%%%%%%%%%%%%%%%%%%%%%%%%%%%%%%%%%%%%%%%%%%%%%%%%%%%%%%%%%%%%%%%%%%%

\subsubsection{Bar interpretation}
Assigning to an $A_\infty$-twisting cohain $\phi:C\to M$ the dg coalgebra map $B(\phi):C\to BM$ and having in mind that $\phi\sim_c\psi$ implies $B(f\phi)\sim_{D(c)}B(\psi)$ we obtain the

\begin{theorem}.
There is a bijection $\tilde{D}(C,M)\leftrightarrow [C,BM]_{DGCoalg}$.
\end{theorem}

Taking $C=B(M')$, the bar construction of an $A_\infty$-algebra $(M',\{m_i'\})$ and having in mind that $[B(M'),B(M)]_{DGCoalg}$ is bijective to  $[M',M]_{A_\infty}$, the set of homotopy classes in the category of $A_\infty$-algebras,
we obtain the
\begin{proposition}.\label{D=[]}
There is a bijection $D(B(M'),M)\leftrightarrow [M',M]_{A_\infty}$.
\end{proposition}

%%%%%%%%%%%%%%%%%%%%%%%%%%%%%%%%%%%%%%%%%%%%%%%%%%%%%%%%%%%%%%%%%%%%%%%%%%%%%%%%%%%%%%%%%%%%%%%%%%%%%%%%%%%%%%%%%%%%%%%%%%%%%%%%%%%%%%%%%%%%%%

\subsubsection{Bijections.}
 The following property of the functor $\tilde{D}$ makes it useful in the homotopy classification of maps.
\begin{theorem}.\label{AinfDwe} (See \cite{Kad86}) (a)
If $\{f_i\}:M\to M'$ is a weak equivalence of $A_\infty$-algebras, then
$$
D(\{f_i\}):D(C,M)\to D(C,M')
$$
is a bijection.

(b)
If $g:C\to C'$ is a weak equivalence of connected dg coalgebras (i.e., a homology isomorphism), then
$$
D(g):D(C',M')\to D(C,M')
$$
is a bijection.
\end{theorem}
Combining we obtain
\begin{proposition}.\label{Bijec}
A weak equivalence of $A_\infty$-algebras $\{f_i\}:M\to M'$ and a weak equivalence of dg coalgebras $g:C\to C'$ induce bijections
$$
D(C,M)\stackrel{D(\{f_i\})}{\to}D(C,M')\stackrel{D(g)}{\leftarrow}D(C',M').
$$
\end{proposition}

%%%%%%%%%%%%%%%%%%%%%%%%%%%%%%%%%%%%%%%%%%%%%%%%%%%%%%%%%%%%%%%%%%%%%%%%%%%%%%%%%%%%%%%%%%%%%%%%%%%
%%%%%%%%%%%%%%%%%%%%%%%%%%%%%%%%%%%%%%%%%%%%%%%%%%%%%%%%%%%%%%%%%%%%%%%%%%%%%%%%%%%%%%%%%%%%%%%%%%%
%%%%%%%%%%%%%%%%%%%%%%%%%%%%%%%%%%%%%%%%%%%%%%%%%%%%%%%%%%%%%%%%%%%%%%%%%%%%%%%%%%%%%%%%%%%%%%%%%%%
%%%%%%%%%%%%%%%%%%%%%%%%%%%%%%%%%%%%%%%%%%%%%%%%%%%%%%%%%%%%%%%%%%%%%%%%%%%%%%%%%%%%%%%%%%%%%%%%%%%
%%%%%%%%%%%%%%%%%%%%%%%%%%%%%%%%%%%%%%%%%%%%%%%%%%%%%%%%%%%%%%%%%%%%%%%%%%%%%%%%%%%%%%%%%%%%%%%%%%%

The following theorem is an analogue of Berikashvilis's Theorem \ref{BerThm} for $A_\infty$-algebras and was proved in \cite{Kad86}.

\begin{theorem}.\label{kadthmDtilde}
Let $(K,d_K,\nabla_K)$ be a dg coalgebra with free $K_i$ and let  $(M,\{m_i\})$ be a connected dg algebra.
If $f=\{f_i\}:(M,\{m_i\})\to (M',\{m_i'\})$ is a weak equivalence of $A_\infty$-algebras  then
$$
D_\infty(f):D_\infty(K,M)\to D_\infty(K,M')
$$
is a bijection.
\end{theorem}

In fact this theorem means that $[K,BM]\to [K,BM']$ is a bijection. The theorem consists of two parts:

\noindent \textbf{Surjectivity.} Any $A_\infty$-twisting cochain $\phi:K\to M'$ can be lifted to an $A_\infty$-twisting cochain $\psi:K\to M$ so that $\phi \sim f\psi$.

\noindent \textbf{Injectivity.} If $\psi\sim \psi'\in T_\infty(K,M)$ then $f\psi\sim f\psi\in T_\infty(K,M')$.

%%%%%%%%%%%%%%%%%%%%%%%%%%%%%%%%%%%%%%%%%%%%%%%%%%%%%%%%%%%%%%%%%%%%%%%%%%%%%%%%%%%%%%%%%%%%%%%%%%%
%%%%%%%%%%%%%%%%%%%%%%%%%%%%%%%%%%%%%%%%%%%%%%%%%%%%%%%%%%%%%%%%%%%%%%%%%%%%%%%%%%%%%%%%%%%%%%%%%%%
%%%%%%%%%%%%%%%%%%%%%%%%%%%%%%%%%%%%%%%%%%%%%%%%%%%%%%%%%%%%%%%%%%%%%%%%%%%%%%%%%%%%%%%%%%%%%%%%%%%
%%%%%%%%%%%%%%%%%%%%%%%%%%%%%%%%%%%%%%%%%%%%%%%%%%%%%%%%%%%%%%%%%%%%%%%%%%%%%%%%%%%%%%%%%%%%%%%%%%%
%%%%%%%%%%%%%%%%%%%%%%%%%%%%%%%%%%%%%%%%%%%%%%%%%%%%%%%%%%%%%%%%%%%%%%%%%%%%%%%%%%%%%%%%%%%%%%%%%%%

\subsubsection{Lifting of $A_\infty$-Twisting Cochains}

Below we'll use the surjectivity part of this theorem whose proof we sketch here.

\begin{theorem}. \label{lifta}
Let $(K,d_K,\nabla_K)$ be a dg coalgebra with free $K_i$ and let $f:(M,\{m_i\})\to (M',\{m'_i\})$ be a weak equivalence of connected $A_\infty$-algebras (that is, $M_0=M'_0=0$). Then for an arbitrary $A_\infty$-twisting cochain
$$
\phi=\phi_2+\phi_3+\ ...\ +\phi_n+\ ...\ :K\to M'
$$
there exists an $A_\infty$-twisting cochain
$$
\psi=\psi_2+\psi_3+\ ...\ +\psi_n+\ ...\ :K\to M
$$
such that $\phi \sim f\psi$.
\end{theorem}

\noindent \textbf{Proof.} Start with an $A_\infty$-twisting cochain
$$
\phi=\phi_2+\phi_3+\ ...\ +\phi_n+\ ...\ :K\to M'.
$$
The condition (\ref{AinfTw}) gives $m_1\phi_2=0$, and since $f:(M,m_1)\to (M'm'_1)$ is homology isomorphism then there exist $\psi_2:K_2\to M_1$ and $c'_2:K_2\to M'_2$ such that $m_1\psi_2=0$ and $f\psi_2=\phi_2+m'_1c'_2$. Perturbing $\phi$ by this $c'_2$ we obtain $F_{c'_2}\phi$ for which
$(F_{c'_2}\phi)_2=\phi_2+m'_1c'_2=f\psi_2$.

%%%%%%%%%%%%%%%%%%%%%%%%%%%%%%%%%%%%%%%%%%%%%%%%%%%%%%%%%%%%%%%%%%%%%%%%%%%%%%%%%%%%%%%%%%%%%%%%%%%%%%%%%%%%%%%%%%%%%%%%%%%%%%%%%%%%%%%%

Assume now that we already have
$\psi_2,\ \psi_3,\  ...\ ,\psi_{n-1}$ which satisfy (\ref{AinfTw}) and (\ref{AinfTwMap}) in appropriate dimensions, that is
\begin{equation}\label{kadtinfk}
m_1\psi_k=\psi_{k-1} d_K+\sum_{i=2}^{\int(k/2)}\sum_{k_1+...k_i=k} m_i(\psi_{k_1}\otimes\ ...\ \otimes\psi_{k_i})\nabla_K^k,\ \ k=2,3,...,n-1
\end{equation}
and
\begin{equation}\label{kadfphk}
\phi_k=\sum_{i=1}^{\int(k/2)}\sum_{k_1+...+k_i=k}f_i(\psi_{k_1}\otimes\ ...\ \otimes \psi_{k_i})\nabla_K^i,\ \ k=1,2,...,n-1.
\end{equation}
We need the next component $\psi_n:K_n\to A_{n-1}$ satisfying the condition (\ref{kadtinfk}) for $k=n$
\begin{equation}\label{kadphn}
m_1\psi_n=\psi_{n-1} d_K+\sum_{i=2}^{\int(n/2)}\sum_{k_1+...k_i=n} m_i(\psi_{k_1}\otimes\ ...\ \otimes\psi_{k_i})\nabla_K^k,
\end{equation}
and $c'_n:K_n\to A'_n$ such that
\begin{equation}\label{kadfphn}
m_1c'_n+\phi_n=\sum_{i=1}^{\int(n/2)}\sum_{k_1+...+k_i=n}f_i(\psi_{k_1}\otimes\ ...\ \otimes \psi_{k_i})\nabla_K^i.
\end{equation}
Then perturbing $\phi$ by $c'_n$ we obtain new $\phi_n$ for which the condition (\ref{kadfphk}) will be satisfied for $k=n$, and this will complete the proof.

%%%%%%%%%%%%%%%%%%%%%%%%%%%%%%%%%%%%%%%%%%%%%%%%%%%%%%%%%%%%%%%%%%%%%%%%%%%%%%%%%%%%%%%%%%%%%%%%%%%%%%%%%%%%%%%%%%%%%%%%%%%%%%%%%%%%%%%%%%%%%%%%%%%%

Let us put
$$
\begin{array}{l}
U_n=\psi_{n-1} d_K+\sum_{i=2}^{int(n/2)}\sum_{k_1+...k_i=n} m_i(\psi_{k_1}\otimes\ ...\ \otimes\psi_{k_i})\nabla_K^k,\\
U'_n=\phi_{n-1} d_K+\sum_{i=2}^{int(n/2)}\sum_{k_1+...k_i=n} m_i(\phi_{k_1}\otimes\ ...\ \otimes\phi_{k_i})\nabla_K^k\\
V_n=\sum_{i=2}^{int(n/2)}\sum_{k_1+...+k_i=n}f_i(\psi_{k_1}\otimes\ ...\ \otimes \psi_{k_i})\nabla_K^i.
\end{array}
$$
Then the needed conditions (\ref{kadphn}) and  (\ref{kadfphn}) become
$$
U_n =m_1\psi_n,\ \ \ m'_1c'_n+\phi_n=f_1\psi_n+V_n.
$$

%%%%%%%%%%%%%%%%%%%%%%%%%%%%%%%%%%%%%%%%%%%%%%%%%%%%%%%%%%%%%%%%%%%%%%%%%%%%%%%%%%%%%%55

First, it is possible to check that $m_1U_n=0$ and $f_1U_n=U'_n+m'_1V_n$. Having in mind that $U'_n=m'_1\phi_n$ the last condition means
$f_1U_n=m'_1(\phi_n+V_n)$.
So $U_n$ is a map to $m_1$-cycles
$$U_n:K_n\to Z(M_{n-2})\subset M_{n-2}$$
and $f_1U_n$ maps $K_n$ to boundaries
$$f_1U_n:K_n\to B(M'_{n-2})\subset M_{n-2}.$$
Then since $f_1:(M,m_1)\to (M',m'_1)$ is a homology isomorphism, there exist $\overline{\psi}_n:K_n\to M_{n-1}$ such that
$m_1\overline{\psi}_n=U_n.$ So our $\overline{\psi_n}$ satisfies (\ref{kadphn}). Now let us perturb this $\overline{\phi}$ in order to catch the condition (\ref{kadfphn}) too.

Using the above equality $f_1U_n=m'_1(\phi_n+V_n)$, we have
$$
m'_1f_1\overline{\psi_n}=f_1m_1\overline{\psi_n}=f_1U_n=m'_1(\phi_n+V_n)
$$
i.e., $m'_1(f_1\overline{\psi_n}-(\phi_n+V_n))=0$. Thus
$$
z'_n=(f_1\overline{\psi_n}-(\phi_n+V_n)):K_n\to Z(M'_{n-1}),
$$
and again, since $f_1:(M,m_1)\to (M',m'_1)$ is a homology isomorphism there exist $z_n:K_n\to Z(M_{n-1})$ and $c'_n:K_n\to M'_n$ such that
$z'_n=f_1z_n-m'_1c'_n$. Let's define $\psi_n=\overline{\psi_n}-z_n$. Then
$$
f_1\psi_n=f_1\overline{\psi_n}-f_1z_n=(z'_n+\phi_n+V_n)-(z'_n-m'_1c'_n)=\phi_n+V_n.
$$

\subsection{$C_\infty$-algebras}

This is the commutative version of the notion of
$A_\infty$-algebra. For an $A_\infty$-algebra $(M,\{m_i\})$ it is
clear what it means for the operation $m_2:M\otimes M\to M$, but what about the commutativity of the higher operations
$m_i:M\otimes ...\otimes M\to M, \ i\geq 3$? We are going to describe this now.

\subsubsection{Shuffle product}

It was mentioned in (\ref{Shuffle}) that the tensor algebra $(T(V),\mu)$ and the tensor coalgebra $(T^c(V),\Delta)$
coincide as graded modules, but the
multiplication $\mu$ of $T(V)$ and the comultiplication $\Delta$ of $T^c(V)$ are
not compatible with each other, so they do not define a graded
bialgebra structure on $T(V)=T^c(V)$.

But, as it was indicated in (\ref{Shuffle}) there exists the {\it shuffle} multiplication
$\mu_{sh}:T^c(V)\otimes T^c(V)\to T^c(V)$ introduced by Eilenberg
and MacLane which turns $(T^c(V),\Delta,\mu_{sh})$ into a graded
bialgebra.
This multiplication is defined as a graded coalgebra map induced
by the universal property of $T^c(V)$ by $\alpha:T^c(V)\otimes
T^c(V)\to V$ given by $\alpha(v\otimes 1)=\alpha(1\otimes v)=v$
and $\alpha=0$ otherwise. This multiplication is associative and
in fact is given by
$$
\mu_{sh}([a_1,...,a_m]\otimes[a_{i+1},...,a_{n}])=\sum \pm
[a_{\sigma(1)},...,a_{\sigma(n)}]),
$$
where the summation is taken over all $(m,n)$-shuffles. That is,
over all permutations of the set $(1,2,...,n+m)$ which satisfy the
condition: $i<j$ if $1\leq \sigma (i)<\sigma (j)\leq n$ or
$n+1\leq \sigma (i)<\sigma (j)\leq n+m$. In particular
$$
[a]*_{sh}[b]=[a,b]\pm[b,a],
$$
$$
[a]*_{sh}[b,c]=[a,b,c]\pm[b,a,c]\pm[b,c,a].
$$

\subsubsection{The notion of $C_\infty$-algebra}
\label{CinftyDef}

Now we can define the notion of $C_\infty$-algebra, which is a
commutative version of Stasheff's notion of $A_\infty$-algebra.

%%%%%%%%%%%%%%%%%%%%%%%%%%%%%%%%%%%%%%%%%%%
\begin{definition}. (\cite{Smi80}, \cite{Kad88},\cite{Markl}, \cite{GJ94}) A $C_{\infty}$-algebra is an
$A_{\infty}$-algebra $(M,\{m_i\})$ which additionally satisfies
the following condition: each operation $m_i$ vanishes on
shuffles, that is, for $a_1,...,a_i\in M$ and $k=1,2,...,i-1$
\begin{equation}
\label{m0onshuffle} m_i(\mu_{sh}((a_1\otimes ...\otimes
 a_k)\otimes(a_{k+1}\otimes ...\otimes a_i)))=0.
\end{equation}
\end{definition}

In particular, this gives
$$
\begin{array}{l}
m_2(a\otimes b)+m_2(b\otimes a)=0, \\

m_3(a\otimes b\otimes c)+m_3(a\otimes c\otimes b)+m_3(c\otimes a\otimes b)=0, \\
...
\end{array}
$$
i.e., the maps $f_i$ vanish on shuffles.

\begin{definition}.
A morphism of $C_{\infty}$-algebras is defined as a morphism of
$A_{\infty}$-algebras $ \{f_i\}: (M,\{m_i\})\to (M',\{m'_i\}) $
whose components $f_i$ vanish on shuffles, that is
\begin{equation}
\label{f0onshuffle} f_i((\mu_{sh}(a_1\otimes ...\otimes
 a_k)\otimes(a_{k+1}\otimes ...\otimes a_i)))=0.
\end{equation}
\end{definition}
The composition is defined as in the $A_{\infty}$-case and the bar
construction argument (see  (\ref{BarCinfty}) below) allows to
show that the composition is a $C_{\infty}$-morphism.

In particular, this gives
$$
\begin{array}{l}
f_2(a\otimes b)+m_2(b\otimes a)=0, \\

f_3(a\otimes b\otimes c)+m_3(a\otimes c\otimes b)+m_3(c\otimes a\otimes b)=0, \\
...
\end{array}
$$
i.e., the maps $f_i$ vanish on shuffles.

In particular, for the operation $m_2$ we have $m_2(a\otimes b\pm
b\otimes a)=0$, so a $C_{\infty}$-algebra of type
$(M,\{m_1,m_2,0,0,...\})$ is a commutative dg algebra (cdga) with
the differential $m_1$ and strictly associative and commutative
multiplication $m_2$. Thus the category of commutative dg algebras is a
subcategory of the category of $C_{\infty}$-algebras.

\subsubsection{Bar interpretation of a $C_{\infty}$-algebra}
\label{BarCinfty}

The notion of $C_{\infty}$-algebra is motivated by the following
observation. It is well known that if a dg algebra $(A,d,\mu)$ is graded commutative
then the differential of the bar construction $BA$ is not only a
coderivation but also a derivation with respect to the shuffle
product, so the bar construction $(BA,d_\beta,\Delta,\mu_{sh})$ of
a cdga is a dg bialgebra, see (\ref{BarCommut}).

By definition the bar construction of an $A_{\infty}$-algebra
$(M,\{m_i\})$ is a dg coalgebra
$\tilde{B}(M)=(T^c(s^{-1}M),d_\beta,\Delta)$.

But if $(M,\{m_i\})$ is a $C_{\infty}$-algebra, then
$\tilde{B}(M)$ becomes a dg bialgebra:
\begin{proposition}.
For an $A_{\infty}$-algebra $(M,\{m_i\})$ the differential of the
bar construction $d_\beta$ is a derivation with respect to the
shuffle product if and only if each operation $m_i$ vanishes on
shuffles, that is $(M,\{m_i\})$ is a $C_{\infty}$-algebra.
\end{proposition}
\noindent {\bf Proof.} The map $ \Phi:T^c(s^{-1}M)\otimes
T^c(s^{-1}M)\to T^c(s^{-1}M) $ defined by
$\Phi=d_\beta\mu_{sh}-\mu_{sh}(d_\beta\otimes id+id\otimes
d_\beta)$ is a coderivation, see the arguments in (\ref{Gcoalg}). Thus, according to the universal property
of the tensor coalgebra (\ref{UnivTcVder}) the map $\Phi$ is trivial if and only if
$p_1\Phi=0$, and the last condition means exactly
(\ref{m0onshuffle}).

\subsubsection{Bar interpretation of a morphism of $C_{\infty}$-algebras}
\label{BarCinftyMor}

A morphism of $C_\infty$ algebras has an analogous interpretation.

\begin{proposition}. Let
$ \{f_i\}: (M,\{m_i\})\to (M',\{m'_i\}) $ be an
$A_{\infty}$-algebra morphism of $C_{\infty}$-algebras. Then the
induced map of bar constructions $ \tilde{B}\{f_i\}$ is a map of
dg bialgebras if and only if each $f_i$ vanishes on shuffles, that
is $\{f_i\}$ is a morphism of $C_{\infty}$-algebras.
\end{proposition}
\noindent {\bf Proof.} The map
$\Psi=\tilde{B}\{f_i\}\mu_{sh}-\mu_{sh}(\tilde{B}\{f_i\}\otimes
\tilde{B}\{f_i\}) $ is a coderivation, see the arguments in (\ref{Gcoalg}). Thus, according to the
universal property of the tensor coalgebra (\ref{UnivTcVder}) the map $\Psi$ is trivial if
and only if $p_1\Psi=0$, and the last condition means exactly
(\ref{f0onshuffle}).

Thus the bar functor maps the subcategory of $C_{\infty}$-algebras
to the category of dg bialgebras.

%%%%%%%%%%%%%%%%%%%%%%%%%%%%%%%%%%%%%%%
%%%%%%%%%%%%%%%%%%%%%%%%%%%%%%%%%%%%%%%
%%%%%%%%%%%%%%%%%%%%%%%%%%%%%%%%%%%%%%%

\subsubsection{Adjunctions}  \label{adj}

The bar and cobar functors
$$
B:DGAlg\to DGCoalg,\  \Omega :DGCoalg\to DGAlg
$$
are adjoint and there exist standard weak equivalences $\Omega
B(A)\to A,\ \ C\to B\Omega C$. So $\Omega B(A)\to A$ is a {\it
free resolution} of a dga $A$.

If $A$ is commutative, the cobar-bar resolution is out of the
category: $\Omega B(A)$ is not commutative.

In this case instead the cobar-bar functors we must use the
adjoint functors $\Gamma,\ {\cal A}$, see \cite{SS93}, which we
describe now.

For a commutative dg algebra the bar construction is a dg
bialgebra, so the restriction of the bar construction is the
functor $B:CDGAlg \to DGBialg$. Furthermore, the functor of
indecomposables $Q:DGBialg\to DGLieCoalg$ maps the category of dg
bialgebras to the category of dg Lie coalgebras. Let $\Gamma$ be
the composition
$$
\Gamma: CDGAlg \stackrel{B}{\to}
DGBialg\stackrel{Q}{\to}DGLieCoalg.
$$
There is the adjoint of $\Gamma$ $ {\cal A}:DGLieCoalg \to CDGAlg,
$ which is dual to the Chevalley-Eilenberg functor. There is the
standard weak equivalence ${\cal A}\Gamma A\to A$. This is the commutative analogue of the
standard weak equivalence $\Omega B(A)\to A$ from (\ref{AdjBarCobar}).

\section{Minimal $A_{\infty}$  and $C_{\infty}$-Algebras and
Hochschild and Harrison Cohomology} \label{hochschild}

Here we present the connection of the notion of minimal $A_\infty$
(resp. $C_\infty$) algebra with the Hochschild (resp. Harrison)
cochain complexes, studied in \cite{Kad88}, see also \cite{Kad07}.

\subsection{Hochschild Cohomology}

Below we describe the classical Hochschild cochain complex \cite{Hoch} and present also some
additional structures on it which will be essential to describe minimal $A_\infty$ algebra structures.

\subsubsection{Hochschild Cochain Complex}

Let $A$ be an algebra and let $M$ be a bimodule over $A$.
The Hochschild cochain complex  is defined as
$$
C^*(A,M)=\sum C^n(A,M),\  \ C^n(A,M)=Hom(\otimes^nA,M),
$$
the coboundary operator $\delta: C^n(A,M)\to C^{n+1}(A,M)$ is given by
$$
\begin{array}{c}
\delta f(a_1\otimes ... \otimes a_{n+1})=
a_1\cdot f(a_2\otimes ...
\otimes a_{n+1})\\ +
\sum_{k}\pm f(a_1\otimes ... \otimes a_k\cdot a_{k+1}\otimes  ...
\otimes a_{n+1})\pm f(a_1\otimes ... \otimes a_{n})\cdot a_{n+1}.
\end{array}
$$
The Hochschild cohomology of $A$ with coefficients in $M$ is defined as the homology
of this cochain complex and is denoted by $Hoch^*(A,M)$.

If $A$ and $M$ are graded then each $C^n(A,M)$ is graded too:
$C^n(A,M)=\sum_{k}C^{n,k}(A,M)$ where $C^{n,k}(A,M)=\Hom^k(\otimes^nA,M)$;
here $\Hom^k$ means homomorphisms of degree $k$. So the Hochschild cochain complex is bigraded in this case.

It is clear that $\delta$ maps $ C^{n,k}(A,M)$ to $C^{n+1,k}(A,M)$ , so
$(C^{*,k}(A,M),\delta) $ is a subcomplex in
$(C^*(A,M),\delta)$.
Thus Hochschild cohomology in this case is bigraded:
$Hoch^n(A,M)=\sum_{k} Hoch^{n,k}(A,M)$,
where $ Hoch^{n,k}(A,M)$ is the $n$-th homology module of
$ (C^{*,k}(A,M),\delta )$.

For our needs instead of $M$ we take the algebra $A$ itself. The complex
$ C^{*,*}(A,A)$ is a bigraded differential algebra with respect to the
following cup product:
for $f\in C^{m,k}(A,A)$ and $g\in C^{n,t}(A,A)$ the product
$f\smile g\in C^{m+n,k+t}(A,A)$ is defined by
$$
f\smile g(a_1\otimes ... \otimes a_{m+n})=f(a_1\otimes ...
\otimes a_{m})\cdot g (a_{n+1}\otimes ... \otimes a_{m+n}),
$$
and the Hochschild differential is a derivative:
\begin{equation}
\label{steen0}
\delta (f\smile g)=\delta f\smile g \pm f\smile \delta g.
\end{equation}

\subsubsection{Gerstenhabers circle product}

Besides this product $ C^{*,*}(A,A)$ there exist much richer and important algebraic operations, finally forming a hGa (homotopy Gerstenhaber algebra) structure. First of all, there is Gerstenhaber's so called circle product \cite{Ger63}  $f\circ g$,
sometimes called Gerstenhabers brace  $f\{g\}$, but let us denote it as $f\smile_1g$
since it has properties very similar to Steenrod's $\smile_1$ product in the cochain complex of a topological space.
Here is the definition of this Gerstenhabers product:  for $f\in C^{m,s}(A,A), \ g\in C^{n,t}(A,A)$ their $\smile_1$ product is defined as
\begin{equation}
\label{Cup1}
\begin{array}{c}
f\smile_1g(a_1\otimes ... \otimes a_{m+n-1})=\\
\sum_{k=0}^{m-1}\pm f(a_1\otimes ...\otimes a_k\otimes g(a_{k+1}\otimes ...
\otimes a_{k+n})\otimes ... \otimes a_{m+n-1})\in \\ C^{m+n-1,s+t}(A,A).
\end{array}
\end{equation}
In \cite{Kad88}, see also \cite{Kad07}, it is shown that Gerstenhabers circle = brace = $\smile_1$ product satisfies
Steenrod's condition
\begin{equation}
\label{steen1}
\delta (f\smile_1 g)=\delta f\smile_1 g \pm f\smile_1 \delta g
\pm f\smile g\pm g\smile f.
\end{equation}

What is more amazing,$\smile_1$ also has the following property:
\begin{equation}
\label{Hirsch}
(f\smile g)\smile_1 h= f\smile (g\smile_1 h) \pm (f\smile_1 h)\smile g,
\end{equation}
which means that $...\smile_1 h$ is a derivation. This is an analogue of Hirsch's formula in the cochain complex of a topological space.

Besides,
although the $\smile _1$ is not associative, it is possible to show that it satisfies the so called pre-Jacobi
identity
$$
f\smile _1(g\smile _1h)-(f\smile _1g)\smile _1h=f\smile _1(h\smile _1g)-(f\smile _1h)\smile
_1g
$$
which guarantees that the commutator%
$$
[f,g]=f\smile _1g-g\smile _1f
$$
satisfies the Jacobi identity. Besides, the condition (\ref{steen1}) implies
that $[\ ,\ ]:C^{m,s}(A,A)\otimes C^{n,t}(A,A)\to C^{m+n-1,s+t}(A,A)$
is a chain map, and this implies that $(C^{*,*}(A,A),\delta, [\ ,\ ])$ is a dg Lie algebra.

By the way, the Hochschild differential $\delta$ can be expressed as
$$
\delta f=f\smile_1 \mu+\mu\smile_1 f=[f,\mu]
$$
where $\mu:A\otimes A\to A$ is the produduct operation of $A$.

%%%%%%%%%%%%%%%%%%%%%%%%%%%%%%%%%%%%%%%%%%%%%%%%%%%%%%%%%%%%%%%%%%%%%%%%%%%%%%%%%%%%%%%%%%%%%%%%%%%%%

%%%%%%%%%%%%%%%%%%%%%%%%%%%%%%%%%%%%%%%%%%%%%%%%%%%%%%%%%%%%%%%%%%%%%%%%%%%%%%%%%%%%%%%%%%%%%%%%%%%%%

\subsubsection{Higher operations in the Hochschild complex} \label{GetzlerKad}

In order to perturb $\smile_1$-twisting cohains (i.e., minimal $A_\infty$ algebra structures, as we'll see below) we need a kind of analog of
the Berikishvilis equivalence of Browns twisting cochains, but now for this nonassociative $\smile_1$ case. Unfortunately for this only Gerstenhaber's $\circ =\smile_1$ is not enough, so
in \cite{Kad88} we have introduced higher order multioperations, the $\smile_1$ products $f\smile_1(g_1,...,g_i)$
of one Hochschild cohain $f$ and a collection  of Hochschild cochains $g_1,...,g_i$. Now these operations are called Getzler-Kadeishvili brace operations and are denoted as $f\{g_1,...,g_i\}$. Here is the definition:

\begin{equation}
\label{cup1gg}
\begin{array}{ll}
f\smile_1(g_1,...,g_i)(a_1\otimes ...\otimes a_{n+n_1+...+n_i-i})=&\\

\sum_{k_1,...,k_i}  f(a_1\otimes ...\otimes a_{k_1}\otimes
g_1(a_{k_1+1}\otimes ... \otimes a_{k_1+n_1})
\otimes &\\
\ \ \ ... \otimes a_{k_2}\otimes g_2(a_{k_2+1}\otimes ... \otimes
a_{k_2+n_2})\otimes a_{k_2+n_2+1}
\otimes ...\\
\ \ \ \otimes a_{k_i}\otimes g_i(a_{k_i+1}\otimes ... \otimes
a_{k_i+n_i})\otimes...\otimes a_{n+n_1+...+n_i-i}).
\end{array}
\end{equation}

These higher braces $f\{g_1,...,g_i\}$ form on the Hochschild cochain complex $C^*(A,A)$ a so called homotopy Gerstenhaber algebra (hGa) structure.

\subsubsection{Hochschild cochain complex as a Homotopy Gerstenhaber Algebra}

Let $H$ be a graded algebra. Consider the Hochshild cochain
complex of this graded algebra with coefficients in itself, $C^{*,*}(,H,H)$ which is bigraded in this case:
$$C^{n,m}(H,H)=\Hom^m(H^{\otimes n},H),$$
where $Hom^m$ means homomorphisms of degree $m$.

 This bigraded complex
carries a structure of homotopy Gerstenhaber algebra, see
\cite{Kad88}, \cite{GJ94}, \cite{GV95}, \cite{Kad07}
$C^{*,*}(H,H),\delta, \smile, \{E_{1,k},\ k=1,2,...\})$, which
consists of the following structure maps:

(i) The Hochschild differential $\delta :C^{n-1,m}(H,H)\to
C^{n,m}(H,H)$ given by
$$
\begin{array}{ll}
\delta f(a_1\otimes ...\otimes a_n)=&a_1\cdot f(a_2\otimes ...\otimes a_n) \\
&+\sum_{k}\pm f(a_1\otimes ...\otimes a_{k-1}\otimes a_k\cdot
a_{k+1}\otimes ..\otimes a_n)\\
&\pm f(a_1\otimes ...\otimes a_{n-1})\cdot a_n;
\end{array}
$$
(ii) The $\smile $ product defined by
$$
f\smile g(a_1\otimes ...\otimes a_{n+m})=f(a_1\otimes ...\otimes
a_n)\cdot g(a_{n+1}\otimes ...\otimes a_{n+m}).
$$
(iii) The {\it brace} operations $f\{g_1,...,g_i\}$ from \cite{Kad88}, which we write here
as
$$
f\{g_1,...,g_i\}=E_{1,i}(f;g_1,...,g_i),
$$
$$ E_{1,i}:C^{n,m}\otimes C^{n_1,m_1}\otimes...\otimes
C^{n_i,m_i}\to C^{n+\sum n_t-i,m+\sum m_t},
$$
given by
\begin{equation}
\label{cup1gg}
\begin{array}{ll}
E_{1,i}(f;g_1,...,g_i)(a_1\otimes ...\otimes
a_{n+n_1+...+n_i-i})&\\
= \sum_{k_1,...,k_i} \pm f(a_1\otimes ...\otimes a_{k_1}\otimes
g_1(a_{k_1+1}\otimes ... \otimes a_{k_1+n_1})
\otimes &\\
\ \ \ ... \otimes a_{k_2}\otimes g_2(a_{k_2+1}\otimes ... \otimes
a_{k_2+n_2})\otimes a_{k_2+n_2+1}
\otimes ...\\
\ \ \ \otimes a_{k_i}\otimes g_i(a_{k_i+1}\otimes ... \otimes
a_{k_i+n_i})\otimes...\otimes a_{n+n_1+...+n_i-i}).
\end{array}
\end{equation}
The first brace operation $E_{1,1}$ has the properties of
Steenrod's $\smile_1$ product, so we use the notation
$E_{1,1}(f,g)=f\smile_1 g$. In fact this is Gerstenhaber's $f\circ
g$ product \cite{Ger63}, \cite{Ger64}.

%%%%%%%%%%%%%%%%%%%%%%%%%%%%%%%%%%%%%%%%%%%%%%%%%%%%%%%%%%%%%%%%%%%%%%%%%%%%%%%%%%%%%%%%%%%%%%%%%%%%%%%%%%%%%%%%%

The name {\it Homotopy G-algebra} is motivated by the fact that this
structure induces on homology $H(A)$ the structure of Gerstenhaber algebra
(G-algebra).

The sequence $\{E_{1,k}\}$ defines a twisting cochain
$$
E:BC^{*,*}(H,H)\otimes BC^{*,*}(H,H)\to C^{*,*}(H,H),
$$
and, consequently, defines a strictly associative product
on the bar construction $BC^{*,*}(H,H)$

$$\mu _E:BC^{*,*}(H,H)\otimes BC^{*,*}(H,H)\rightarrow BC^{*,*}(H,H)$$
which turns it into a DG-bialgebra.

%%%%%%%%%%%%%%%%%%%%%%%%%%%%%%%%%%%%%%%%%%%%%%%%%%%%%%%%%%%%%%%%%%%%%%%%%%%%%%%%%%%%%%%%%%%%%%%%%%%%%%%%%%%%%%%%%

%%%%%%%%%%%%%%%%%%%%%%%%%%%%%%%%%%%%%%%%%%%%%%%%%%%%%%%%%%%%%%%%%%%%%%%%%%%%%%%%%%%%%%

\subsection{Description of minimal $A_{\infty}$-algebra structure as a twisting cochain in Hochschild complex} \label{hochschild}

Here we present the connection of the notion of minimal $A_\infty$-algebra with Hochschild
cochain complexes, studied in \cite{Kad88}, see also \cite{Kad07}.

%%%%%%%%%%%%%%%%%%%%%%%%%%%%%%%%%%%%%%%%%%%%%%%%%%%%%%%%%%%%%%%%%%%%%%%%%%%%%%%%%%%%%%%%%%%%%%%%%%%%%%%%%%%%

%%%%%%%%%%%%%%%%%%%%%%%%%%%%%%%%%%%%%%%%%%%%%%%%%%%%%%%%%%%%%%%%%%%%%%%%%%%%%%%%%%%%%%%%%%%%%%%%%%%%%
\subsubsection{Minimal $A_\infty$ structure as a Hochschild twisting cochain}
\label{AinfTwist}

Now let
$$
(H,\{m_1=0, m_2, m_3, \ ...\ ,m_n,\ ...\ \})
$$
be a minimal $A_{\infty}$-algebra. So
$(H,m_2)$ is an associative graded algebra with multiplication
$a\cdot b=m_2(a\otimes b)$ and we can consider the Hochschild cochain complex $C^{*,*}(H,H)$.

Each operation $m_i$ can be considered as a Hochschild cochain
$m_i\in C^{i,2-i}(H,H)$. Let $m=m_3+m_4+...\in C^{*,2-*}(H,H)$.
Stasheff's defining condition of $A_{\infty}$-algebra (\ref{ainf})
\begin{equation}
\label{ainfalg}
\begin{array}{l}
\sum_{k=0}^{n-1} \sum_{j=1}^{n-k}  \\
m_{n-j+1}(a_1\otimes . . . \otimes a_k\otimes  m_j(a_{k+1}\otimes
. . . \otimes a_{k+j})\otimes . . . \otimes a_n)=0 .
\end{array}
\end{equation}
in our case looks like this:
$$
\begin{array}{c}
m_2(a_1\otimes m_{n-1}(a_2\otimes \ ...\ \otimes a_n))+
m_2(m_{n-1}(a_1\otimes \ ...\ \otimes a_{n-1})\otimes a_n)\\ +
\sum_k m_{n-1}(a_1\otimes \ ...\
\otimes a_k\otimes m_2(a_{k+1}\otimes a_{k+2})\otimes a_{k+3}\otimes \ ...\ \otimes a_n) \\ =

\sum_{k=0}^{n-1} \sum_{j=3}^{n-k}  \\
m_{n-j+1}(a_1\otimes . . . \otimes a_k\otimes  m_j(a_{k+1}\otimes
. . . \otimes a_{k+j})\otimes . . . \otimes a_n).
\end{array}
$$
This can be rewritten as
$$
\delta m_n=\sum_{j=3,...,n-2} m_{n-j+1}\smile_1 m_j.
$$
And this means
exactly that $\delta m=m\smile_1 m$ for  $m=m_3+m_4+...,\ \   m_i\in C^{i,2-i}(H,H)$.

So a minimal $A_{\infty}$-algebra
structure on $H$ is, in fact, a {\it twisting cochain}
$$
m=m_3+m_4+...\in C^{*,2-*}(H,H)
$$
in the
Hochschild complex with respect to the (nonassociative) $\smile_1$ product.

%%%%%%%%%%%%%%%%%%%%%%%%%%%%%%%%%%%%%%%%%%%%%%%%%%%%%%%%%%%%%%%%%%%%%%%%%%%%%%%%%%%%%%%%%%%%%%%%%%%%%%%%%%%%

%%%%%%%%%%%%%%%%%%%%%%%%%%%%%%%%%%%%%%%%%%%%%%%%%%%%%%%%%%%%%%%%%%%%%%%%%%%%%%%%%%%%%%%%%%%%%%%%%%%%%%%%%%%%%%%%%

%%%%%%%%%%%%%%%%%%%%%%%%%%%%%%%%%%%%%%%%%%%%%%%%%%%%%%%%%%%%%%%%%%%%%%%%%%%%%%%%%%%%%%%%%%%%%%%%%%%%%%%%%%%%%%%%%

%%%%%%%%%%%%%%%%%%%%%%%%%%%%%%%%%%%%%%%%%%%%%%%%%%%%%%%%%%%%%%%%%%%%%%%%%%%%%%%%%%%%%%%%%%%%%%%%%%%%%%%%%%%%%%%%%

\subsubsection{Perturbations of minimal $A_\infty$ algebras}

Now let $(H,\{m_i\}$ be a minimal $A_{\infty}$-algebra, so
$(H,m_2)$ is an associative graded algebra with multiplication
$a\cdot b=m_2(a\otimes b)$.

As it was mentioned in section \ref{AinfTwist}, each operation $m_i$ can be considered as a Hochschild cochain
$m_i\in C^{i,2-i}(H,H)$. Let $m=m_3+m_4+...\in C^{*,2-*}(H,H)$.
The defining condition of an $A_{\infty}$-algebra (\ref{ainf}) means
exactly $\delta m=m\smile_1 m$. So a minimal $A_{\infty}$-algebra
structure on $H$ is, in fact, a {\it twisting cochain} in the
Hochschild complex with respect to the $\smile_1$ product.

There is the notion of equivalence of such twisting cochains:
$m\sim m'$ if there exists
$p=p^{2,-1}+p^{3,-2}+...+p^{i,1-i}+...,\ ,\ p^{i,1-i}\in
C^{i,1-i}(H,H)$ such that
\begin{equation}
\label{eqv}
\begin{array}{l}

m-m'= \delta p+ p\smile p+p\smile_1m+ \\
 m'\smile_1p+
E_{1,2}(m';p,p)+ E_{1,3}(m';p,p,p)+ ...\ .
\end{array}
\end{equation}
\begin{proposition}. Twisting cochains $m, m'\in C^{*,2-*}(H,H)$ are
equivalent if and only if $(H,\{m_i\})$ and $(H',\{m'_i\})$ are
isomorphic $A_{\infty}$-algebras.
\end{proposition}
\noindent{\bf Proof.} Indeed,
$$
\{p_i\}:(H,\{m_i\})\to (H,\{m'_i\})
$$
with $p_1=id,\ p_i=p^{i,1-i}$ is the needed isomorphism: the
condition (\ref{eqv}) coincides with the defining condition
(\ref{morphism}) of a morphism of $A_\infty$-algebras and the
Proposition \ref{mapiso} implies that this morphism is an
isomorphism.

This gives the possibility of perturbing twisting cochains
without changing their equivalence class:
\begin{proposition}.
Let $m$ be a twisting cochain (i.e., a minimal $A_{\infty}$-algebra
structure on $H$) and let $ p\in C^{n,1-n}(H,H)$ be an arbitrary
cochain. Then there exists a twisting cochain $\bar{m}$,
equivalent to $m$, such that $m_i=\bar{m}_i$ for $i\leq n$ and $
\bar{m}_{n+1}=m_{n+1}+\delta p$.
\end{proposition}
\noindent{\bf Proof.} The twisting cochain $\bar{m}$ can be solved
inductively from the equation (\ref{eqv}).

\begin{theorem}.\label{hoch3} Suppose for a graded algebra $H$ that its Hochschild cohomology groups are trivial,
$Hoch^{n,2-n}(H,H)=0$ for $n\geq 3$. Then each $m\sim 0$, that is
each minimal $A_{\infty}$-algebra structure on $H$ is degenerate, i.e.,
$(H,\{m_i\})$ is isomorphic to a (strictly associative) $A_\infty$ algebra $(H,\{m_1=0,m_2,m_3=0,...\})$.

.
\end{theorem}

\noindent {\bf Proof.} From the equality $\delta m=m\smile_1 m$ in
dimension $4$ we obtain $\delta m_3=0$, that is, $m_3$ is a cocycle.
Since $Hoch^{3,-1}(H,H)=0$ there exists $p^{2,-1}$ such that
$m_3=\delta p^{2,-1}$. Perturbing our twisting cochain $m$ by
$p^{2,-1}$ we we obtain a new twisting cochain
$\bar{m}=\bar{m}_3+\bar{m}_4+...$ equivalent to $m$ and with
$\bar{m}_3=0$. Now the component $\bar{m}_4$ becomes a cocycle,
which can be killed using $Hoch^{4,-2}(H,H)=0$, etc.

%%%%%%%%%%%%%%%%%%%%%%%%%%%%%%%%%%%%%%%%%%%%%%%%%%%%%%%%%%
%%%%%%%%%%%%%%%%%%%%%%%%%%%%%%%%%%%%%%%%%%%%%%%%%%%%%%%%%%
%%%%%%%%%%%%%%%%%%%%%%%%%%%%%%%%%%%%%%%%%%%%%%%%%%%%%%%%%%
%%%%%%%%%%%%%%%%%%%%%%%%%%%%%%%%%%%%%%%%%%%%%%%%%%%%%%%%%%
%%%%%%%%%%%%%%%%%%%%%%%%%%%%%%%%%%%%%%%%%%%%%%%%%%%%%%%%%%

\subsubsection{Minimal $C_{\infty}$-algebra structure and Harrison cohomology}  \label{Harrison}

Suppose now $(H,\mu)$ is a commutative graded algebra. The
Harrison cochain complex $\bar{C}^*(H,H)$ is defined as a
subcomplex of the Hochschild complex consisting of cochains which
disappear on shuffles. If $(H,\{m_i\})$ is a $C_{\infty}$-algebra
then the twisting element $m=m_3+m_4+...$ belongs to Harrison
subcomplex $\bar{C}^*(H,H)\subset C^*(H,H)$ and we have the
\begin{theorem}.\label{harr3} Suppose for a graded commutative algebra $H$
Harrison cohomology $Harr^{n,2-n}(H,H)=0$ for $n\geq 3$. Then each
$m\sim 0$, that is each minimal $C_{\infty}$-algebra structure on
$H$ is degenerate.
\end{theorem}

\section{Minimality theorem, $A_{\infty}$-algebra structure in homology}

%%%%%%%%%%%%%%%%%%%%%%%%%%%%%%%%%%%%%%%%%%%%%%%%%%%%%
%%%%%%%%%%%%%%%%%%%%%%%%%%%%%%%%%%%%%%%%%%%%%%%%%%%%%
%%%%%%%%%%%%%%%%%%%%%%%%%%%%%%%%%%%%%%%%%%%%%%%%%%%%%

\subsection{$A_{\infty}$-algebra structure in homology}

Let $(A,d,\mu )$ be a dg algebra and $(H(A),\mu ^{*})$ be its
homology algebra.

Although the dg algebras $(A,d,\mu )$ and $(H(A),d=0,\mu^* )$ have same homology algebras, the homology algebra $H(A)$ carries less information that the dg algebra $A$.

Generally, there does not exist a map of dg algebras $H(A)\to A$ which induces an isomorphism of homology algebras.
Of course, stepping from $A$ to the smaller object $H(A)$ one looses part of the information. To compensate this loss, it is natural to enrich the algebraic structure on the smaller object $H(A)$. The classical examples of such enrichments are Steenrod squares, Massey products, ... .

\subsubsection{Minimality Theorem}

Below we present one sort of such additional algebraic structure, namely $A_\infty$-algebra structure on the cohomology algebra, the so called minimality theorem.

It was mentioned above (\ref{Ainf}) Stasheff's $A_\infty$-algebras are sort of Strong Homotopy Associative Algebras, the operation $m_3$
is a homotopy which measures the nonassociativity of the product $m_2$.
So its existence on homology the $H(A)$, which is strictly associative
looks a bit strange, but although the product on $H(A)$ is associative,
there appears a structure of a (generally nondegenerate) \emph{minimal}
$A_{\infty}$-algebra, which can be considered as an $A_{\infty}$
\emph{deformation} of $(H(A),\mu^*)$, \cite{Kad07} . Namely, in
\cite{Kad76}, \cite {Kad80} the following \emph{minimalty theorem} was proved:

\begin{theorem}.
\label{Minimality}
Suppose that for a dg algebra $(A,d,\mu)$ all homology modules $H_i(A)$ are
free.

(a) Then there exist a structure of minimal $A_{\infty}$-algebra on $H(A)$
$$
(H(A),\{m_1=0,\ m_2=\mu^*,\ m_3,\ ...\ ,\ m_i\})
$$
and a weak equivalence of
$A_{\infty}$-algebras
$$
\{f_i\}:(H(A),\{m_i\})\rightarrow (A,\{d,\mu ,0,0,...\})
$$
such that $m_1=0$, $m_2=\mu ^{*}$, $f_1^{*}=id_{H(A)}$.

(b) Furthermore, for a dga map $g:A\to A'$ there exists a morphism of
$A_{\infty}$-algebras $\{g_i\}:(H(A),\{m_i\})\to (H(A'),\{m'_i\})$
with $g_1=g^*$ and such that the diagram
$$
\begin{array}{ccc}
(H(A),\{m_i\})& \stackrel{\{g_i\}}{\longrightarrow } & (H(A'),\{m'_i\}) \\

\{f_i\}\downarrow &    &  \downarrow \{f'_i\} \\

A & \stackrel{g}{\longrightarrow } & A'

\end{array}
$$
commutes up to homotopy in the category of $A_\infty$ algebras.

%%%%%%%%%%%%%%%%%%%%%%%%%%%%%%%%%%%%%%%%%%%%%%%%%%%%%%%%%%%%%
%%%%%%%%%%%%%%%%%%%%%%%%%%%%%%%%%%%%%%%%%%%%%%%%%%%%%%%%%%%%%
%%%%%%%%%%%%%%%%%%%%%%%%%%%%%%%%%%%%%%%%%%%%%%%%%%%%%%%%%%%%%
%%%%%%%%%%%%%%%%%%%%%%%%%%%%%%%%%%%%%%%%%%%%%%%%%%%%%%%%%%%%%

(c) Such a
structure is unique up to isomorphism in the category of
$A_{\infty}$-algebras: if $(H(A),\{m_i\})$ and $(H(A),\{m'_i\})$
are two minimal $A_{\infty}$-algebra structures on $H(A)$ which satisfy the conditions from (a)
then there exists an isomorphism of $A_\infty$-algebras $\{g_i\}:(H(A),\{m_i\})\to
(H(A),\{m'_i\})$ with $g_1=id$.

\end{theorem}
%%%%%%%%%%%%%%%%%%%%%%%%%%%%%%%%%%%%%%%%%%%%%%%%%%%%%%%%%%%%%
%%%%%%%%%%%%%%%%%%%%%%%%%%%%%%%%%%%%%%%%%%%%%%%%%%%%%%%%%%%%%
%%%%%%%%%%%%%%%%%%%%%%%%%%%%%%%%%%%%%%%%%%%%%%%%%%%%%%%%%%%%%
%%%%%%%%%%%%%%%%%%%%%%%%%%%%%%%%%%%%%%%%%%%%%%%%%%%%%%%%%%%%%

\noindent {\bf Proof of (a).} We are going to construct the components $f_i,\ m_i$ inductively satisfying the
defining condition of an $A_\infty$-morphism
\begin{equation}
\label{Amor}
\begin{array}{c}

df_n(a_1\otimes...\otimes a_{n})=f_1m_n(a_1\otimes...a_n)\\ +
\sum_{j=2}^{n-1}\sum_{k=0}^{n-j}f_{n-j+1}(a_1\otimes ...\otimes
a_k\otimes m_j(a_{k+1}\otimes ...\otimes a_{k+j})\otimes ...\otimes a_n)\\ +
\sum_{k=1}^{n-1}f_{k}(a_1\otimes...\otimes a_{k})\cdot f_{n-k}(a_{k+1}\otimes...\otimes a_{n}).
\end{array}
\end{equation}

Let us start with a cycle-choosing
homomorphism $f_1:H(A)\to A$, that is $f_1(a)\in a\in H(A)$ which can be constructed using the freeness of the modules $H_i(A)$. This map is not
multiplicative but $ f_1(a_1\cdot a_2)-f_1(a_1)\cdot f(a_2)\sim 0\in A $.
So there exists $f_2:H(A)\otimes H(A)\to A$ s.t.
$f_1(a_1\cdot
a_2)-f_1(a_1)\cdot f(a_2)=d f_2(a_1\otimes _2)$, and this, assuming $m_2(a_1\otimes a_2)=a_1\cdot a_2$,
is exactly the condition (\ref{Amor}) for $n=2$.

Let us denote
$$
\begin{array}{l}
U_3(a_1\otimes a_2\otimes a_3)= \\
f_1(a_1)\cdot f_2(a_2\otimes a_3)+ f_2(a_1\cdot a_2\otimes a_3)+
f_2(a_1\otimes a_2\cdot a_3)+ f_2(a_1\otimes a_2)\cdot f_1(a_3).
\end{array}
$$

Note that the main condition (\ref{Amor}) for $n=3$ looks as follows
$$
df_3(a_1\otimes a_2\otimes a_3)=f_1m_3(a_1\otimes a_2\otimes a_3)-U_3(a_1\otimes a_2\otimes a_3).
$$

Direct calculation shows that $dU_3(a_1\otimes a_2\otimes a_3)=0$, so $U_3(a_1\otimes a_2\otimes a_3)$ is a cycle.
Let us define
$m_3(a\otimes b\otimes c)\in H(A)$ as the homology class of this cycle
\begin{equation}
\label{m3}
\begin{array}{c}
m_3(a\otimes b\otimes c)=\{U_3(a\otimes b\otimes c)\}=\\
f_1(a_1)\cdot f_2(a_2\otimes a_3)+ f_2(a_1\cdot a_2\otimes a_3)+
f_2(a_1\otimes a_2\cdot a_3)+ f_2(a_1\otimes a_2)\cdot f_1(a_3).
\end{array}
\end{equation}
Then, since $f_1$ is a \emph{cycle choosing} homomorphism,
$$
f_1m_3(a_1\otimes a_2\otimes a_3)-U_3(a_1\otimes a_2\otimes a_3)
$$
is homological to zero. Thus, again using the freeness of $H_i(A)$, it is possible to construct
a homomorphism $f_3:H(A)\otimes H(A)\otimes H(A) \to A$ such that
$$
df_3(a_1\otimes a_2\otimes a_3)=f_1m_3(a_1\otimes a_2\otimes a_3)-U_3(a_1\otimes a_2\otimes a_3),
$$
and this is exactly the condition (\ref{Amor}) for $n=3$.

Assume now that $f_i,\ m_i, \ i\leq n-1$ are already constructed and they satisfy (\ref{Amor}).

Let us denote
$$
\begin{array}{c}
U_n(a_1\otimes...\otimes a_{n})\\=
\sum_{j=2}^{n-1}\sum_{k=0}^{n-j}f_{n-j+1}(a_1\otimes ...\otimes
a_k\otimes m_j(a_{k+1}\otimes ...\otimes a_{k+j})\otimes ...\otimes a_n)\\ +
\sum_{k=1}^{n-1}f_{k}(a_1\otimes...\otimes a_{k})\cdot f_{n-k}(a_{k+1}\otimes...\otimes a_{n}).
\end{array}
$$
Then the defining condition (\ref{Amor}) can be rewritten as
\begin{equation}\label{Un}
df_n(a_1\otimes...\otimes a_{n})=f_1m_n(a_1\otimes...a_n)+U_n(a_1\otimes...\otimes a_n).
\end{equation}
Direct calculation shows that $dU_n=0$, so we define
$m_n(a_1\otimes...\otimes a_n)$ as the homology class of the cycle $U_n(a_1\otimes...\otimes a_n)$. Then using the freeness of modules the $H_i(A)$ one can construct a homomorphism $f_n:H_*(A)\otimes...\otimes H(A)\to A$ satisfying (\ref{Un}), and this completes the proof of (a).

%%%%%%%%%%%%%%%%%%%%%%%%%%%%%%%%%%%%%%%%%%%%%%%%%%%%%%%%%%%%%
%%%%%%%%%%%%%%%%%%%%%%%%%%%%%%%%%%%%%%%%%%%%%%%%%%%%%%%%%%%%%
%%%%%%%%%%%%%%%%%%%%%%%%%%%%%%%%%%%%%%%%%%%%%%%%%%%%%%%%%%%%%

%%%%%%%%%%%%%%%%%%%%%%%%%%%%%%%%%%%%%%%%%%%%%%%%%%%%%%%%%%%%%
%%%%%%%%%%%%%%%%%%%%%%%%%%%%%%%%%%%%%%%%%%%%%%%%%%%%%%%%%%%%%
%%%%%%%%%%%%%%%%%%%%%%%%%%%%%%%%%%%%%%%%%%%%%%%%%%%%%%%%%%%%%

\noindent {\bf Proof of (b).}
The morphism of $A_\infty$-algebras
$(H(A),\{m_i\}) \stackrel{\{f_i\}}{\longrightarrow }A \stackrel{g}{\longrightarrow A'}$ induces a twisting cochain
$\phi:B(H(A),\{m_i\}) \to A$.

Furthermore, since
$$
(H(A'),\{m'_i\}) \stackrel{\{f'_i\}}{\longrightarrow }A'
$$
is a weak equivalence of $A_\infty$-algebras, by the lifting theorem (\ref{lifta}) $\phi$ can be
lifted to the $A_\infty$-twisting cochain $\psi:B(H(A),\{m_i\})\to (H(A'),\{m'_i\})$ which in fact represents an $A_\infty$-morphism
$\{g_i\}:(H(A),\{m_i\})\to (H(A'),\{m'_i\})$,
such that
$\{f'_i\}\circ \psi \sim \phi$. This equivalence means that
$A_\infty$-morphisms $\{f'_i\}\circ \{g_i\}$ and
$g\circ \{g_i\}$ are homotopic. This completes the proof.

%%%%%%%%%%%%%%%%%%%%%%%%%%%%%%%%%%%%%%%%%%%%%%%%%%%%%%%%%%%%%
%%%%%%%%%%%%%%%%%%%%%%%%%%%%%%%%%%%%%%%%%%%%%%%%%%%%%%%%%%%%%
%%%%%%%%%%%%%%%%%%%%%%%%%%%%%%%%%%%%%%%%%%%%%%%%%%%%%%%%%%%%%

\noindent {\bf Proof of (c).} Suppose for a dg algebra $(A,d,\mu)$ we have two homology
$A_{\infty}$-algebras $(H(A),\{m_i\})$ and $(H(A),\{m'_i\})$, i.e., there exist $A_\infty$ weak equivalences
$\{f_i\}:(H(A),\{m_i\})\rightarrow A, \ \ \{f'_i\}:(H(A),\{m_i\})\rightarrow A$ with $(f_1)^{*}=(f'_1)^{*}=id_{H(A)}$.
Then, using the part (b) for $d=id_A$ we obtain $\{g_i\}:(H(A),\{m_i\})\to (H(A),\{m'_i\})$ with $g_1=g^*=id_{H(A)}$. So, $\{g_i\}$ is an isomorphism of $A_\infty$ algebras by (\ref{mapweak}).

%%%%%%%%%%%%%%%%%%%%%%%%%%%%%%%%%%%%%%%%%%%%%%%%%%%%%%%%%%%%%
%%%%%%%%%%%%%%%%%%%%%%%%%%%%%%%%%%%%%%%%%%%%%%%%%%%%%%%%%%%%%
%%%%%%%%%%%%%%%%%%%%%%%%%%%%%%%%%%%%%%%%%%%%%%%%%%%%%%%%%%%%%
%%%%%%%%%%%%%%%%%%%%%%%%%%%%%%%%%%%%%%%%%%%%%%%%%%%%%%%%%%%%%

\begin{corollary}.
\label{Cor1} The mapping of differential coalgebras
$$
B(\{f_i\}): \tilde{B}(H(A ),\{m_i\})\to B(A )
$$
induces an isomorphism in homology.
\end{corollary}

\subsubsection{Connection with Massey Products.} The first new operation $m_3(a_1\otimes a_2\otimes a_3)\in H(A)$ was defined
as the homology class of the
cycle
$$
\begin{array}{l}
U_3(a_1\otimes a_2\otimes a_3)= \\
f_1(a_1)\cdot f_2(a_2\otimes a_3)+ f_2(a_1\cdot a_2\otimes a_3)+
f_2(a_1\otimes a_2\cdot a_3)+ f_2(a_1\otimes a_2)\cdot f_1(a_3).
\end{array}
$$

This description immediately induces the connection of $m_3$
with the Massey product: If $a,b,c\in H(A)$  is a  Massey triple, i.e.,
if $a_1\cdot a_2=a_2\cdot a_3=0$, then
$$
U_3(a_1\otimes a_2\otimes a_3)=f_1(a_1)\cdot f_2(a_2\otimes a_3)+ f_2(a_1\otimes b_2)\cdot f_1(a_3)
$$
and this is exactly the combination which defines the Massey product $\large<a_1,a_2,a_3\large>$, thus
$m_3(a_1\otimes a_2\otimes a_3)\in <a_1,a_2,a_3>$.
This gives examples of dg
algebras with essentially nontrivial homology $A_\infty$-algebras.

\subsubsection{Special cases}

Taking $(A,d,\mu)=C^*(X)$, the cochain complex of a topological space $X$, the theorem defines
on the cohomology algebra $H^*(X)$ the structure of a
minimal $A_\infty$-algebra $(H^*(X),\{m_i\})$,  which carries
more information about $X$ than just the graded algebra structure.
In particular, we shall show later that this cohomology $A_\infty$-algebra $(H^*(X),\{m_i\})$ determines the cohomology modules of the loop space $\Omega X$ whereas the bare the cohomology algebra $(H^*(X),m_2)$ does not.
s

Furthermore, taking $(A,d,\mu)=C_*(G)$, the chain complex of a topological group, the theorem defines on the
Pontriagin ring $H_*(G)$ the structure
of a minimal $A_\infty$-algebra $(H_*(G),\{m_i\})$  which will carry
more information about $G$ than just the ring structure.
In particular, we shall show later that this homology $A_\infty$-algebra $(H_*(G),\{m_i\})$ determines
homology modules of classifying space $B_G$ whereas the bare Pontriagin ring $(H_*(G),m_2)$ does not.

\subsubsection{Minimal $A_\infty$-algebra structure on the homology of an $A_\infty$-algebra}

The Minimality Theorem is true also when, instead of a dg
algebra $(A,d,\mu)$ we take
an arbitrary $A(\infty)$-algebra $(M,\{m_i\})$.

In this case, minimal $A_\infty$-algebra structure appears on the homology $H(M)$ of the dg module $(M,m_1:M\to M)$, see \cite{Kad82}.

This structure had applications in string theory, see for example \cite{Aspin}.

%%%%%%%%%%%%%%%%%%%%%%%%%%%%%%%%%%%%%%%%%%%%%%%%%%%%%%%%%%%%%%%%%%%%%%%%%%%%%%%%%%%%%%%%%%%%%%%%%%%
%%%%%%%%%%%%%%%%%%%%%%%%%%%%%%%%%%%%%%%%%%%%%%%%%%%%%%%%%%%%%%%%%%%%%%%%%%%%%%%%%%%%%%%%%%%%%%%%%%%
%%%%%%%%%%%%%%%%%%%%%%%%%%%%%%%%%%%%%%%%%%%%%%%%%%%%%%%%%%%%%%%%%%%%%%%%%%%%%%%%%%%%%%%%%%%%%%%%%%%
%%%%%%%%%%%%%%%%%%%%%%%%%%%%%%%%%%%%%%%%%%%%%%%%%%%%%%%%%%%%%%%%%%%%%%%%%%%%%%%%%%%%%%%%%%%%%%%%%%%
%%%%%%%%%%%%%%%%%%%%%%%%%%%%%%%%%%%%%%%%%%%%%%%%%%%%%%%%%%%%%%%%%%%%%%%%%%%%%%%%%%%%%%%%%%%%%%%%%%%
%%%%%%%%%%%%%%%%%%%%%%%%%%%%%%%%%%%%%%%%%%%%%%%%%%%%%%%%%%%%%%%%%%%%%%%%%%%%%%%%%%%%%%%%%%%%%%%%%%%
%%%%%%%%%%%%%%%%%%%%%%%%%%%%%%%%%%%%%%%%%%%%%%%%%%%%%%%%%%%%%%%%%%%%%%%%%%%%%%%%%%%%%%%%%%%%%%%%%%%
%%%%%%%%%%%%%%%%%%%%%%%%%%%%%%%%%%%%%%%%%%%%%%%%%%%%%%%%%%%%%%%%%%%%%%%%%%%%%%%%%%%%%%%%%%%%%%%%%%%
%%%%%%%%%%%%%%%%%%%%%%%%%%%%%%%%%%%%%%%%%%%%%%%%%%%%%%%%%%%%%%%%%%%%%%%%%%%%%%%%%%%%%%%%%%%%%%%%%%%
%%%%%%%%%%%%%%%%%%%%%%%%%%%%%%%%%%%%%%%%%%%%%%%%%%%%%%%%%%%%%%%%%%%%%%%%%%%%%%%%%%%%%%%%%%%%%%%%%%%
%%%%%%%%%%%%%%%%%%%%%%%%%%%%%%%%%%%%%%%%%%%%%%%%%%%%%%%%%%%%%%%%%%%%%%%%%%%%%%%%%%%%%%%%%%%%%%%%%%%
%%%%%%%%%%%%%%%%%%%%%%%%%%%%%%%%%%%%%%%%%%%%%%%%%%%%%%%%%%%%%%%%%%%%%%%%%%%%%%%%%%%%%%%%%%%%%%%%%%%
%%%%%%%%%%%%%%%%%%%%%%%%%%%%%%%%%%%%%%%%%%%%%%%%%%%%%%%%%%%%%%%%%%%%%%%%%%%%%%%%%%%%%%%%%%%%%%%%%%%
%%%%%%%%%%%%%%%%%%%%%%%%%%%%%%%%%%%%%%%%%%%%%%%%%%%%%%%%%%%%%%%%%%%%%%%%%%%%%%%%%%%%%%%%%%%%%%%%%%%
%%%%%%%%%%%%%%%%%%%%%%%%%%%%%%%%%%%%%%%%%%%%%%%%%%%%%%%%%%%%%%%%%%%%%%%%%%%%%%%%%%%%%%%%%%%%%%%%%%%

 \section{Applications of the  Minimality Theorem}

\subsection{Application: Cohomology $A_\infty$ algebra of a space and Cohomology modules of Loop Space.}
Taking $A=C^*(X)$, the cochain dg algebra of a 1-connected space
$X$, we obtain an $A_{\infty}$-algebra structure
$(H^*(X),\{m_i\})$ on the cohomology algebra $H^*(X)$.

The cohomology algebra equipped with this additional structure carries
more information than just the cohomology algebra. Some
applications of this structure are given in \cite{Kad80} ,
\cite{Kad93}. For example the cohomology $A_{\infty}$-algebra
$(H^*(X),\{m_i\})$ determines the cohomology of the loop space
$H^*(\Omega X)$ whereas the bare algebra $(H^*(X),m_2)$ does not:

\begin{theorem}. $H(B(H^*(X),\{m_i\}))=H^*(\Omega X)$.
\end{theorem}

%%%%%%%%%%%%%%%%%%%%%%%%%%%%%%%%%%%%%%%%%%%%%%%%%%%%%%%%%%%%%%%%%%%%%%%%%%%%%%%%%%%%%%%%%%%%%%%%%%%
%%%%%%%%%%%%%%%%%%%%%%%%%%%%%%%%%%%%%%%%%%%%%%%%%%%%%%%%%%%%%%%%%%%%%%%%%%%%%%%%%%%%%%%%%%%%%%%%%%%
%%%%%%%%%%%%%%%%%%%%%%%%%%%%%%%%%%%%%%%%%%%%%%%%%%%%%%%%%%%%%%%%%%%%%%%%%%%%%%%%%%%%%%%%%%%%%%%%%%%
%%%%%%%%%%%%%%%%%%%%%%%%%%%%%%%%%%%%%%%%%%%%%%%%%%%%%%%%%%%%%%%%%%%%%%%%%%%%%%%%%%%%%%%%%%%%%%%%%%%
%%%%%%%%%%%%%%%%%%%%%%%%%%%%%%%%%%%%%%%%%%%%%%%%%%%%%%%%%%%%%%%%%%%%%%%%%%%%%%%%%%%%%%%%%%%%%%%%%%%
%%%%%%%%%%%%%%%%%%%%%%%%%%%%%%%%%%%%%%%%%%%%%%%%%%%%%%%%%%%%%%%%%%%%%%%%%%%%%%%%%%%%%%%%%%%%%%%%%%%
%%%%%%%%%%%%%%%%%%%%%%%%%%%%%%%%%%%%%%%%%%%%%%%%%%%%%%%%%%%%%%%%%%%%%%%%%%%%%%%%%%%%%%%%%%%%%%%%%%%
%%%%%%%%%%%%%%%%%%%%%%%%%%%%%%%%%%%%%%%%%%%%%%%%%%%%%%%%%%%%%%%%%%%%%%%%%%%%%%%%%%%%%%%%%%%%%%%%%%%
%%%%%%%%%%%%%%%%%%%%%%%%%%%%%%%%%%%%%%%%%%%%%%%%%%%%%%%%%%%%%%%%%%%%%%%%%%%%%%%%%%%%%%%%%%%%%%%%%%%
%%%%%%%%%%%%%%%%%%%%%%%%%%%%%%%%%%%%%%%%%%%%%%%%%%%%%%%%%%%%%%%%%%%%%%%%%%%%%%%%%%%%%%%%%%%%%%%%%%%
%%%%%%%%%%%%%%%%%%%%%%%%%%%%%%%%%%%%%%%%%%%%%%%%%%%%%%%%%%%%%%%%%%%%%%%%%%%%%%%%%%%%%%%%%%%%%%%%%%%
%%%%%%%%%%%%%%%%%%%%%%%%%%%%%%%%%%%%%%%%%%%%%%%%%%%%%%%%%%%%%%%%%%%%%%%%%%%%%%%%%%%%%%%%%%%%%%%%%%%
%%%%%%%%%%%%%%%%%%%%%%%%%%%%%%%%%%%%%%%%%%%%%%%%%%%%%%%%%%%%%%%%%%%%%%%%%%%%%%%%%%%%%%%%%%%%%%%%%%%
%%%%%%%%%%%%%%%%%%%%%%%%%%%%%%%%%%%%%%%%%%%%%%%%%%%%%%%%%%%%%%%%%%%%%%%%%%%%%%%%%%%%%%%%%%%%%%%%%%%
%%%%%%%%%%%%%%%%%%%%%%%%%%%%%%%%%%%%%%%%%%%%%%%%%%%%%%%%%%%%%%%%%%%%%%%%%%%%%%%%%%%%%%%%%%%%%%%%%%%

Firstly, according to the minimality theorem the homology of the dg algebra $(C^*(X),d,\smile )$ carries
a structure of minimal $A_\infty$-algebra $(H^*(X),\{m_i\})$, such that there is  and a weak equivalence
of $A_\infty$-algebras
$$
\{f_i\}:(H^*(X),\{m_i\})\to (C^*(X),\{m_1=d,\ m_2=\smile,\ m_3=0,\ ... \ \},
$$
which induces weak equivalences of their bar constructions and an isomorphism of graded modules
$$
H(B(H^*(X),\{m_i\})) \approx H(BC^*(X))\approx H^*(\Omega X).
$$
Thus the object $(H^*(X),\{m_i\})$, which is called the cohomology $A_\infty$-algebra of $X$, determines the cohomology \emph{modules} of the loop space $H^*(\Omega X)$. But not the cohomology \emph{algebra} $H^*(\Omega X)$.

%%%%%%%%%%%%%%%%%%%%%%%%%%%%%%%%%%%%%%%%%%%%%%%%%%%%%%%%%%%%%%%%%%%%%%%%%%%%%%%%%%%%%%%%%%%%%%%%%%%
%%%%%%%%%%%%%%%%%%%%%%%%%%%%%%%%%%%%%%%%%%%%%%%%%%%%%%%%%%%%%%%%%%%%%%%%%%%%%%%%%%%%%%%%%%%%%%%%%%%
%%%%%%%%%%%%%%%%%%%%%%%%%%%%%%%%%%%%%%%%%%%%%%%%%%%%%%%%%%%%%%%%%%%%%%%%%%%%%%%%%%%%%%%%%%%%%%%%%%%
%%%%%%%%%%%%%%%%%%%%%%%%%%%%%%%%%%%%%%%%%%%%%%%%%%%%%%%%%%%%%%%%%%%%%%%%%%%%%%%%%%%%%%%%%%%%%%%%%%%
%%%%%%%%%%%%%%%%%%%%%%%%%%%%%%%%%%%%%%%%%%%%%%%%%%%%%%%%%%%%%%%%%%%%%%%%%%%%%%%%%%%%%%%%%%%%%%%%%%%
%%%%%%%%%%%%%%%%%%%%%%%%%%%%%%%%%%%%%%%%%%%%%%%%%%%%%%%%%%%%%%%%%%%%%%%%%%%%%%%%%%%%%%%%%%%%%%%%%%%
%%%%%%%%%%%%%%%%%%%%%%%%%%%%%%%%%%%%%%%%%%%%%%%%%%%%%%%%%%%%%%%%%%%%%%%%%%%%%%%%%%%%%%%%%%%%%%%%%%%
%%%%%%%%%%%%%%%%%%%%%%%%%%%%%%%%%%%%%%%%%%%%%%%%%%%%%%%%%%%%%%%%%%%%%%%%%%%%%%%%%%%%%%%%%%%%%%%%%%%
%%%%%%%%%%%%%%%%%%%%%%%%%%%%%%%%%%%%%%%%%%%%%%%%%%%%%%%%%%%%%%%%%%%%%%%%%%%%%%%%%%%%%%%%%%%%%%%%%%%
%%%%%%%%%%%%%%%%%%%%%%%%%%%%%%%%%%%%%%%%%%%%%%%%%%%%%%%%%%%%%%%%%%%%%%%%%%%%%%%%%%%%%%%%%%%%%%%%%%%
%%%%%%%%%%%%%%%%%%%%%%%%%%%%%%%%%%%%%%%%%%%%%%%%%%%%%%%%%%%%%%%%%%%%%%%%%%%%%%%%%%%%%%%%%%%%%%%%%%%
%%%%%%%%%%%%%%%%%%%%%%%%%%%%%%%%%%%%%%%%%%%%%%%%%%%%%%%%%%%%%%%%%%%%%%%%%%%%%%%%%%%%%%%%%%%%%%%%%%%
%%%%%%%%%%%%%%%%%%%%%%%%%%%%%%%%%%%%%%%%%%%%%%%%%%%%%%%%%%%%%%%%%%%%%%%%%%%%%%%%%%%%%%%%%%%%%%%%%%%
%%%%%%%%%%%%%%%%%%%%%%%%%%%%%%%%%%%%%%%%%%%%%%%%%%%%%%%%%%%%%%%%%%%%%%%%%%%%%%%%%%%%%%%%%%%%%%%%%%%
%%%%%%%%%%%%%%%%%%%%%%%%%%%%%%%%%%%%%%%%%%%%%%%%%%%%%%%%%%%%%%%%%%%%%%%%%%%%%%%%%%%%%%%%%%%%%%%%%%%

\subsection{Application: Homology modules of the Classifying Space of a Topological Group.}

Taking $A=C_*(G)$, the chain dg algebra of a topological group
$G$, we obtain an $A_{\infty}$-algebra structure
$(H_*(G),\{m_i\})$ on the Pontriagin algebra $H_*(G)$.
 The homology
$A_{\infty}$-algebra $(H_*(G),\{m_i\})$ determines the homology of the
classifying space $H_*(B_G)$ when just the Pontryagin algebra
$(H_*(G),m_2)$ does not:

\begin{theorem}. $H(B(H_*(G),\{m_i\}))=H_*(B_G)$.
\end{theorem}

%%%%%%%%%%%%%%%%%%%%%%%%%%%%%%%%%%%%%%%%%%%%%%%%%%%%%%%%%%%%%%%%%%%%%%%%%%%%%%%%%%%%%%%%%%%%%%%%%%%
%%%%%%%%%%%%%%%%%%%%%%%%%%%%%%%%%%%%%%%%%%%%%%%%%%%%%%%%%%%%%%%%%%%%%%%%%%%%%%%%%%%%%%%%%%%%%%%%%%%
%%%%%%%%%%%%%%%%%%%%%%%%%%%%%%%%%%%%%%%%%%%%%%%%%%%%%%%%%%%%%%%%%%%%%%%%%%%%%%%%%%%%%%%%%%%%%%%%%%%
%%%%%%%%%%%%%%%%%%%%%%%%%%%%%%%%%%%%%%%%%%%%%%%%%%%%%%%%%%%%%%%%%%%%%%%%%%%%%%%%%%%%%%%%%%%%%%%%%%%
%%%%%%%%%%%%%%%%%%%%%%%%%%%%%%%%%%%%%%%%%%%%%%%%%%%%%%%%%%%%%%%%%%%%%%%%%%%%%%%%%%%%%%%%%%%%%%%%%%%

\subsection{Application: $A_\infty$-model of a Fibre Bundle}

%%%%%%%%%%%%%%%%%%%%%%%%%%%%%%%%%%%%%%%%%%%%%%%%%%%%%%%%%%%%%%%%%%%%%%%%%%%%%%%%%%%%%%%%%%%%%%%%%%%
%%%%%%%%%%%%%%%%%%%%%%%%%%%%%%%%%%%%%%%%%%%%%%%%%%%%%%%%%%%%%%%%%%%%%%%%%%%%%%%%%%%%%%%%%%%%%%%%%%%
%%%%%%%%%%%%%%%%%%%%%%%%%%%%%%%%%%%%%%%%%%%%%%%%%%%%%%%%%%%%%%%%%%%%%%%%%%%%%%%%%%%%%%%%%%%%%%%%%%%
%%%%%%%%%%%%%%%%%%%%%%%%%%%%%%%%%%%%%%%%%%%%%%%%%%%%%%%%%%%%%%%%%%%%%%%%%%%%%%%%%%%%%%%%%%%%%%%%%%%
%%%%%%%%%%%%%%%%%%%%%%%%%%%%%%%%%%%%%%%%%%%%%%%%%%%%%%%%%%%%%%%%%%%%%%%%%%%%%%%%%%%%%%%%%%%%%%%%%%%

The minimality theorem (\ref{Minimality}) and the theorem (\ref{lifta}) about lifting of twisting cochains allow
to construct effective models of fibre bundles. Actually this model and  higher operations $\{m_i\}$ and $\{p_i\}$ were constructed in \cite{Kad76}. Later we have recognized that they form Stasheff's $A_\infty$-structures and the model in these terms was presented in \cite{Kad80}. A similar model was also presented in \cite{Smi80}.

\noindent{\bf Topological level.} Let $\xi=(X,p,B,G)$ be a principal $G$-fibration. Let $F$ be a $G$-space. Then the action $G\times F\to F$ determines the \emph{associated fibre bundle} $\xi(F)=(E,p,B,F,G)$ with fiber $F$. Thus $\xi$ and the action $G\times F\to F$ on the topological level determine $E$.

\noindent{\bf Chain level.}
Let $K=C_*(B),\ A=C_*(G),\ P=C_*(F)$. The classical result of E. Brown \cite{Brown} states that the principal fibration $\xi$ determines a twisting cochain $\phi:K=C_*(B)\to A=C_*(G)$, and the action on the chain level $C_*(G)\otimes C_*(F)\to C_*(F)$ defines the twisted tensor product
$K\otimes_\phi P=C_*(B)\otimes_\phi C_*(F)$ which gives the homology modules of the total space $H_*(E)$.
Thus $\xi$ and the action on the chain level $C_*(G)\otimes C_*(F)\to C_*(F)$ determine $H_*(E)$.

The twisting cochain $\phi$ is not uniquely determined and it can be perturbed by the above equivalence relations
for computational reasons.

\noindent{\bf Homology level.} Nodar Berikashvili stated the problem to lift the previous "chain level" model of associated fibration to "homology level", i.e., to construct a "twisted differential" on $C_*(B)\otimes H_*(F)$.
Investigation had shown that the principal fibration $\xi$ and the action of of Pontryagin ring $H_*(G)$ on $H_*(F)$, that is the pairing
$H_*(G)\otimes H_*(F)\to H_*(F)$ \emph{do not determine} $H_*(E)$. But by the minimality theorem it appeared that $H_*(G)$ carries not only the Pontryagin product
$H_*(G)\otimes H_*(G)\to H_*(G)$ but also richer algebraic structure, namely the structure of a minimal $A_\infty$-algebra $(H_*(G),\{m_i\})$.

Furthermore, the action $G\times F\to F$ induces not only the pairing  $H_*(G)\otimes H_*(F)\to H_*(F)$
but, by the modular analog of the Minimality Theorem \cite{Kad80}
 also the structure of minimal $A_\infty$-module $(H_*(F),\{p_i\})$,
 $$
 p_i:H_*(G)\otimes ... ((i-1)\ times)...\otimes H_*(G)\otimes H_*(F)\to H_*(F),
 $$
and all these operations allow to define the correct differential on $C_*(B)\otimes H_*(F)$:
there is a weak equivalence, homology isomorphism
$$
C_*(B)\otimes_\psi H_*(F) = K\otimes_\psi H(P)\to  K\otimes_\phi P = C_*(B)\otimes_\phi C_*(F)\sim C_*(E).
$$

%%%%%%%%%%%%%%%%%%%%%%%%%%%%%%%%%%%%%%%%%%%%%%%%%%%%%%%%%%%%%%%%%%%%%%%%%%%%%%%%%%%%%%%%%%%%%%%%%%%
%%%%%%%%%%%%%%%%%%%%%%%%%%%%%%%%%%%%%%%%%%%%%%%%%%%%%%%%%%%%%%%%%%%%%%%%%%%%%%%%%%%%%%%%%%%%%%%%%%%
%%%%%%%%%%%%%%%%%%%%%%%%%%%%%%%%%%%%%%%%%%%%%%%%%%%%%%%%%%%%%%%%%%%%%%%%%%%%%%%%%%%%%%%%%%%%%%%%%%%
%%%%%%%%%%%%%%%%%%%%%%%%%%%%%%%%%%%%%%%%%%%%%%%%%%%%%%%%%%%%%%%%%%%%%%%%%%%%%%%%%%%%%%%%%%%%%%%%%%%
%%%%%%%%%%%%%%%%%%%%%%%%%%%%%%%%%%%%%%%%%%%%%%%%%%%%%%%%%%%%%%%%%%%%%%%%%%%%%%%%%%%%%%%%%%%%%%%%%%%

%%%%%%%%%%%%%%%%%%%%%%%%%%%%%%%%%%%%%%%%%%%%%%%%%%%%%%%%%%
%%%%%%%%%%%%%%%%%%%%%%%%%%%%%%%%%%%%%%%%%%%%%%%%%%%%%%%%%%
%%%%%%%%%%%%%%%%%%%%%%%%%%%%%%%%%%%%%%%%%%%%%%%%%%%%%%%%%%
%%%%%%%%%%%%%%%%%%%%%%%%%%%%%%%%%%%%%%%%%%%%%%%%%%%%%%%%%%
%%%%%%%%%%%%%%%%%%%%%%%%%%%%%%%%%%%%%%%%%%%%%%%%%%%%%%%%%%

  \section{$C_{\infty}$-algebra structure in the homology
  of a commutative dg algebra, applications in Rational Homotopy Theory}

There is a commutative version of the above Minimality Theorem,
see \cite{Kad88h}, \cite{Kad93}, \cite{Markl}:

\begin{theorem}.\label{cinf}
Suppose that for a commutative dg algebra $A$ all homology $R$-modules
$H^i(A)$ are free.

Then there exist: a structure of minimal $C_{\infty}$-algebra
$(H(A),\{m_i\})$ on $H(A)$ and a weak equivalence of
$C_{\infty}$-algebras
$$
\{f_i\}:(H(A),\{m_i\})\rightarrow (A,\{d,\mu ,0,0,...\})
$$
such that $m_1=0$, $m_2=\mu ^{*}$, $f_1^{*}=id_{H(A)}$.

Furthermore, for a cdga map $f:A\to A'$ there exists a morphism of
$C_{\infty}$-algebras $\{f_i\}:(H(A)\{m_i\})\to (H(A')\{m'_i\})$
with $f_1=f^*$.
\end{theorem}
 Such a
structure is unique up to isomorphism in the category of
$C_{\infty}$-algebras.

Below we present some applications of this $C_\infty$-algebra
structure in rational homotopy theory.

%%%%%%%%%%%%%%%%%%%%%%%%%%%%%%%%%%%%%%%%%%%%%%%%%%%%%%%%%%%%%%
%%%%%%%%%%%%%%%%%%%%%%%%%%%%%%%%%%%%%%%%%%%%%%%%%%%%%%%%%%%%%%
%%%%%%%%%%%%%%%%%%%%%%%%%%%%%%%%%%%%%%%%%%%%%%%%%%%%%%%%%%%%%%
%%%%%%%%%%%%%%%%%%%%%%%%%%%%%%%%%%%%%%%%%%%%%%%%%%%%%%%%%%%%%%
%%%%%%%%%%%%%%%%%%%%%%%%%%%%%%%%%%%%%%%%%%%%%%%%%%%%%%%%%%%%%%

%%%%%%%%%%%%%%%%%%%%%%%%%%%%%%%%%%%%%%%%%%%%%%%%%%%%%%%%%%%%%%
%%%%%%%%%%%%%%%%%%%%%%%%%%%%%%%%%%%%%%%%%%%%%%%%%%%%%%%%%%%%%%
%%%%%%%%%%%%%%%%%%%%%%%%%%%%%%%%%%%%%%%%%%%%%%%%%%%%%%%%%%%%%%
%%%%%%%%%%%%%%%%%%%%%%%%%%%%%%%%%%%%%%%%%%%%%%%%%%%%%%%%%%%%%%
%%%%%%%%%%%%%%%%%%%%%%%%%%%%%%%%%%%%%%%%%%%%%%%%%%%%%%%%%%%%%%
%%%%%%%%%%%%%%%%%%%%%%%%%%%%%%%%%%%%%%%%%%%%%%%%%%%%%%%%%%%%%%

%\section{Applications in Rational Homotopy Theory}

\subsection{Classification of rational homotopy types}

Let $X$ be a 1-connected space. In the case of rational
coefficients there exist Sullivan's {\it commutative} cochain
complex $A(X)$ of $X$. It is well known that the weak equivalence
type of a cdg algebra $A(X)$ determines the rational homotopy type
of $X$: 1-connected $X$ and $Y$ are {\it rationally homotopy
equivalent} if and only if $A(X)$ and $A(Y)$ are weakly
equivalent cdg algebras. Indeed, in this case $A(X)$ and $A(Y)$
have {\it isomorphic} minimal models $M_X\approx M_Y$, and this
implies that $X$ and $Y$ are rationally homotopy equivalent. This
is the key geometrical result of Sullivan which we are going to
exploit below.

Now we take $A=A(X)$ and apply Theorem \ref{cinf}. Then we
obtain on $H(A)=H^*(X,Q)$ a structure of minimal $C_{\infty}$
algebra $(H^*(X,Q),\{m_i\})$ which we call the {\it rational
cohomology $C_{\infty}$-algebra of $X$}.

Generally, an isomorphism of rational cohomology algebras $H^*(X,Q)$
and $ H^*(Y,Q)$ does not imply a homotopy equivalence $X\sim Y$, not even
rationally. We claim that $(H^*(X,Q),\{m_i\})$ is a {\it complete}
rational homotopy invariant:

\begin{theorem}. 1-connected $X$ and $X'$ are rationally
homotopy equivalent if and only if
 $(H^*(X,Q),\{m_i\})$ and $(H^*(X',Q),\{m'_i\})$ are
isomorphic as $C_{\infty}$-algebras.
\end{theorem}
\noindent{\bf Proof.} Suppose $X\sim X'$. Then $A(X)$ and $A(X')$
are weakly equivalent, that is there exists a cgda $A$ and weak
equivalences $A(X)\leftarrow A\to A(X')$. This implies weak
equivalences of the corresponding homology $C_{\infty}$-algebras
$$
(H^*(X,Q),\{m_i\})\leftarrow (H^*(A),\{m_i\})\to
(H^*(X',Q),\{m'_i\}),
$$
which by of minimality are both isomorphisms.

Conversely, suppose $(H^*(X,Q),\{m_i\})\approx
(H^*(X',Q),\{m'_i\})$. Then

$${\cal A}
QB(H^*(X,Q),\{m_i\})\approx {\cal A} QB(H^*(X',Q),\{m'_i\}).$$
Denote this cdga as $A$. Then we have weak equivalences of CGD
algebras  (see (\ref{adj}))
$$
A(X)\leftarrow {\cal A} \Gamma A(X)\leftarrow A\to {\cal A} \Gamma
A(X') \to A(X').
$$

This theorem in fact classifies rational homotopy types with given
cohomology algebra $H$ as all possible minimal
$C_{\infty}$-algebra structures on $H$ modulo $C_{\infty}$
isomorphisms.

\vspace{3mm}

\noindent{\bf Example.} Here we describe an example which we will
use to illustrate the results of this and forthcoming sections.

We consider the following commutative graded algebra. Its
underlying graded $Q$-vector space has the following generators: a generator
$e$ of dimension 0, generators $x,\ y$ of dimension 2, and a
generator $z$ of dimension 5, so
\begin{equation}\label{example}
H^*=\{H^0=Q_e,\ 0,\ H^2=Q_x\oplus Q_y,\ 0,\ 0,\ H^5=Q_z,\ 0,\ 0,\
...\ \},
\end{equation}
and the multiplication is trivial by dimensional reasons, with
unit $e$. In fact
$$
H^*=H^*(S^2\vee S^2\vee S^5,Q).
$$
This example was considered in \cite{HS}. It was shown there that
there are just two rational homotopy types with such cohomology
algebra.

The same result can be obtained from our classification.

What minimal $C_\infty$-algebra structures are possible on $H^*$?

For dimensional reasons only one nontrivial operation $
m_3:H^2\otimes H^2\otimes H^2\to H^5 $ is possible.

The specific condition of a $C_\infty$-algebra, namely the
disappearance on shuffles implies that
$$
m_3(x,x,x)=0,\ m_3(y,y,y)=0,\ m_3(x,y,x)=0,\ m_3(y,x,y)=0
$$
and
$$
m_3(x,x,y)=m_3(y,x,x),\ \ m_3(x,y,y)=m_3(y,y,x).
$$
Thus each $C_\infty$-algebra structure on $H^*$ is characterized
by a pair of rational numbers $p,\ q$ and
$$
m_3(x,x,y)=pz,\ \ m_3(x,y,y)=qz.
$$
So let us write an arbitrary minimal $C_\infty$-algebra structure
on $H^*$ as a column vector
$\left(\begin{array}{c}p\\q\end{array}\right)$.

Now let us look at the structure of an isomorphism of
$C_\infty$-algebras
$$
\{f_i\}:(H^*,m_3)\to (H^*,m'_3).
$$
Again for dimensional reasons just one component $f_1:H^*\to H^* $
is possible, which in its turn consists of two isomorphisms
$$
f_1^2:H^2=Q_x\oplus Q_y\to H^2=Q_x\oplus Q_y,\ \ f_1^5:H^5=Q_z\to
H^5=Q_z.
$$
The first one is represented by a nondegenerate matrix
$A=\left(\begin{array}{cc}a&c\\b&d\end{array}\right)$,
$$
f^2_1(x)=ax\oplus by,\ \ f^2_1(y)=cx\oplus dy,
$$
and the second one by a nonzero rational number $r$,
$f^5_1(z)=rz$.

A calculation shows that the condition $f^5_1m_3=m'_3(f^2_1\otimes
f^2_1\otimes f^2_1)$, to which the defining condition
of an $A_\infty$-algebra morphism (\ref{morphism})degenerates looks as followes
$$
r\left(\begin{array}{c}p\\q\end{array}\right)=det\
A\left(\begin{array}{cc}a&b\\c&d\end{array}\right)
\left(\begin{array}{c}p'\\q'\end{array}\right).
$$
This condition shows that two minimal $C_\infty$-algebra
structures $m_3=\left(\begin{array}{c}p\\q\end{array}\right)$ and
$m'_3=\left(\begin{array}{c}p'\\q'\end{array}\right)$ are
isomorphic if and only if they are related by a nondegenerate linear
transformation.

Thus there exist just two isomorphism classes of minimal
$C_\infty$-algebras on $H^*$: the trivial one $(H^*,m_3=0)$ and the
nontrivial one $(H^*,m_3\neq 0)$. So we have just two rational
homotopy types whose rational cohomology is $H^*$. We denote them by
$X$ and $Y$, respectively and analyze them in the next sections.

\vspace{5mm}

 Below we give some applications of the cohomology $C_\infty$-algebra in various
 problems of rational homotopy theory.

\subsection{Formality}

%%%%%%%%%%%%%%%%%%%%%%%%%%%%%%%%%%%%%%%%%%%%%%%%%%%%%%%%%%%%%%%%%
Among rational homotopy types with given cohomology algebra, there
is one called {\it formal} which is a "formal consequence of its
cohomology algebra" (Sullivan). Explicitly, this is the type whose
minimal model $M_X$ is isomorphic to the minimal model of the
cohomology $H^*(X,Q)$.

Our $C_{\infty}$ model implies the following criterion of
formality:

\begin{theorem}.\label{formal} X is formal if and only if its cohomology
$C_{\infty}$-algebra is degenerate, i.e., it is $C_{\infty}$
isomorphic to one with $m_{\geq 3}=0$.
\end{theorem}

Below we deduce some known results about
formality using this criterion.

\vspace{3mm}

\noindent 1. A commutative graded 1-connected algebra $H$ is
called {\it intrinsically formal} if there is only one homotopy type
with this cohomology algebra $H$, of course the formal one.

The above Theorem \ref{harr3} immediately implies the following
sufficient condition for formality due to Tanre \cite{Tan85}:

\begin{theorem}.  If for a 1-connected graded $Q$-algebra $H$ one has
$$Harr^{k,k-2}(H,H)=0,\ k=3,4,...$$ then
$H$ is intrinsically formal, that is there exists only one
rational homotopy type with $H^*(X,Q)\approx H$.
\end{theorem}

\noindent 2. The following theorem of Halperin and Stasheff from
\cite{HS} is an immediate result of our criterion:
\begin{theorem}.\label{halsta} A commutative graded $Q$-algebra of
type
$$
H=\{H^0=Q,0,0,...,0,H^n,H^{n+1},...,H^{3n-2},0,0,...\}
$$
is intrinsically formal
\end{theorem}
\noindent{\bf Proof.} Since $deg\ m_i=2-i$ there is no room for
operations $m_{i>2}$, indeed the shortest range is $m_3:H^n\otimes
H^n\otimes H^n\to H^{3n-1}=0. $

\vspace{5mm}

\noindent 3. Theorem \ref{formal} easily implies the
\begin{theorem}. Any  1-connected commutative graded algebra $H$ with
$H^{2k}=0$ for all $k$ is intrinsically formal.
\end{theorem}
\noindent {\bf Proof.} Any $A_{\infty}$-operation $m_i$ has degree
$2-i$, thus
$$
m_i:H^{2k_1+1}\otimes ...\otimes H^{2k_i+1}\to
H^{2(k_1+...+k_i+1)}=0.
$$
Thus any $C_{\infty}$-operation is trivial too.

This implies a result of Baues: any space whose even-dimensional
cohomologies are trivial has the rational homotopy type
of a wedge of spheres. Indeed, such an algebra is realized by wedge
of spheres and because of intrinsical formality this is the only
homotopy type.

\vspace{3mm}

%%%%%%%%%%%%%%%%%%%%%%%%%%%%%%%%%%%%%
%%%%%%%%%%%%%%%%%%%%%%%%%%%%%%%%%%%%%
%%%%%%%%%%%%%%%%%%%%%%%%%%%%%%%%%%%%%

\noindent{\bf Example.} The algebra $H^*$ from the example of the
previous section is not intrinsically formal since there are two
homotopy types, $X$ and $Y$, with $H^*(X,Q)=H^*=H^*(Y)$. The space
$X$  is formal (and actually $X=S^2\vee S^2\vee S^5$), since its
cohomology $C_\infty$-algebra $(H^*,m_3=0)$ is trivial. But the
space $Y$ is not: its cohomology $C_\infty$-algebra $(H^*,m_3\neq
0)$ is not degenerate.

We remark here that the formal type is represented by $X=S^2\vee
S^2\vee S^5$ and it is possible to show that the nonformal one is
represented by $Y=S^2\vee S^2\cup_{f:S^4\to S^2\vee S^2} e^5$,
where the attaching map $f$ is a nontrivial element from
$\pi_4(S^2\vee S^2)\otimes Q$.

%%%%%%%%%%%%%%%%%%%%%%%%%%%%%%%%%%%%%%%%%%%%%%%%%%%%%%%%%%%%%%%%%

%%%%%%%%%%%%%%%%%%%%%%%%%%%%%%%%%%%%%%%%%%%%%%%%%%%%%%%%%
\subsection{Rational homotopy groups}

Since the cohomology $C_{\infty}$-algebra $(H^*(X,Q),\{m_i\})$
determines the rational homotopy type it must determine the
rational homotopy groups $\pi_i(X)\otimes Q$ too. We present a
chain complex whose homology is $\pi_i(X)\otimes Q$. Moreover, the
Lie algebra structure is determined as well.

For a cohomology $C_{\infty}$-algebra $\ \ (H^*(X,Q),\{m_i\})$ the bar
construction $\ \ B(H^*(X,Q),\{m_i\})$ is a dg bialgebra, see (\ref{BarCinfty}). Acting on this
bialgebra by the functor $Q$ of indecomposables we obtain a dg Lie
coalgebra.

On the other hand the rational homotopy groups $\pi_*(\Omega X)\otimes
Q$ form a graded Lie algebra with respect to Whiethead product.
Thus its dual cohomotopy groups $\pi^*(\Omega X,Q)=(\pi_*(\Omega
X)\otimes Q)^*$ form a graded Lie coalgebra.

\begin{theorem}. Homology of a dg Lie coalgebra $QB(H^*(X,Q),\{m_i\})$ is isomorphic to the
cohomotopy Lie coalgebra $\pi^*(\Omega X,Q)$.
\end{theorem}
\noindent{\bf Proof.} The theorem follows from the sequence of
graded Lie coalgebra isomorphisms:
$$
\begin{array}{ll}
\pi^*(\Omega X,Q)\approx  &(\pi_*(\Omega X,Q))^*\approx
 (PH_*(\Omega
x,Q)^* \approx QH^*(\Omega X,Q) \\ \approx
&QH(B(A(X)) \approx QH(\tilde{B}(H^*(X,Q),\{m_i\}) \\ \approx
 &H(Q\tilde{B}(H^*(X,Q),\{m_i\}).
\end{array}
$$

\noindent{\bf Example.} For the algebra $H^*$ from the previous
examples the complex $QB(H^*)$ in low dimensions looks as
$$
0\to Q_x\oplus Q_y\stackrel{0}{\to}Q_{x\otimes x}\oplus
Q_{x\otimes y}\oplus Q_{y\otimes y}\stackrel{0}{\to} Q_{x\otimes
x\otimes y}\oplus Q_{x\otimes y\otimes y}\stackrel{d=m_3}{\to}
Q_z\oplus ... \ .
$$
The differential $d=m_3$ is trivial for the formal space $X$ and
is nontrivial for $Y$. Thus for both rational homotopy types we
have
$$
\pi^2=H^1(QB(H^*))=2Q,\ \ \pi^3=H^2(QB(H^*))=3Q,
$$
and
$$
\begin{array}{l}
\pi^4(X)=H^3(QB(H^*),d=0)=2Q,\\
\pi^4(Y)=H^3(QB(H^*),d\neq 0)=Ker\
d=Q.
\end{array}
$$

%%%%%%%%%%%%%%%%%%%%%%%%%%%%%%%%%%%%%%%%%%%%%%%%%%%%%%%%%

\subsection{Realization of homomorphisms}

Let $G:H^*(X,Q)\to H^*(Y,Q)$ be a homomorphism of cohomology
algebras. When this homomorphism is realizable as a map of
rationalizations $g:Y_Q\to X_Q ,\ \ g^*=G$? In the case when $G$
is an isomorphism this question was considered in \cite{HS}. It
was considered also in \cite{Vig81}. The following theorem gives
the complete answer:

\begin{theorem}.\label{realizat} A homomorphism $G$ is realizable if and only if it is extendable to a
$C_{\infty}$-map
$$
\{g_1=G,g_2,g_3,...\}:(H^*(X,Q),\{m_i\})\to
(H^*(Y,Q),\{m'_i\}).
$$
\end{theorem}
\noindent{\bf Proof.} One direction of this is consequence of the last part
of Theorem \ref{cinf}.

To show the other direction we use Sullivan's minimal models $M_X$ and
$M_Y$ of $A(X)$ and $A(Y)$. It is enough to show that the
existence of $\{g_i\}$ implies the existence of cdg algebra map
$g:M_Y\to M_X$.

So we have $C_{\infty}$-algebra maps
$$
M_X\stackrel{\{f_i\}}{\leftarrow}
(H^*(X,Q),\{m_i\})\stackrel{\{g_i\}}{\to}
(H^*(y,Q),\{m'_i\})\stackrel{\{f'_i\}}{\to}M_Y.
$$
Recall the following property of a minimal cdg algebra $M$: for a
weak equivalence of cdg algebras $\phi:A\to B$ and a cdg algebra
map $f:M\to B$ there exists a cdg algebra map $F:M\to A$ such that
$\phi F$ is homotopic to $f$. Using this property it is easy to
show the existence of a cdga map $\beta:M_X\to {\cal A}QB(M_X)$,
the right inverse of the standard map $\alpha:{\cal A}QB(M_X)\to
M$. Composing this map with ${\cal A}QB({\{f'_i\}}){\cal
A}QB({\{g_i\}})$ we obtain a cdga map
$$
{\cal A}QB({\{f'_i\}}){\cal A}QB({\{g_i\}})\beta:M_X\to M_Y.
$$
This theorem immediately implies the
{\begin{corollary} For
formal $X$ and $Y$ each $G:H^*(X,Q)\to H^*(Y,Q)$ is realizable.
\end{corollary}
\noindent{\bf Proof.} In this case $\{G,0,0,...\}$ is a
$C_{\infty}$-extension f $G$.

\noindent{\bf Example.} Consider the homomorphism
$$
G:H^*(X)=H^*(Y)\to H^*(S^5)
$$
induced by the standard imbedding $g:S^5\to X=S^2\vee S^2\vee
S^5$. Of course $G$ is realizable as $g:S^5\to X$ but not as
$S^5\to Y$. Indeed, for such realizability, according to Theorem
\ref{realizat}, we need a $C_\infty$-algebra morphism
$$
\{g_i\}:(H^*,\{0,0,m_3,0,...\}\to (H^5(S^5,Q),\{0,0,0,...\})
$$
with $g_1=G$. For dimensional reasons all the components
$g_2,g_3,...$ all are trivial, so this morphism has the form
$\{G,0,0,...\}$. But this collection is not a morphism of
$C_\infty$-algebras since the condition $Gm_3=0$, to which
 the defining condition (\ref{morphism}) of an
$A_\infty$-algebra morphism degenerates, is not satisfied.

A. Razmadze Mathematical Institute,

1, M. Alexidze Str., Tbilisi, 0193, Georgia

tornike.kadeishvili@tsu.ge

\end{document}